\documentclass{article}
\usepackage{amssymb,amsmath,bm,geometry,graphics,graphicx,url,color}
\usepackage{amsfonts}
\usepackage{graphicx}
\usepackage{multirow}

\def\pdt2{\partial_t^2}
\def\pdx2{\partial_x^2}

\newcommand{\normmm}[1]{{\left\vert\kern-0.25ex\left\vert\kern-0.25ex\left\vert #1
    \right\vert\kern-0.25ex\right\vert\kern-0.25ex\right\vert}}

\def\RR{{\mathbb{R}}}
\def\CC{{\mathbb{C}}}

\def\pdt2{\partial_t^2}
\def\pdx2{\partial_x^2}

\newcommand{\abs}[1]{\left\vert#1\right\vert}

\def\RR{{\mathbb{R}}}
\def\CC{{\mathbb{C}}}

\def\e{{\mathrm e}}

\def\eps{\varepsilon}

\DeclareMathOperator{\sinc}{sinc}

\newtheorem{theo}{Theorem}[section]

\newtheorem{rem}[theo]{Remark}

\newtheorem{algorithm}[theo]{Algorithm}

\usepackage[all]{xy}

\def\no{\noindent}
\title{Structure-preserving algorithms with uniform error bound and long-time energy conservation for
highly oscillatory Hamiltonian systems}

\author{Bin Wang\,\footnote{School of Mathematics and Statistics,
Xi'an Jiaotong University, Xi'an, Shannxi   710049, P.R.China.
E-mail:~{\tt wangbinmaths@xjtu.edu.cn}} \and Yaolin Jiang
\thanks{School of Mathematics and Statistics,
Xi'an Jiaotong University, Xi'an, Shannxi   710049, P.R.China. E-mail:~{\tt
yljiang@mail.xjtu.edu.cn}}}

\begin{document}
\maketitle

\begin{abstract}
Structure-preserving algorithms and algorithms with uniform error bound
{have constituted} two interesting  classes of numerical
methods. In this paper, we {blend} these two kinds of
methods for solving {nonlinear} Hamiltonian systems with
highly oscillatory solution, and the {blended} algorithms
{inherit and respect} the advantage of each method.
 Two kinds of algorithms are presented to preserve the symplecticity and energy of the Hamiltonian systems, {respectively.}
 Moreover, the proposed algorithms are shown to have uniform error bound for the highly oscillatory structure.   A numerical experiment
 is  carried out to support the theoretical results {established in this paper} by showing the performance of the {blended} algorithms.

\medskip

\no{Keywords:} Hamiltonian system, highly
oscillatory solution, symplectic algorithms, energy-preserving algorithms, uniform error bound, long-time conservation

\medskip
\no{MSC:} 65L05, 65L20, 65L70, 65P10.

\end{abstract}
\section{Introduction}\label{intro}
It is known that  nonlinear Hamiltonian systems are ubiquitous  in
science and engineering applications. In  numerical simulation of
evolutionary problems,  one of the most difficult problems   is to
deal with highly oscillatory problems, {since they cannot
be solved efficiently using conventional methods. The crucial point
is that standard methods need a very small stepsize and hence a long
runtime to reach an acceptable accuracy \cite{hairer2006}. In this paper we are
concerned with efficient algorithms for} the following damped
second-order differential equation
 \begin{equation}\label{charged-particle sts-cons}
\begin{array}[c]{ll}
\ddot{x}(t)=\frac{1}{\varepsilon}\tilde{B} \dot{x}(t)   +F(x(t)), \quad
x(0)=x_0,\quad \dot{x}(0)=\dot{x}_0,\quad t\in[0,T],
\end{array}
\end{equation}
where  $x(t)\in \RR^d$, $\tilde{B}$ is a $d\times d$ skew symmetric
matrix, and  $F(x) = -\nabla_x U(x)$ is the negative gradient of a
real-valued function $U(x)$ whose second derivatives are continuous.
In this work, we focus on the {case where}
$0<\varepsilon\ll 1$. {This implies}
   that the solution of this dynamic   is {\emph{highly oscillatory}}.
For the dimension $d$, it is required that $d\geq2$ {since
$\tilde{B} $ is a zero matrix once $d=1$, and then} the system
\eqref{charged-particle sts-cons} {reduces to} a
second-order ODE $\ddot{x}(t)=F(x(t))$ without highly oscillatory
{solutions.} Denote by $v=\dot{x}$ and then the energy of
this dynamic  is given by
\begin{equation}\label{energy of cha}
E(x,v)=\frac{1}{2}\abs{v}^2+U(x),
\end{equation}
which is exactly conserved along the solutions, i.e.
$$E(x(t),v(t))=E(x(0),v(0))\  \textmd{for any}\ t\in[0,T].$$
We note further that
with  $ p = v
-\frac{1}{2\varepsilon}\tilde{B}x, $ the equation
\eqref{charged-particle sts-cons} can be transformed into a
Hamiltonian system with the non-separable Hamiltonian
 \begin{equation}\label{H}H(x,p)=\frac{1}{2} \abs{p+\frac{1}{2\varepsilon}\tilde{B}x}^2 + U(x).\end{equation}

Hamiltonian systems with highly oscillatory {solutions}
frequently occur in physics and engineering such as charged-particle
dynamics, Vlasov equations, classical and quantum mechanics, and molecular dynamics.
 Their numerical computation   contains numerous enduring challenges.  In the recent few decades, \emph{geometric numerical integration}  also called as structure-preserving algorithm for
differential equations has received more and more attention. This kind of algorithms is designed to respect the structural invariants and geometry
of the considered system. This idea has been used  by many researchers to derive  different  structure-preserving algorithms (see, e.g.
\cite{FENG,LUBICH,hairer2006,wubook2018}).   For the Hamiltonian system \eqref{H},  there are two remarkable features: the symplecticity of its flow and the conservation of the Hamiltonian.  Consequently, for a numerical algorithm, these two
features should be respected as much as possible in the spirit of
geometric numerical integration.

One typical example of system \eqref{charged-particle sts-cons} is Vlasov equations or
charged-particle dynamics in a  strong and uniform magnetic field,
which has been studied by many researchers. {In
practice,} the numerical methods used to treat this system can be
summarized in the following three categories.

a) The primitive numerical methods usually depend on the knowledge of
certain other characteristics of the solution besides high-frequency
oscillation and structure preservation such as the Boris method \cite{Boris} as well as its further researches \cite{Hairer2017-1,Qin2013}.
This method does not perform well for highly oscillatory systems and cannot preserve any structure of the system.

b) Some recent methods are devoted to the structure preservation
such as the volume-preserving algorithms \cite{He2015},  symplectic
methods \cite{He2017,Tao2016}, symmetric methods \cite{Hairer2017-2}
and energy-preserving methods \cite{2020JCAM,Li-ANM,Li-AML,2020JCP}.
In \cite{Hairer2018}, the {long-time} near-conservation
property of a variational integrator was analyzed under
{the condition} $0<\eps\ll1$.
Very recently, some integrators with  large stepsize and their long term behaviour were studied in \cite{LUBICH21} for charged-particle dynamics.
 All of these methods can
preserve or nearly preserve some structure of the considered system.
However, these methods mentioned above do not pay attention to the
high-frequency oscillation, {and then}
  the convergence of these methods is not uniformly accurate for $\varepsilon$.
  Their error constant usually increases when $\eps$ decreases.

c)  Accuracy is often an important consideration  for  highly  oscillatory systems over
long-time intervals. Some new methods with uniform accuracy for $\varepsilon$ have been proposed and analysed recently. The authors in \cite{lubich19} improved asymptotic behaviour of the Boris method and derived
a filtered Boris algorithm   under a maximal ordering scaling. Some   multiscale schemes have been proposed such as the asymptotic preserving schemes \cite{VP4,VP5} and the uniformly accurate schemes \cite{VP1,Zhao}. Although these powerful numerical methods have very good performance in accuracy,  structure (nearly) preservation usually cannot be achieved.

Based on the above points, {a natural question to ask is}
whether one can design a numerical method for
\eqref{charged-particle sts-cons} such that it has uniform error bound
for $\varepsilon$ and can exactly preserve some structure
simultaneously.  A kind of energy-preserving method without
convergent analysis was given in \cite{Wang2020}. It will be shown
in this paper that this method has a uniform error bound which has not
been studied in \cite{Wang2020}. Very recently,    the authors in
\cite{zhao2020} presented some splitting methods with first-order uniform error bound in $x$ and
energy or volume preservation. However, only first-order methods are
proposed there and  {higher-order} ones with energy or
other structure preservation  have not been
{investigated.} A numerical method combining high-order
uniform error bound and structure preservation has more challenges and
importance.

In this paper, we will derive two kinds of algorithms to preserve the symplecticity and energy, respectively. For symplectic algorithms, their near energy conservation over long times will be analysed. Moreover, all the structure-preserving algorithms will be shown to have second-order uniform error bound for $0<\varepsilon\ll 1$ in $x$. Meanwhile, an algorithm with first-order uniform error bound in both $x$ and $v$ will be proposed.
 The remainder of this paper is organised as follows. In Section \ref{sec: methods}, we formulate two kinds of algorithms.
 The main results of these algorithms are given in Section \ref{sec: main results} and a numerical experiment is carried out there to numerically show the performance of the algorithms. The proofs of the main results are presented in Sections
  \ref{sec: symp methods}-\ref{sec: UA peoperty} one by one. The last section includes some  concluding remarks.

\section{Numerical algorithms}\label{sec: methods}
Before  deriving effective algorithms for the system
\eqref{charged-particle sts-cons}, we first present {the
implicit expression of} its exact solution  as follows.

\begin{theo} \label{VOC-SPE} (See \cite{lubich19}.) The  exact solution of
system \eqref{charged-particle sts-cons}  can be expressed as
\begin{equation}\label{VOC}
\begin{array}[c]{ll}&x(t_n+h)=x(t_n)+ h\varphi_1(h \Omega)
v(t_n)+h^2 \int_{0}^1(1-\tau)
\varphi_1((1-\tau)h \Omega) F(x(t_n+h\tau))  d\tau,\\
&v(t_n+h)=\varphi_0(h\Omega)v(t_n)+h \int_{0}^1
\varphi_0((1-\tau)h \Omega) F(x(t_n+h\tau))  d\tau,
 \end{array}
\end{equation}
where $\Omega=\frac{1}{\varepsilon} \tilde{B}$, $h$ is a stepsize, $t_n=nh$ and the
$\varphi$-functions are defined by (see
\cite{Hochbruck2010})
\begin{equation*}
 \varphi_0(z)=e^{z},\ \ \varphi_k(z)=\int_{0}^1
e^{(1-\sigma)z}\frac{\sigma^{k-1}}{(k-1)!}d\sigma, \ \ k=1,2,\ldots.
\end{equation*}
\end{theo}

In what follows, we present two kinds of algorithms which  will correspond to symplectic algorithms and energy-preserving algorithms, respectively.

\begin{algorithm}
\label{scheme 2}
By denoting
the numerical solution $x_n\approx x(t_n),\, v_n\approx v(t_n)$ with $n=0,1,\ldots,$
  an $s$-stage {adaptive} exponential algorithm
{applied with stepsize $h$ is defined by:}
\begin{equation}\label{AEI-new}
\begin{array}[c]{ll}X_{i}=x_{n}+ c_ih\varphi_1(c_ih\Omega) v_{n}+h^2 \sum\limits_{j=1}^{s}{\alpha}_{ij}(h\Omega)F (X_j),\ \ i=1,2,\ldots,s,\\
x_{n+1}=x_{n}+ h\varphi_1(h\Omega) v_{n}+h^2
\sum\limits_{i=1}^{s}\beta_i(h\Omega)F (X_i),\\
v_{n+1}=\varphi_0(h\Omega)v_{n}+h
\sum\limits_{i=1}^{s}\gamma_i(h\Omega)F
(X_i),
\end{array}
\end{equation}
where ${\alpha}_{ij}(h\Omega), \beta_i(h\Omega), \gamma_i(h\Omega)$ are bounded functions of $h\Omega$. As some practical examples, we present five explicit algorithms.

For constant  $F\equiv F_0\in\RR^3$, the variation-of-constants formula \eqref{VOC} reads \begin{equation*}
\begin{array}[c]{ll}&x(t_n+h)=x(t_n)+ h\varphi_1(h \Omega)
v(t_n)+h^2
\varphi_2(h \Omega) F_0,\\
&v(t_n+h)=\varphi_0(h\Omega)v(t_n)+h
\varphi_1(h \Omega) F_0.
 \end{array}
\end{equation*}
Based on this, we consider the following algorithm
\begin{equation*}
\begin{aligned}
&x_{n+1}=x_{n}+ h\varphi_1(h\Omega) v_{n}+h^2
\varphi_2(h\Omega)F (x_{n}),\\
&v_{n+1}=\varphi_0(h\Omega)v_{n}+h
\varphi_1(h\Omega)F (x_{n}).
\end{aligned}
\end{equation*}
 which means that in \eqref{AEI-new}
 $$s=1,\ \ c_{1}=0,\ \   \alpha_{11}=0,\ \   \beta_{1}=\varphi_2(h\Omega),\ \  \gamma_{1}=\varphi_1(h\Omega).$$
 This method is referred to  M1. This method can be verified to be non-symmetric. We modify it to be a symmetric method by considering
 \begin{equation*}
\begin{array}[c]{llll}
s=2,\ \ c_1=0,\ \ c_2=1,& \alpha_{11}=\alpha_{12}=\alpha_{22}=0,\\
\beta_{1}=\varphi_2(h \Omega) ,\ \  \beta_{2}=0, & \gamma_{1}=\frac{\varphi_2(h \Omega)}{\varphi_1(-h \Omega)},
 \ \gamma_{2}=\frac{2e^{h \Omega}\varphi_2(-h
\Omega)}{\varphi_1(h \Omega)},
\end{array}
\end{equation*}
 and denote this method as M2.

 It is noted that the next three methods are formulated based on  the conditions \eqref{17} of symplecticity given blew.
The coefficients are obtained by considering the $s$-stage {adaptive} exponential algorithm \eqref{AEI-new}
with the coefficients for $i=1,\ldots,s,\  j=1,\ldots,i,$
\begin{equation}\label{aeicoe}
\begin{array}[c]{ll}
\alpha_{ij}=a_{ij}(c_{i}-c_{j})\varphi_1((c_{i}-c_{j})h\Omega),\
\beta_{i}=b_{i}(1-c_{i})\varphi_1((1-c_{i})h\Omega),\
 \gamma_i=b_{i}e^{(1-c_{i})h\Omega},
\end{array}
\end{equation}
where $c=(c_1,\ldots,c_{s}),b=(b_{1}$,\ \
$\ldots,b_{s})$ and $A=(a_{ij})_{s\times s}$ are coefficients of an
$s$-stage diagonal implicit RK method. It can be checked easily that if   this RK method is chosen as a symplectic method, then the corresponding coefficients \eqref{aeicoe} satisfy
the conditions \eqref{17}. We omit the details of calculations for brevity.
We first consider
$$s=1,\ \ c_{1}=\frac{1}{2},\ \   b_{1}=1.$$
The {adaptive} exponential algorithm  whose coefficients are given by this choice and  \eqref{aeicoe} is denoted by SM1.
 For $s=2$, choosing
\begin{equation*}
\begin{array}[c]{llll}
c_1=0,\ c_2=1,\ a_{21}=\frac{1}{2},\
\beta_{1}=\frac{1}{2},\  b_{2}=1
\end{array}
\end{equation*}
 yields another method, which is called as SM2.
If we consider
\begin{equation*}
\begin{array}[c]{llll}
c_{1}=\frac{1}{4},\ c_{2}=\frac{3}{4},\ a_{21}=\frac{1}{2},\ b_{1}=b_{2}=\frac{1}{2},
\end{array}
\end{equation*}
 {then the corresponding} method is referred to SM3.
\end{algorithm}

The following algorithm is devoted to the energy-preserving methods which are designed based on the variation-of-constants formula \eqref{VOC} and the idea of continuous-stage methods.
\begin{algorithm}
An $s$-degree continuous-stage {adaptive} exponential
algorithm {applied with stepsize $h$} is defined
{by}
\begin{equation}\label{CSAEI}
\begin{aligned} &X_{\tau}=x_{n}+ hC_{\tau}(h\Omega) v_{n}+h^2 \int_{0}^{1}{A}_{\tau\sigma}(h\Omega)F (X_\sigma)d\sigma,\ \ \ 0\leq\tau\leq1,\\
&x_{n+1}=x_{n}+ h\varphi_1(h\Omega) v_{n}+h^2
\int_{0}^{1}\bar{B}_\tau(h\Omega)F (X_\tau)d\tau,\\
&v_{n+1}=\varphi_0(h\Omega)v_{n}+h
\int_{0}^{1}B_{\tau}(h\Omega)F
(X_\tau)d\tau,
\end{aligned}
\end{equation}
where $ X_{\tau} $ is a polynomial of degree $ s $ with respect to $
\tau $ satisfying 
 $X_{0}=x_{n},\
X_{1}=x_{n+1}.$
\
  $C_{\tau} $, $ \bar{B}_{\tau},
B_{\tau} $ and $A_{\tau,\sigma}$ are polynomials which depend on $
h\Omega $. The $ C_{\tau}(h\Omega) $ satisfies $
C_{c_{i}}(h\Omega)=c_{i}\varphi_{1}(c_{i}h\Omega),$ where $ c_{i} $
{for} $ i=1, \ldots , s+1 $ are the fitting nodes, and
one of them is required to be one.

As an {illustrative} example,  we consider $s=1,\ c_1=0,\
c_2=1$ and choose \begin{equation*}\label{one}
 \begin{aligned}  &C_{\tau}=(1-\tau)I+ \tau \varphi_{1}(h\Omega),\ A_{\tau\sigma}=\tau \varphi_{2}(h\Omega),\
\bar{B}_{\tau}=\varphi_{2}(h\Omega),\    B_{\tau}=\varphi_{1}(h\Omega).
\end{aligned}
\end{equation*}
This obtained algorithm can be rewritten as
\begin{equation}\label{EAVF}
\begin{array}[c]{ll}x_{n+1}=x_{n}+
h\varphi_1(h\Omega) v_{n}+h^2
\varphi_2(h\Omega)  \int_{0}^1
F\big(x_{n}+\sigma(x_{n+1}-x_{n}) \big) d\sigma,\\
v_{n+1}=\varphi_0(h\Omega)v_{n}+h \varphi_{1}(
h\Omega) \int_{0}^1
F\big(x_{n}+\sigma(x_{n+1}-x_{n}) \big) d\sigma,
\end{array}
\end{equation}
which is denoted by EM1.
\end{algorithm}
\begin{rem}
It is noted that EM1 has been given  in \cite{Wang2020} and it was shown to be energy-preserving. However, its
convergence has not been studied there. In this paper, we will
analyse the convergence of each algorithm. It will be shown that
M1 has a first-order uniform error bound in both $x$ and $v$ and
the others are of order two and have a uniform convergence  in $x$  for
$0<\eps\ll 1$. In contrast, many  classical methods such as Euler methods, Runge-Kutta (-Nystr\"{o}m) methods often show non-uniform error bounds in both $x$ and $v$, where the error constant is usually proportional to $1/\eps^k$ for some $k>0$.
\end{rem}

%

%
\begin{rem}
It is remarked that the
following  integrators for solving \eqref{charged-particle sts-cons} has been given in \cite{lubich19}
\begin{equation*}
\begin{array}[c]{ll} x_{n+1}=x_{n}+
h\varphi_1(h \Omega) v_{n}+\frac{1}{2}h^2\Psi(h \Omega) F_n  ,\\
v_{n+1}=\varphi_0(h \Omega)v_{n}+\frac{1}{2}h\big(\Psi_0(h \Omega) F_n
+\Psi_1(h \Omega) F_{n+1} \big),
\end{array}
\end{equation*}
where $F_n=F(x_n)$ and  $\Psi, \Psi_0,\Psi_1$ are matrix-valued and
bounded functions of $h \Omega$ satisfying
$\Psi(0)=\Psi_0(0)=\Psi_1(0)=1$.   However, only convergence is researched there and the structure preservation such as symplecticity or energy conservation has not been discussed.

\end{rem}
\section{Main results and a numerical test}\label{sec: main results}
\subsection{Main results}
The main results of this paper are given by the following four
theorems. The first three theorems are about structure preservations
and the last one {concerns} uniform error bound.

\begin{theo}\label{thm: 1} (\textbf{Symplecticity of SM1-SM3}.)
Consider the methods SM1-SM3 {where} $p_{n+1} = v_{n+1}
-\frac{1}{2\varepsilon} \tilde{B}x_{n+1}$. {In this
case,} for the non-separable Hamiltonian \eqref{H}, the map
$(x_n,p_n ) \rightarrow (x_{n+1},p_{n+1})$ determined by these
methods is symplectic.
\end{theo}

\begin{theo}\label{thm: 2} (\textbf{Energy preservation of EM1} \cite{Wang2020}.)
{The method EM1 preserves} the energy \eqref{energy of
cha} exactly, i.e. $E(x_{n+1},v_{n+1})=E(x_{n},v_{n}), \
n=0,1,\ldots.$
\end{theo}

\begin{theo}\label{thm: 3} (\textbf{Long time energy conservation of M2 and SM1-SM3}.)
Consider the following assumptions.
\begin{itemize}
\item  It is assumed
that  the initial values $x_{0}$ and $v_0:=\dot{x}_0$ are bounded such that the
energy $E$ is bounded independently of $\varepsilon$ along the
solution.

\item {Suppose} that the considered numerical solution    stays in a compact
set.

\item A lower bound on the  stepsize $h/\varepsilon\geq c_0 > 0$ is
required.

\item Assume that the numerical non-resonance condition is true
%
\begin{equation*}
|\sin(\frac{h}{2}(k\cdot \tilde{\Omega}))| \geq c \sqrt{h}\ \
\mathrm{for} \ \ k \in \mathbb{Z}^l\backslash \mathcal{M}
   \ \   \mathrm{with} \ \  |k|\leq N
\end{equation*}
for some $N\geq2$ and $c>0$. The notations used here are referred to the last part of Section \ref{sec: conservations}.
\end{itemize}
For the   methods M2 and SM1-SM3, it holds that
\begin{equation}\label{long EP}
\begin{aligned}
\abs{E(x_{n},v_{n})-E(x_{0},v_{0})}\leq C h
\end{aligned}
\end{equation}
for $0\leq nh\leq h^{-N+1}.$ The constant $C$ is independent of $n, h,\varepsilon$, but depends
on $N, T$ and the constants in the assumptions.
\end{theo}
\begin{rem}
It is noted that  M1  does not have the
above energy conservation property. The reason is that it is not a
symmetric method. {It will be seen from the proof  given
in Section \ref{sec: conservations}}  that symmetry plays an
important role in the analysis.
\end{rem}

\begin{theo}\label{thm: 4}
 (\textbf{Convergence}.)
 For the  methods M1-M2 and the energy-preserving method EM1, under the condition that $h\leq C_1\eps$, the global errors are bounded by
\begin{subequations}\label{thm err}
\begin{align}
\textmd{M1:}\  \ &\abs{x_{n} -x(t_n)}+\abs{v_{n} -v(t_n)}\lesssim  h, \label{con M1}\\
\textmd{M2 and EM1:}\  \ &\abs{x_{n} -x(t_n)}\lesssim  h^2,\ \abs{v_{n} -v(t_n)} \lesssim  h^2/\varepsilon, \label{con M2}
\end{align}
\end{subequations}
where $0<nh\leq T$.   Here we denote $A\lesssim B$ for $A\leq CB$ with a generic constant $C>0$ independent of $h$ or $n$ or $\eps$ but depends on $T$ and $C_1$.

For the symplectic methods SM1-SM3,  under the conditions of Theorem \ref{thm: 3}, the global errors are
\begin{equation}
\textmd{SM1-SM3:}\  \ \ \abs{x_{n} -x(t_n)}\lesssim  h^2,\ \abs{v_{n} -v(t_n)} \lesssim  h^2/\varepsilon, \label{con SM}
\end{equation}
where the error constants are independent of $n, h,\varepsilon$, but depend
on $T$ and the constants in the assumptions of Theorem \ref{thm: 3}.
\end{theo}

\begin{table}[t!]
\newcommand{\tabincell}[2]{\begin{tabular}{@{}#1@{}}#2\end{tabular}}  
\begin{tabular}{|c|c|c|c|c|c|c|}
\hline
\text{Methods}& \text{Symplecticity} & \text{Symmetry} &\text{Energy property} &\text{Convergence}        \\
\hline
  \text{M1} & \text{No}& \text{No}   &\text{No}               &   \tabincell{c}{$\abs{x_{n} -x(t_n)} \lesssim  h$\\ $\abs{v_{n} -v(t_n)}\lesssim  h$}\cr
  \hline
  \text{M2} & \text{No} & \text{Yes}& \text{Near conservation}   &   \tabincell{c}{$\abs{x_{n} -x(t_n)} \lesssim  h^2$\\ $\abs{v_{n} -v(t_n)}\lesssim  h^2/\varepsilon$}\cr
  \hline
  \text{SM1} & \text{Yes}& \text{Yes}   &\text{Near conservation}               &   \tabincell{c}{$\abs{x_{n} -x(t_n)} \lesssim  h^2$\\ $\abs{v_{n} -v(t_n)}\lesssim  h^2/\varepsilon$}\cr
  \hline
  \text{SM2} & \text{Yes} & \text{Yes}& \text{Near conservation}   &   \tabincell{c}{$\abs{x_{n} -x(t_n)} \lesssim  h^2$\\ $\abs{v_{n} -v(t_n)}\lesssim  h^2/\varepsilon$}\cr
  \hline
 \text{SM3} & \text{Yes}& \text{Yes} &\text{Near conservation}     &    \tabincell{c}{$\abs{x_{n} -x(t_n)} \lesssim  h^2$\\ $\abs{v_{n} -v(t_n)}\lesssim  h^2/\varepsilon$}\cr
  \hline
 \text{EM1} & \text{No} & \text{Yes} &\text{Exact conservation}     &    \tabincell{c}{$\abs{x_{n} -x(t_n)} \lesssim  h^2$\\ $\abs{v_{n} -v(t_n)}\lesssim  h^2/\varepsilon$}\cr
  \hline
\end{tabular}
\caption{Properties of the obtained methods.} \label{praERKN}
\end{table}
For the six methods presented in this paper, concerning the symmetry \cite{hairer2006}, it is easy to check that all of them except M1 are symmetric. Their properties are
summarized in Table \ref{praERKN}. The main observation of the paper is that M1 has a first uniform (in $\eps$) error bound in both $x$ and $v$ and
 symmetric or symplectic or energy-preserving methods show second-order uniform  error bound in $x$. Moreover, SM1 can preserve the energy exactly and symplectic or symmetric  methods have a good near conservations of energy  over long times. All of these observations will be numerically illustrated by a test given blew.

\subsection{A numerical test}
 As  an illustrative numerical experiment, we consider the charged particle system
of  \cite{Hairer2017-1} with an additional factor $1/\varepsilon$
and a constant magnetic field. The system can be expressed by
\eqref{charged-particle sts-cons} with $d=3$, where the potential
$U(x)=x_{1}^{3}-x_{2}^{3}+x_{1}^{4}/5 +x_{2}^{4}+x_{3}^{4}$ and
$\tilde{B}= \begin{pmatrix}
                     0 & 0.2 & 0.2 \\
                     -0.2 & 0 & 1 \\
                     -0.2 & -1 & 0 \\
                  \end{pmatrix}.$ The initial
values are chosen as $x(0)=(0.6,1,-1)^{\intercal}$ and
$v(0)=(-1,0.5,0.6)^{\intercal}.$
 \begin{figure}[t!]
\centering
\includegraphics[width=4.2cm,height=4.2cm]{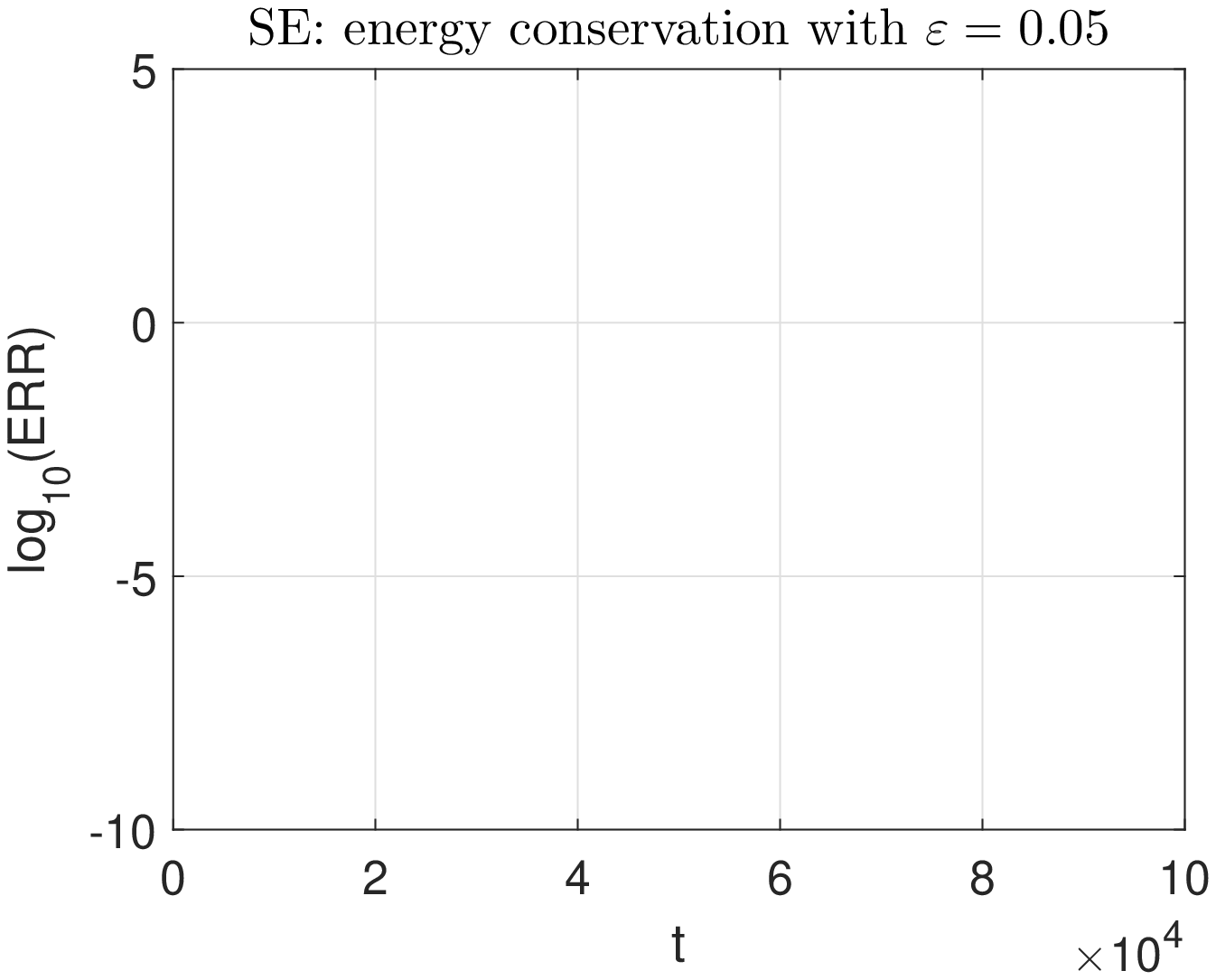}
\includegraphics[width=4.2cm,height=4.2cm]{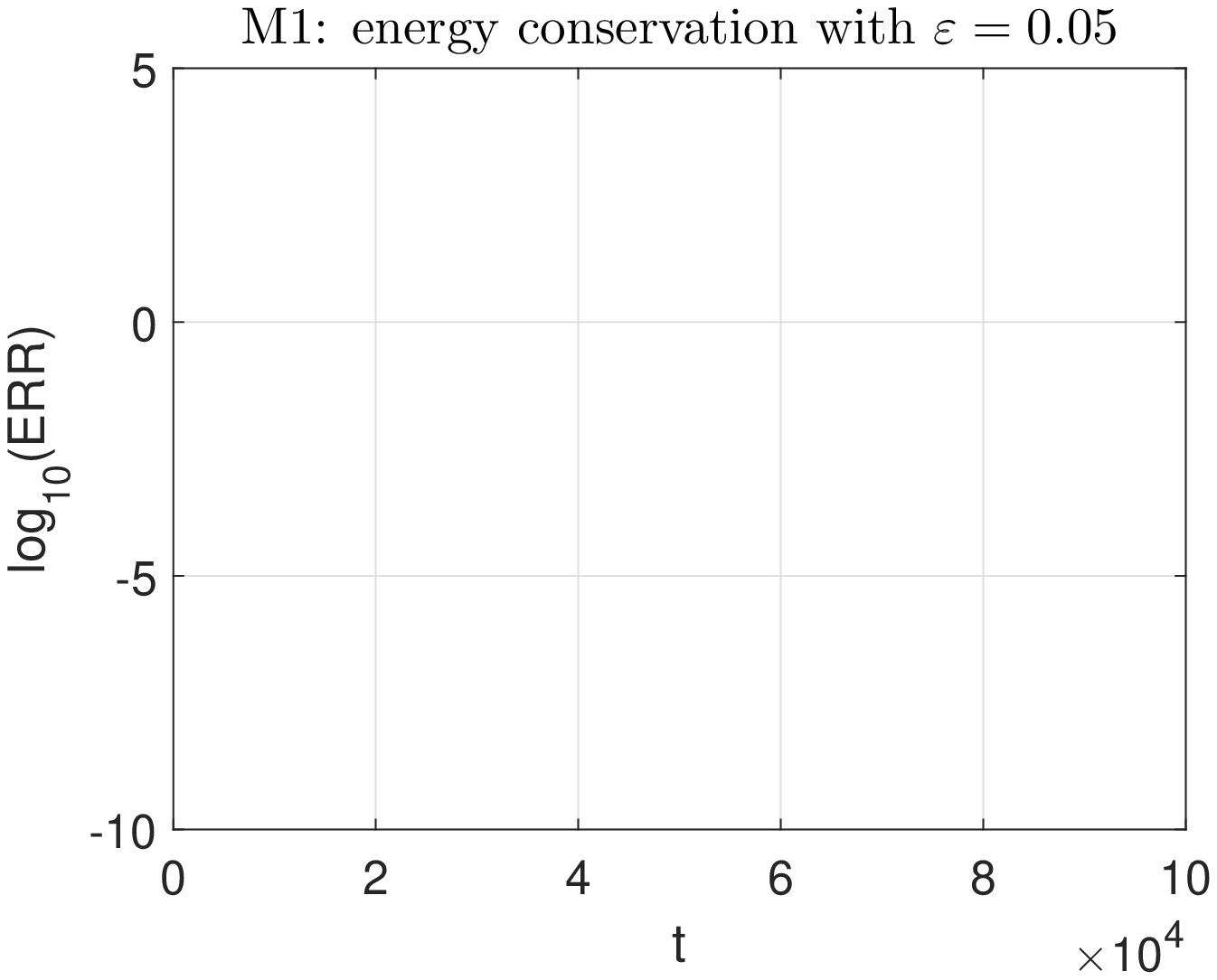}
\includegraphics[width=4.2cm,height=4.2cm]{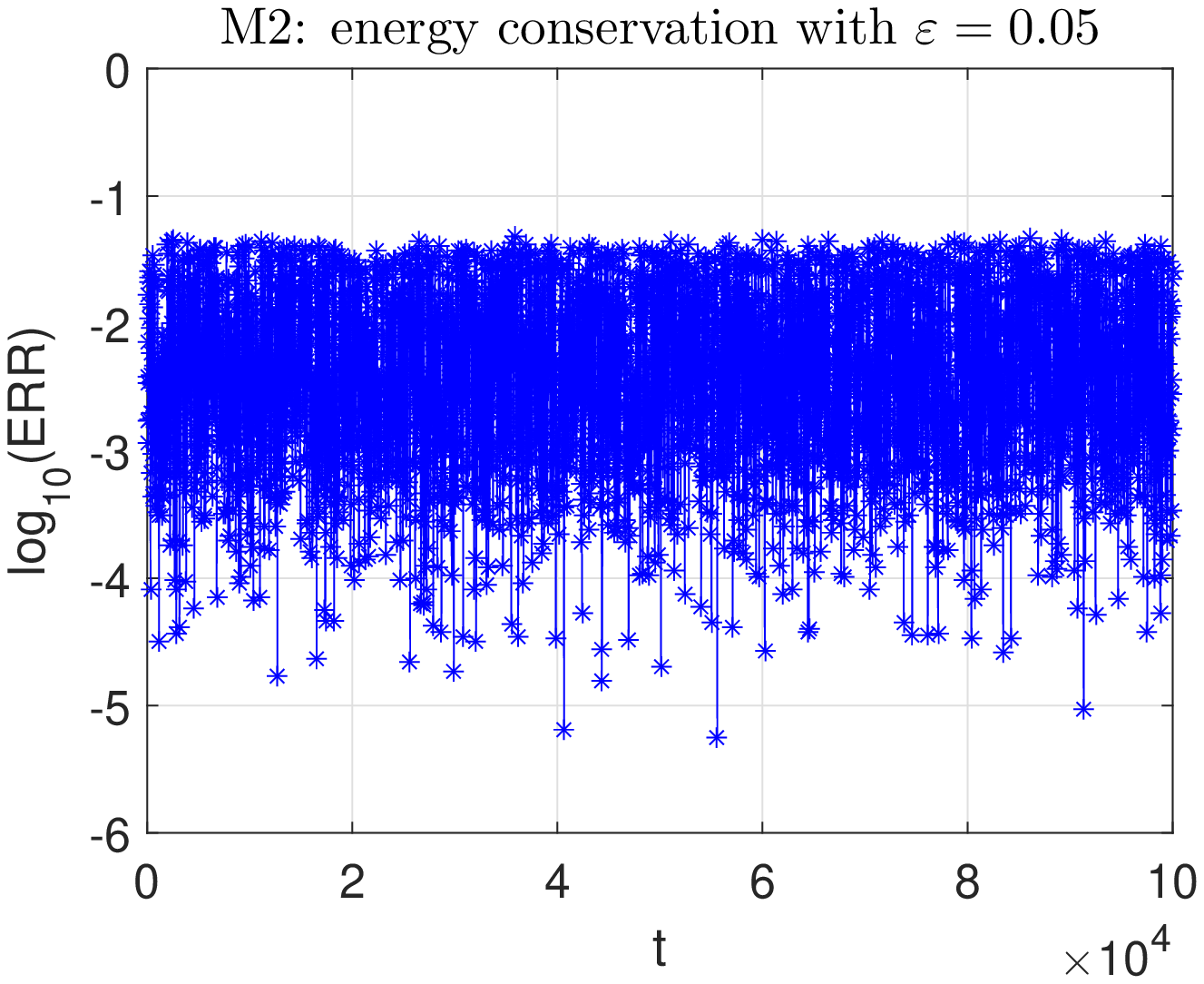}
\includegraphics[width=4.2cm,height=4.2cm]{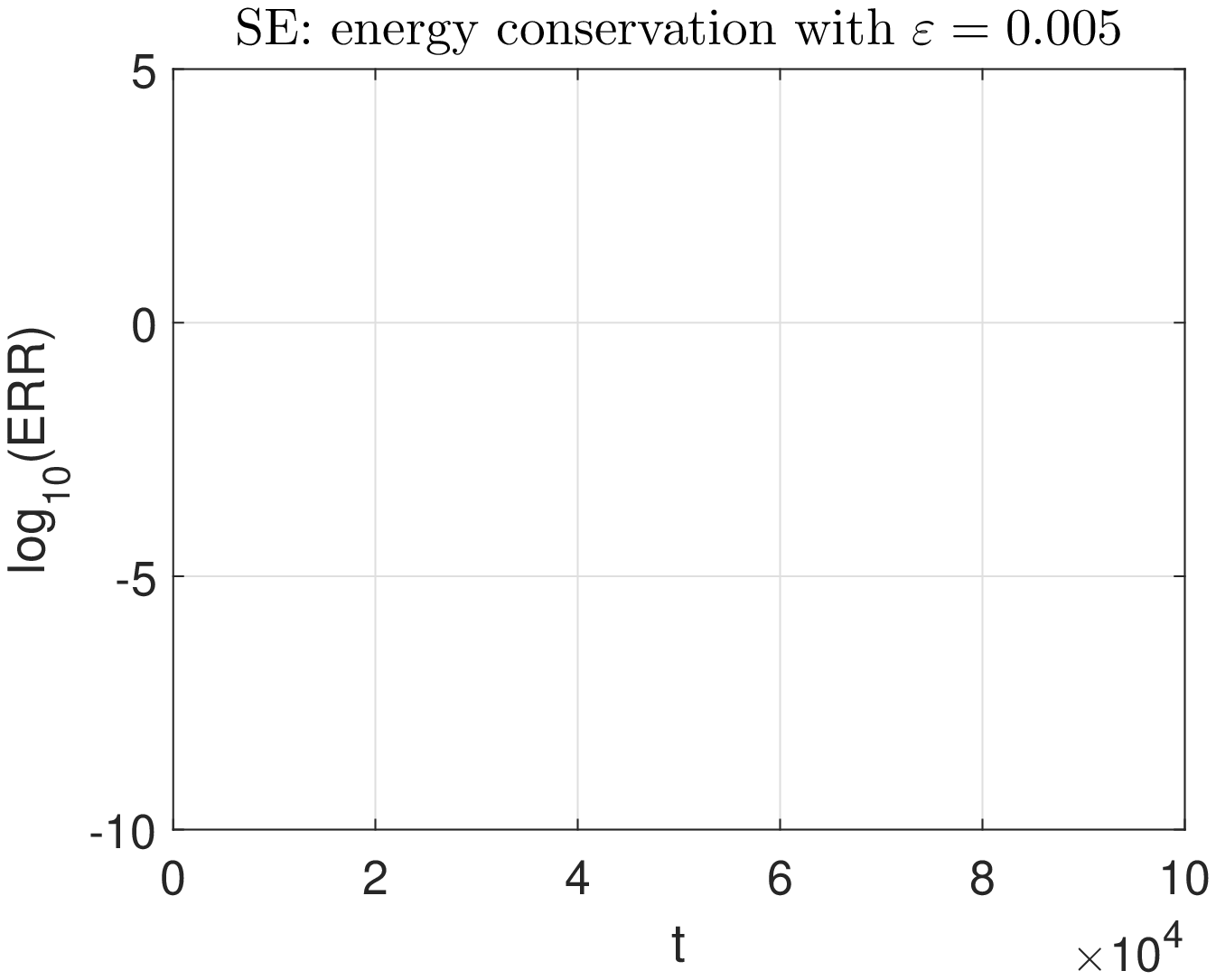}
\includegraphics[width=4.2cm,height=4.2cm]{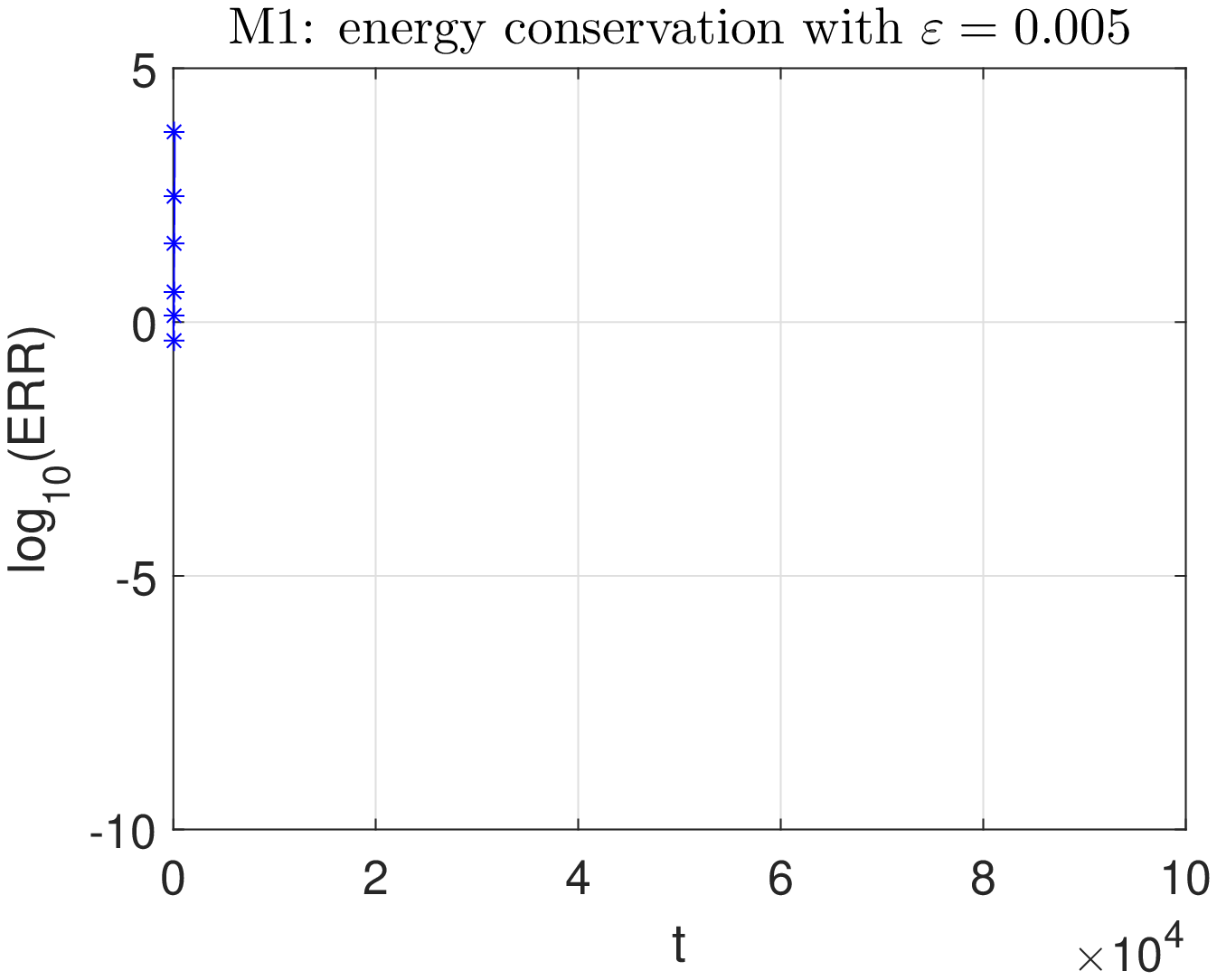}
\includegraphics[width=4.2cm,height=4.2cm]{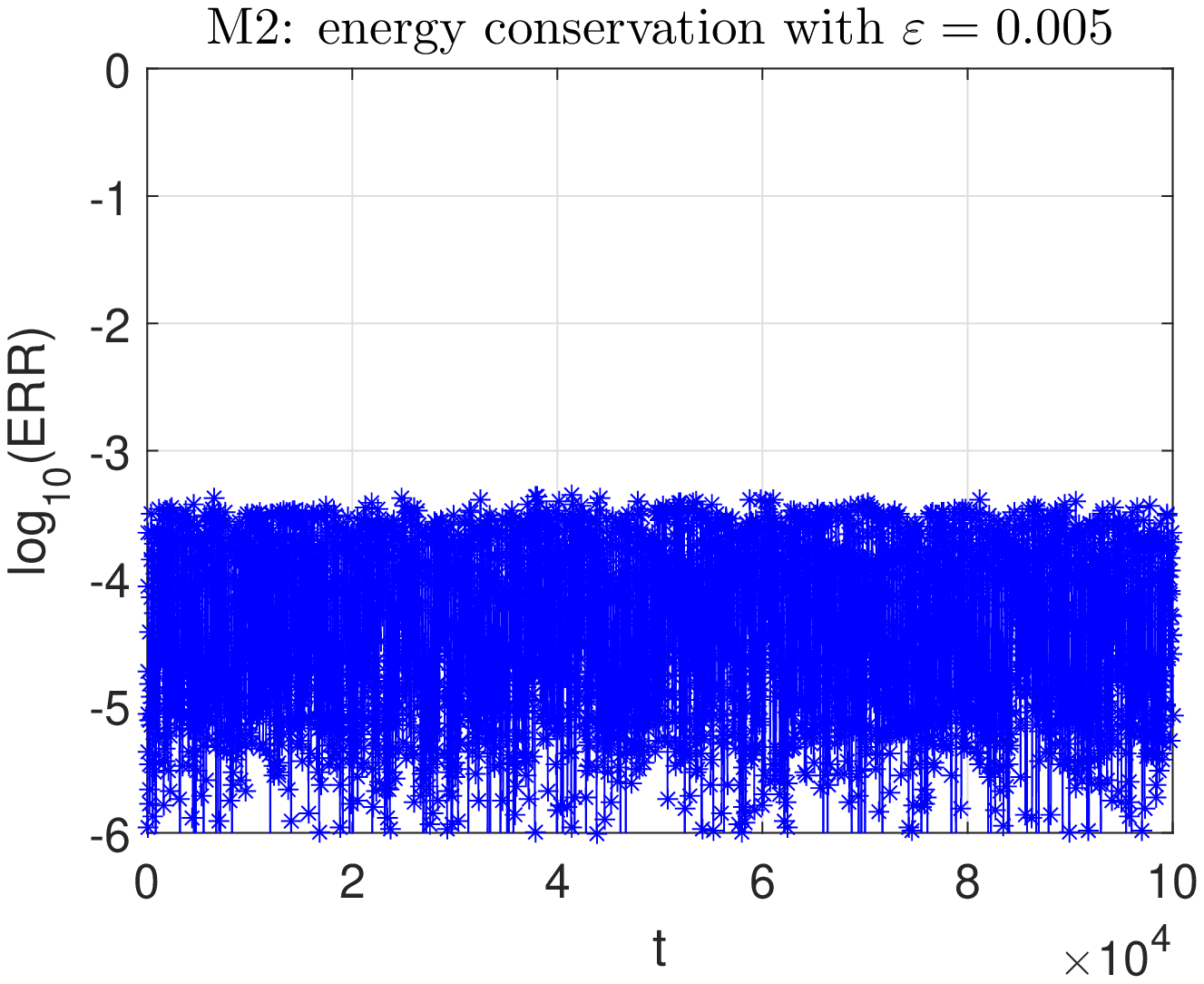}
\caption{The relative energy errors (ERR) against $t$ for SE and M1-M2.} \label{p1}
\end{figure}
 \begin{figure}[t!]
\centering
\includegraphics[width=4.2cm,height=4.2cm]{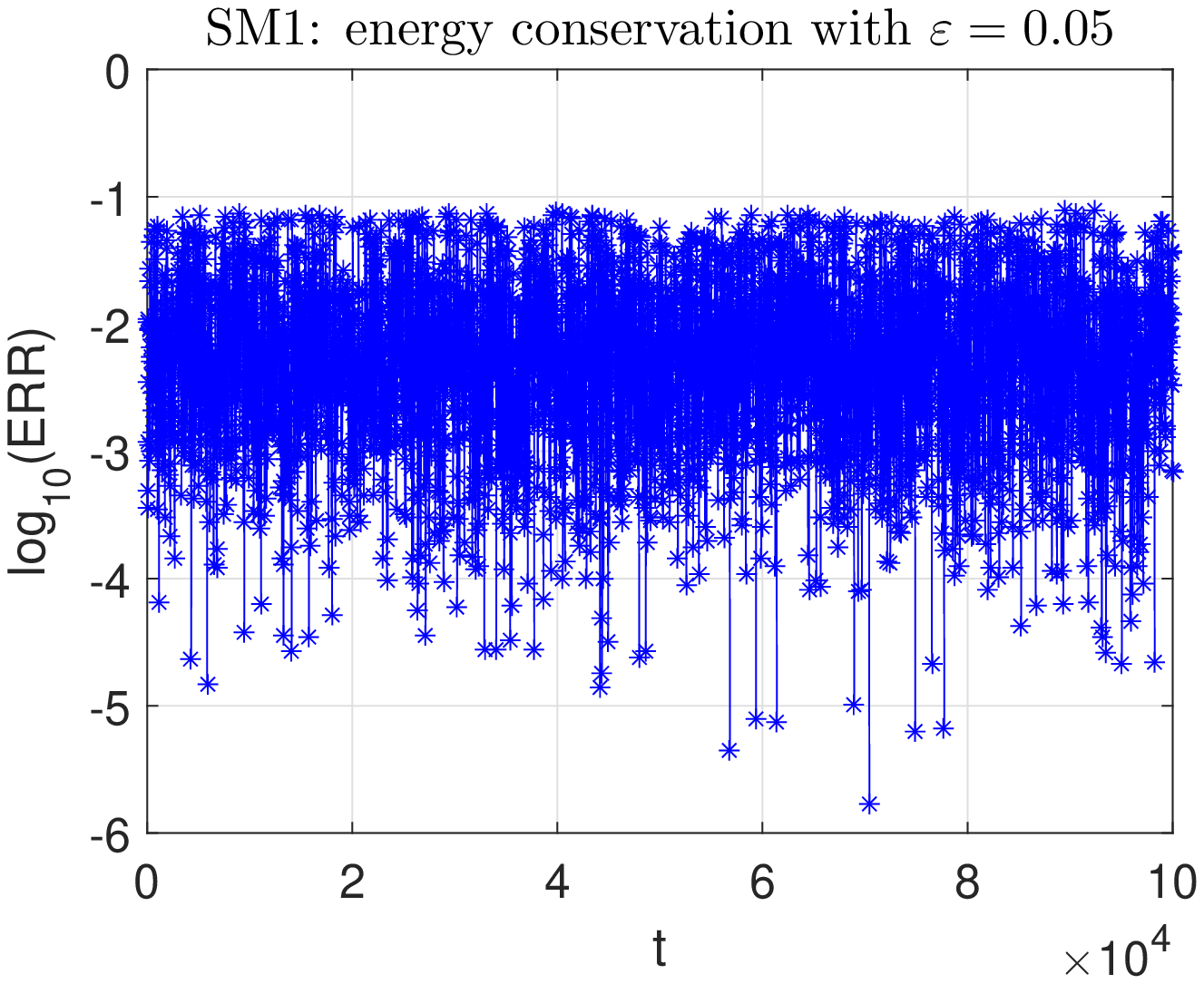}
\includegraphics[width=4.2cm,height=4.2cm]{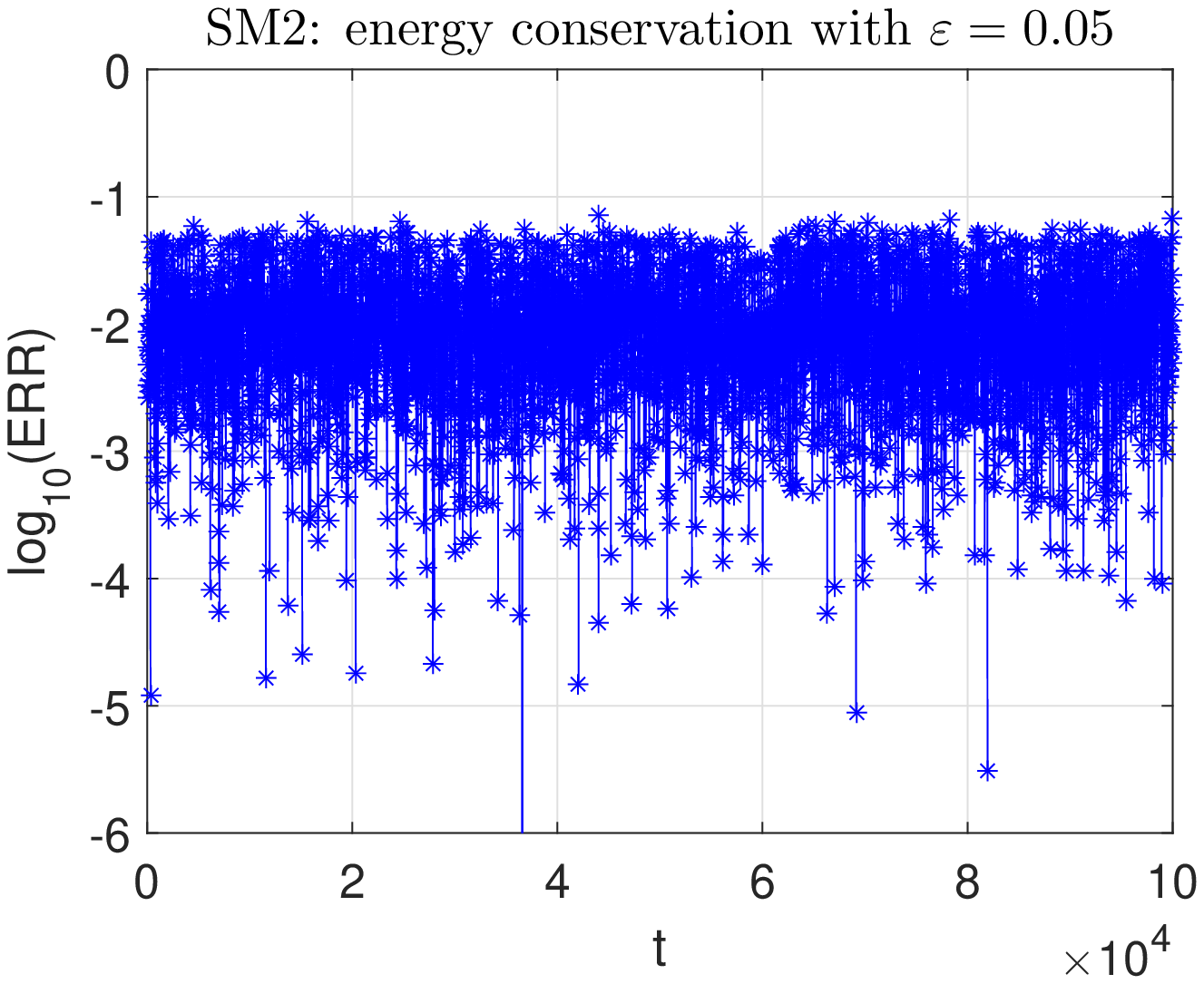}
\includegraphics[width=4.2cm,height=4.2cm]{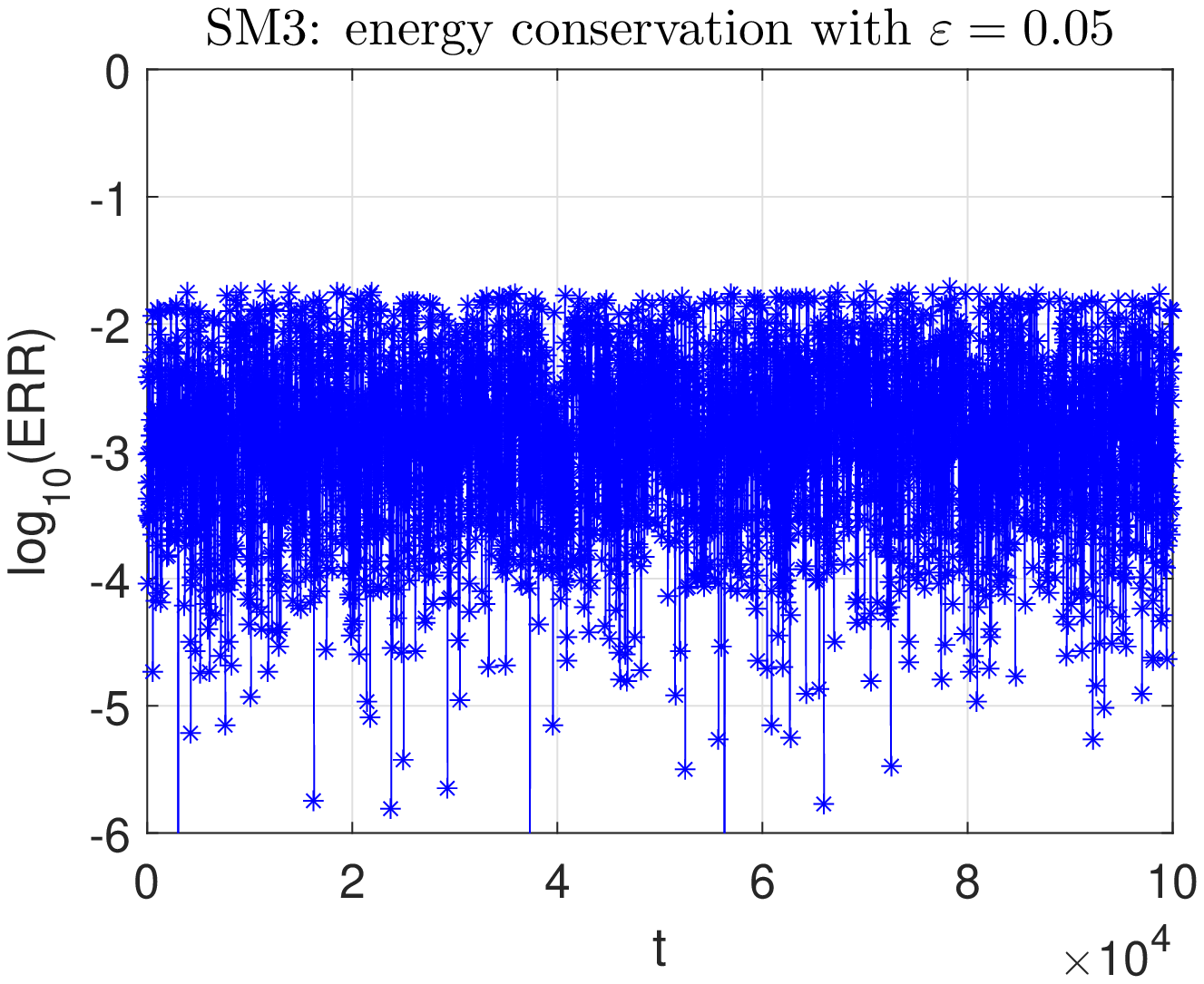}
\includegraphics[width=4.2cm,height=4.2cm]{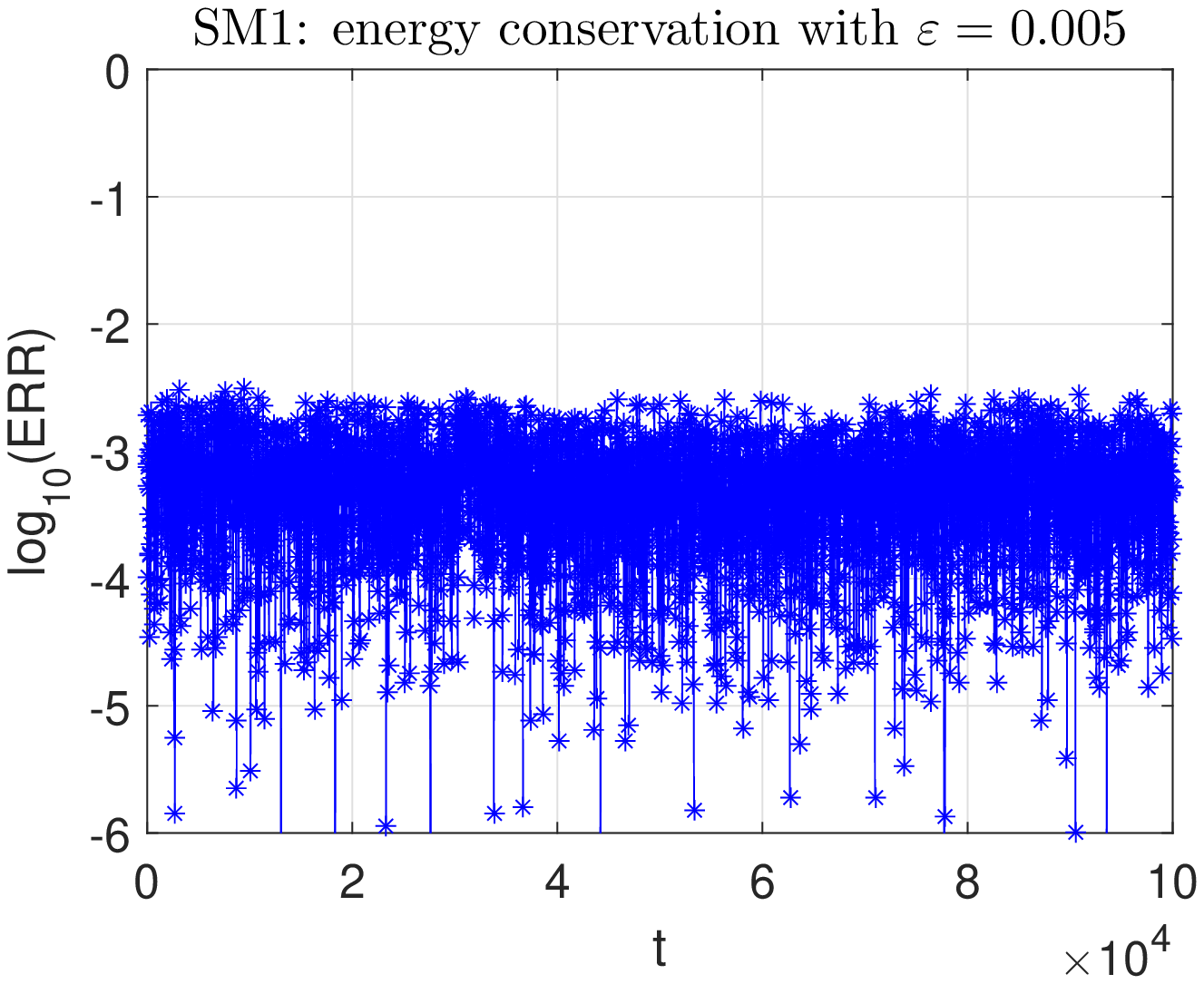}
\includegraphics[width=4.2cm,height=4.2cm]{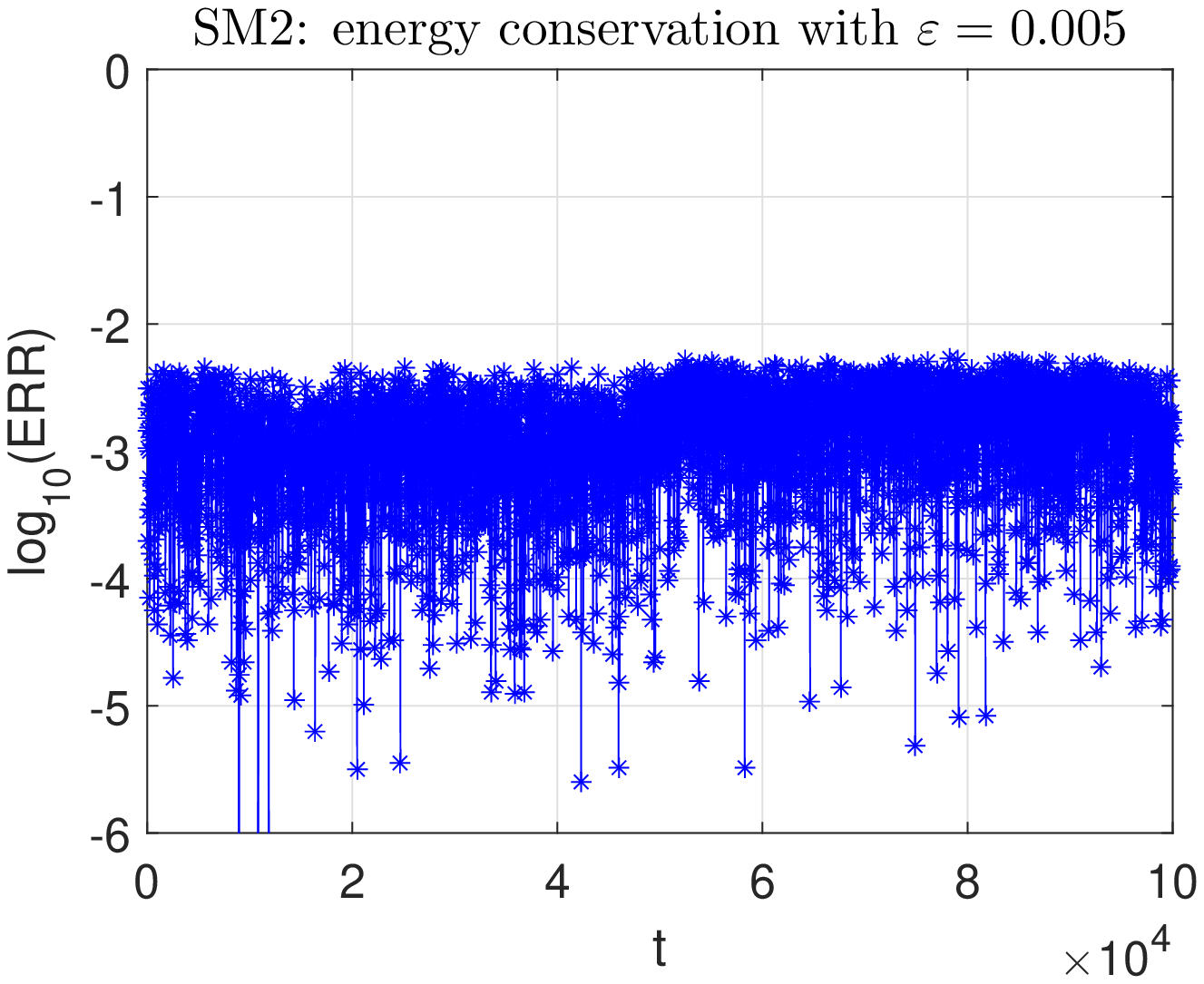}
\includegraphics[width=4.2cm,height=4.2cm]{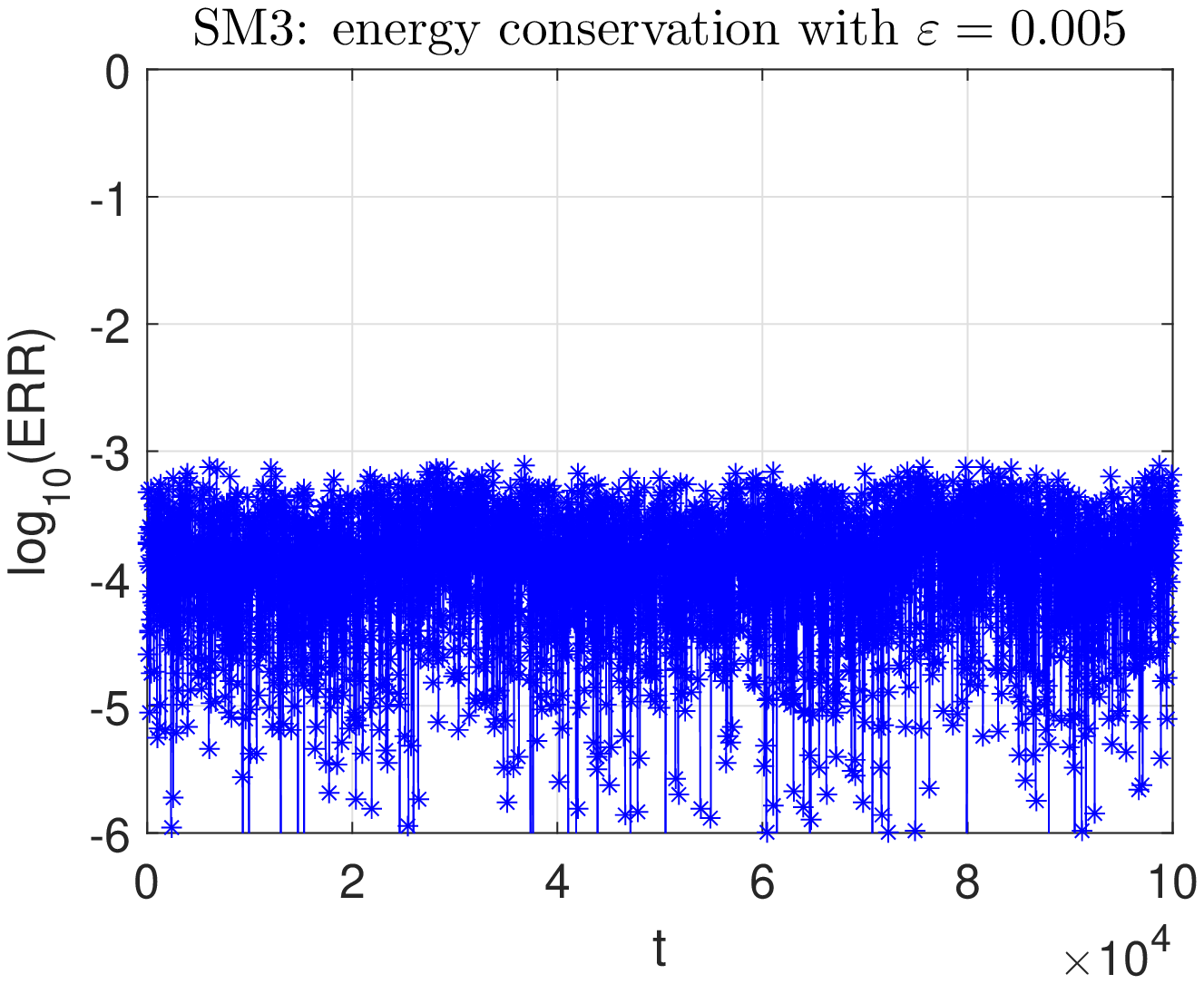}
\caption{The relative energy errors (ERR) against $t$ for symplectic SM1-SM3.} \label{p2}
\end{figure}

 \begin{figure}[t!]
\centering
\includegraphics[width=4.2cm,height=4.2cm]{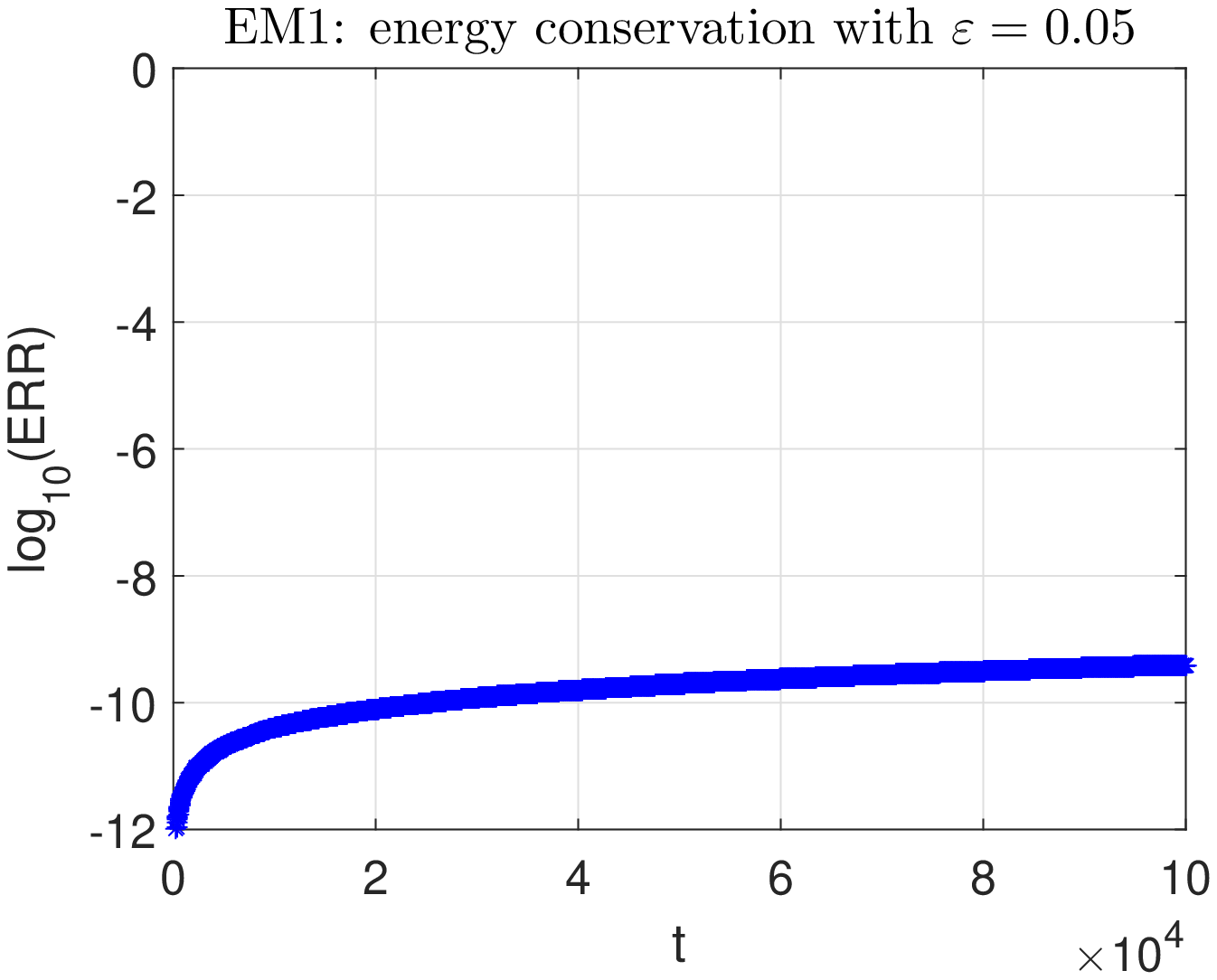}
\includegraphics[width=4.2cm,height=4.2cm]{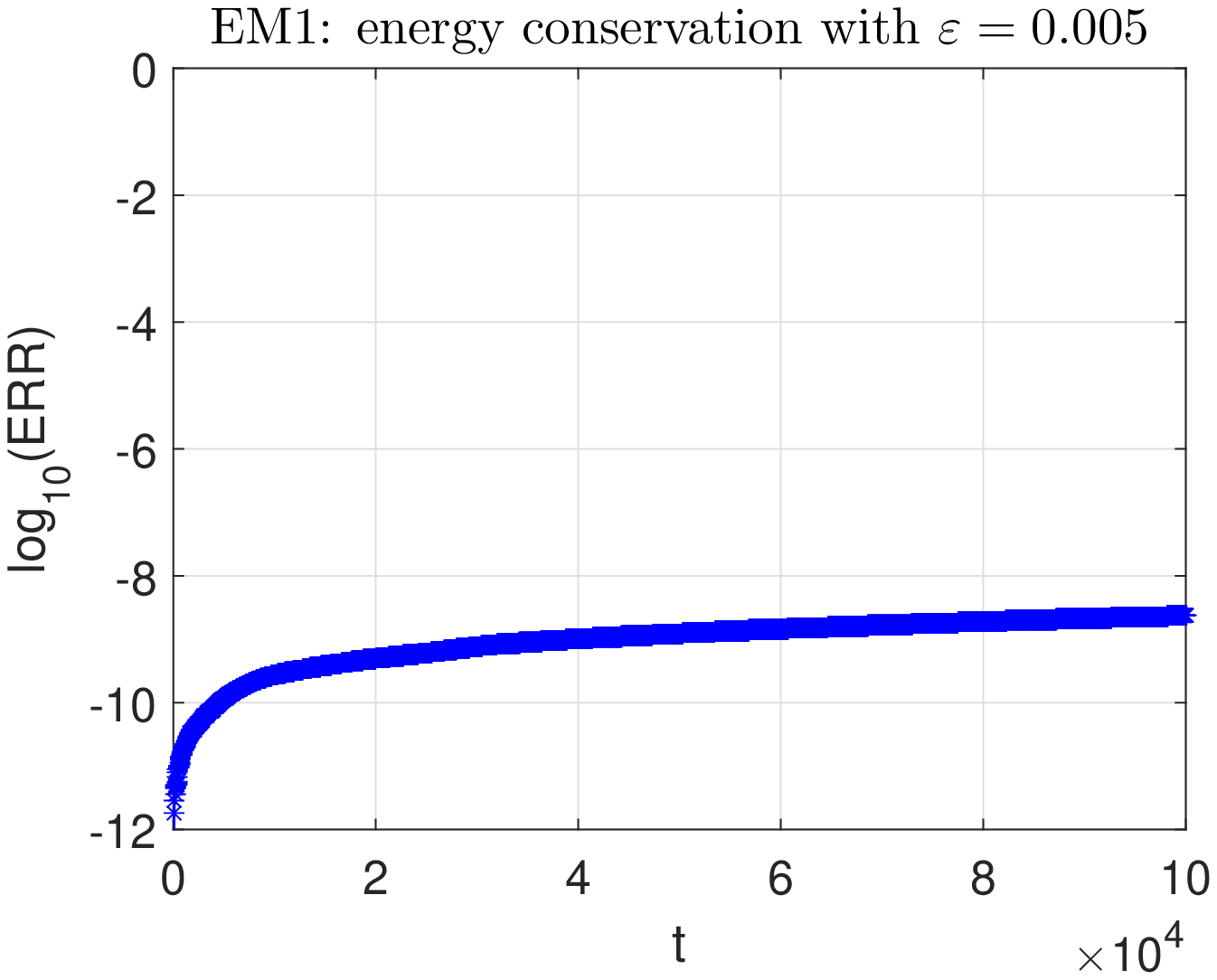}
\caption{The relative energy errors (ERR) against $t$ for energy-preserving  EM1.} \label{p3}
\end{figure}
\subsubsection{Energy conservation} We take $\varepsilon=0.05,0.005$ and apply our six methods as well as the symplectic Euler method (denoted by SE)
to this problem on $[0,100000]$ with $h=\eps$. The  standard fixed point iteration is used for EM1 and we set $10^{-16}$ as the error tolerance and $10$ as the maximum number of iterations.
The
relative errors $ERR:=(E(x_{n},v_{n})-E(x_{0},v_{0}))/E(x_{0},v_{0})$ of the energy
 are
displayed in Figs.  \ref{p1}-\ref{p3}. We do not show the figure if the error is too large. According to  these results, we have the following
observations. M2 and SM1-SM3 (Figs.  \ref{p1}-\ref{p2}) have near energy conservation over long times, EM1 preserves the energy very well (Fig. \ref{p3}) but SE and M1 show a bad energy conservation (Fig. \ref{p1}).

\begin{figure}[t!]
\centering
\includegraphics[width=4.2cm,height=4.2cm]{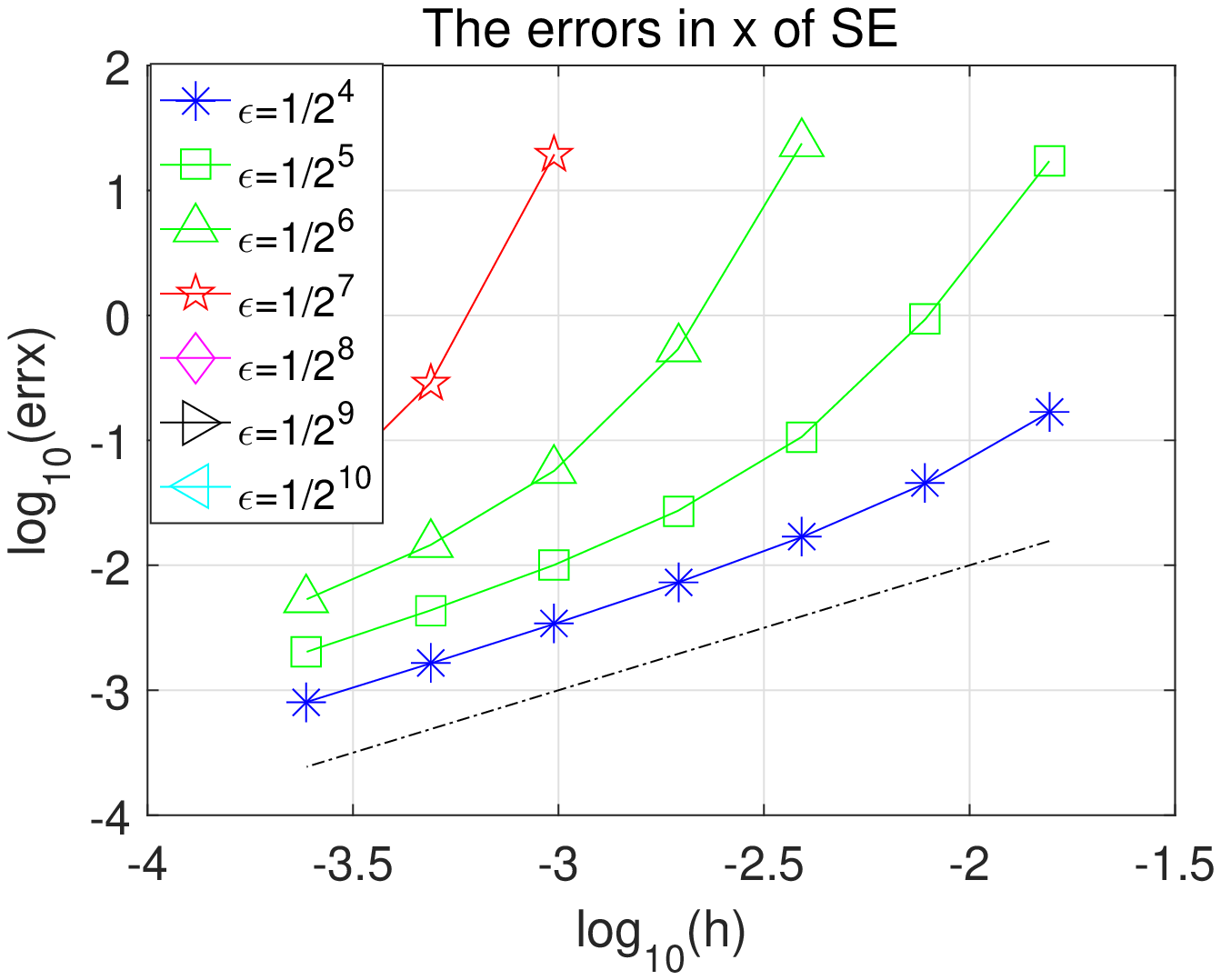}
\includegraphics[width=4.2cm,height=4.2cm]{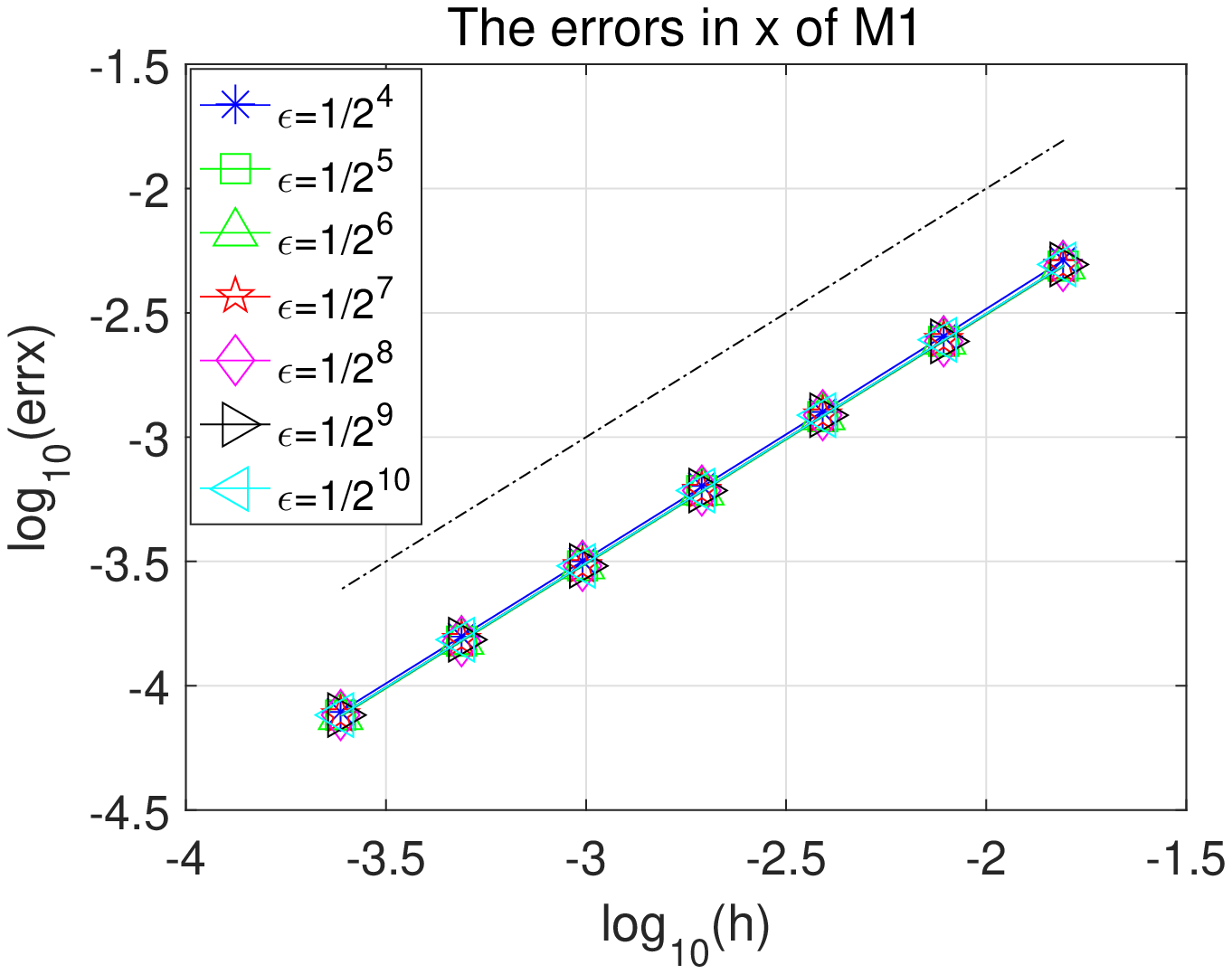}
\includegraphics[width=4.2cm,height=4.2cm]{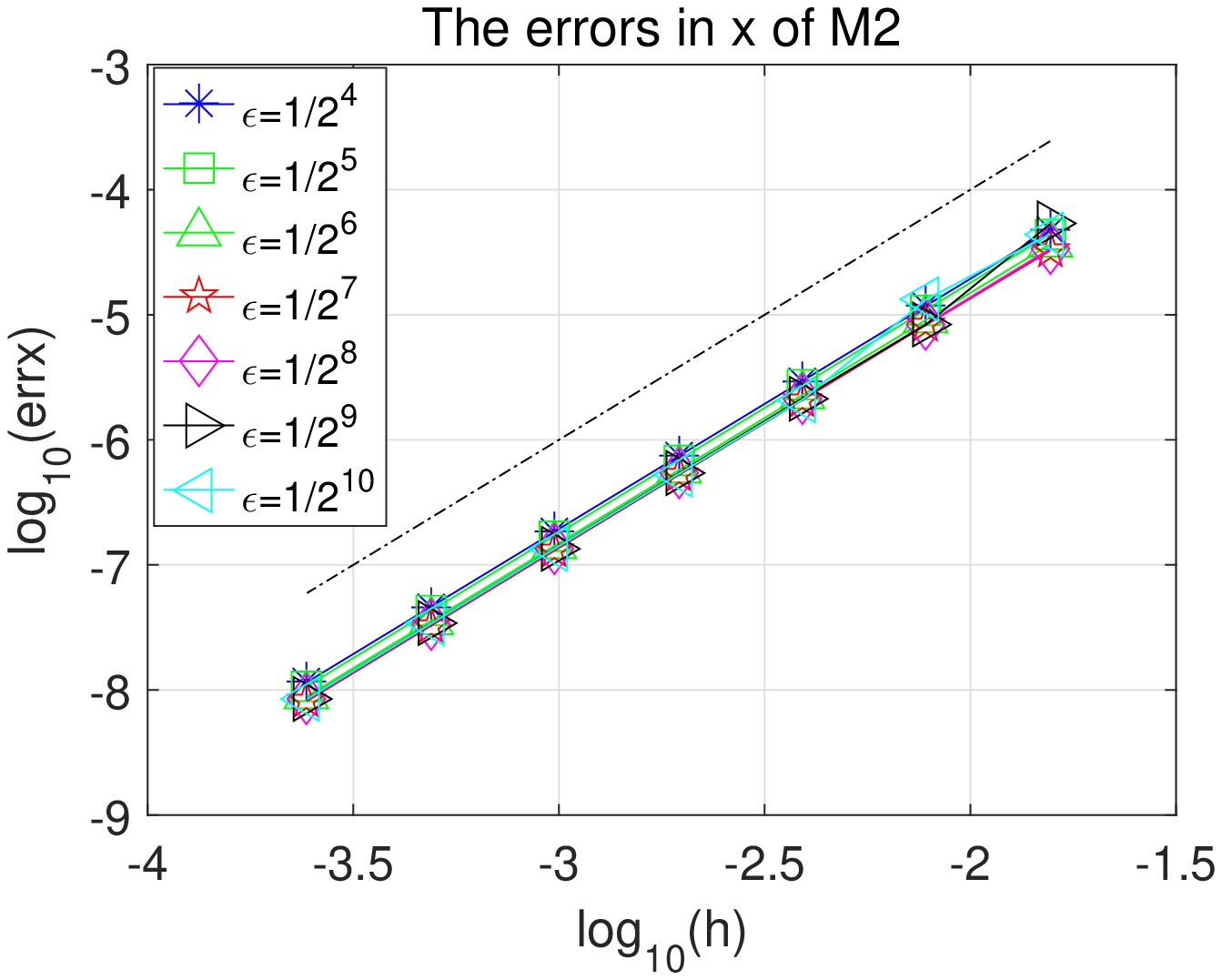}
\includegraphics[width=4.2cm,height=4.2cm]{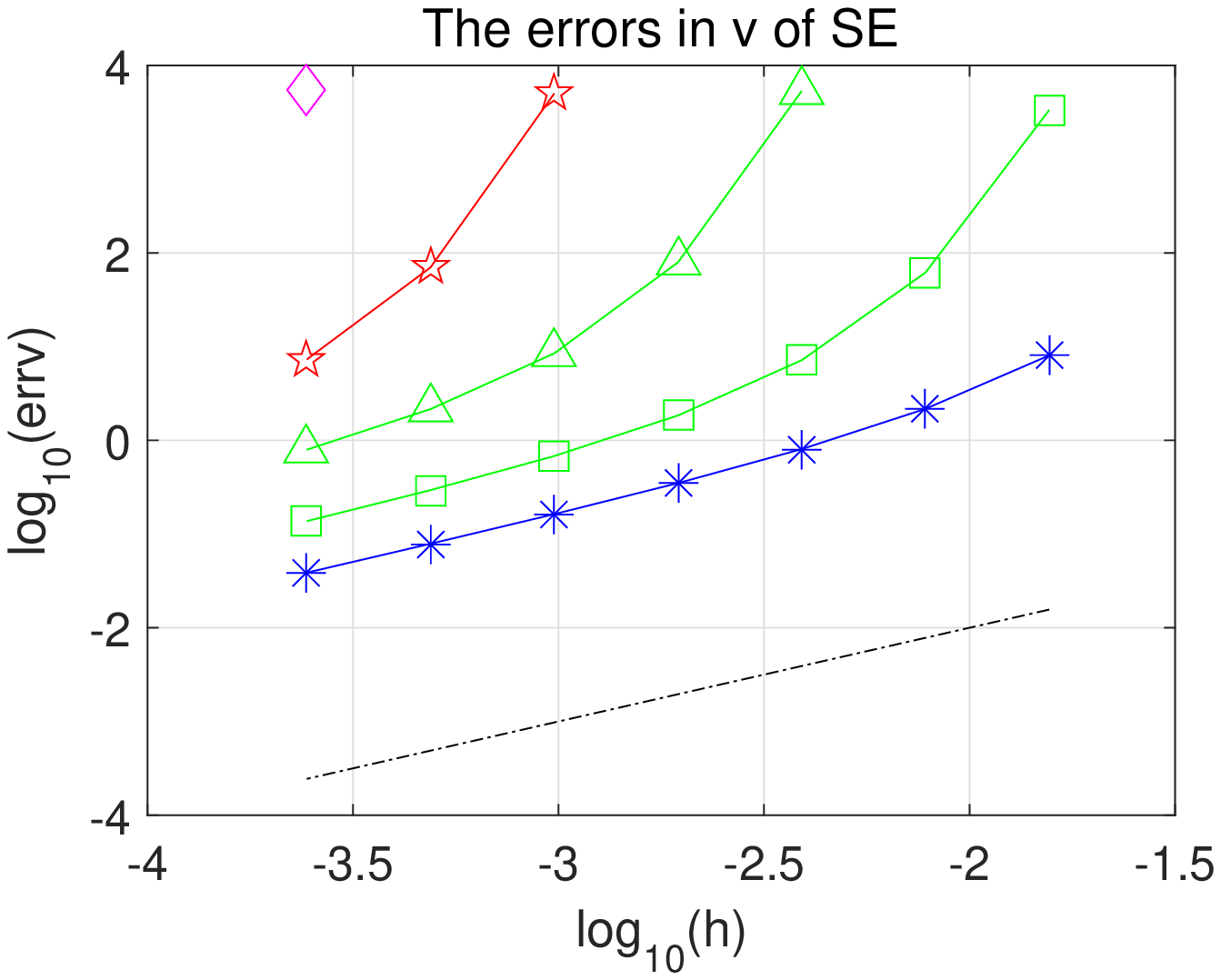}
\includegraphics[width=4.2cm,height=4.2cm]{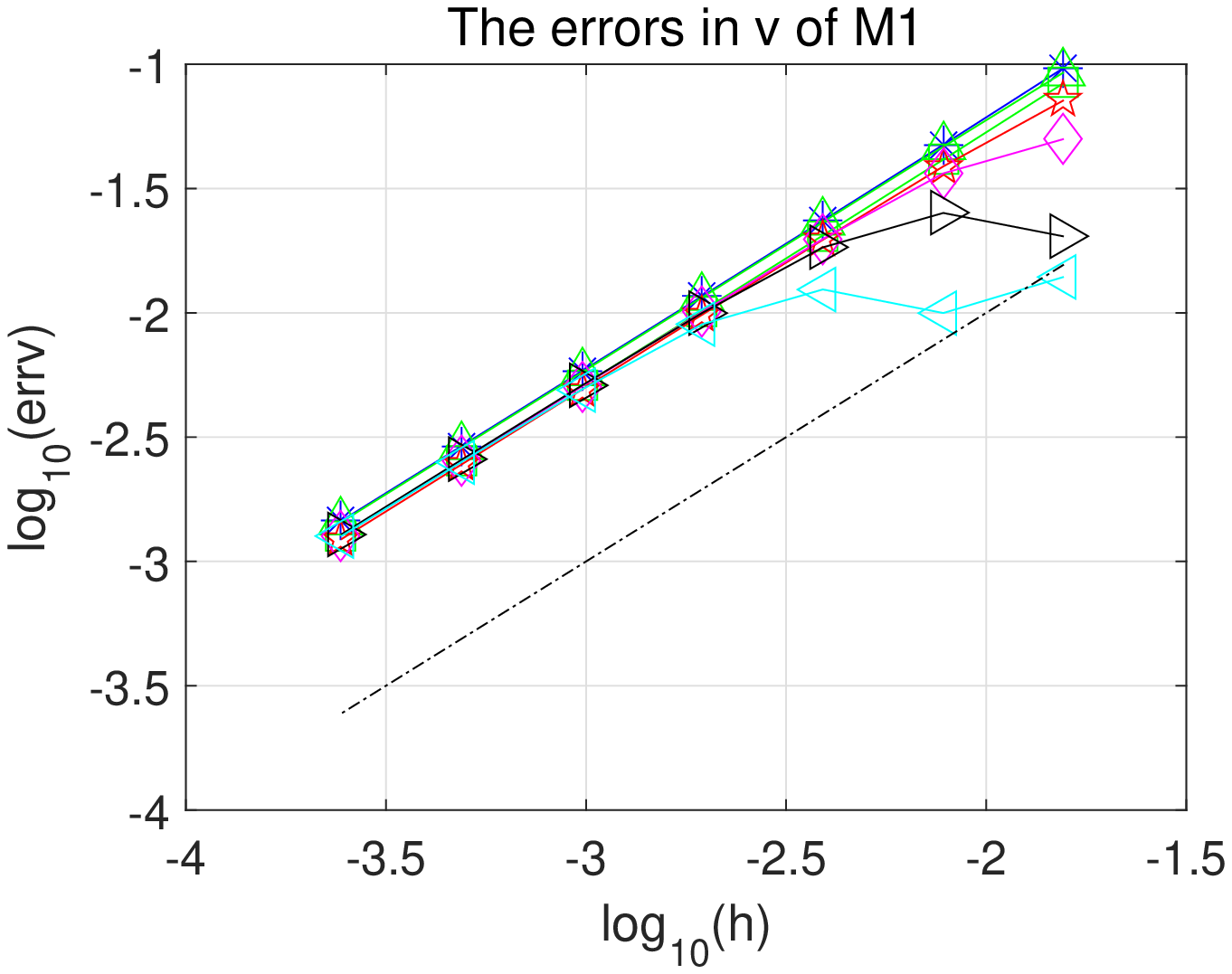}
\includegraphics[width=4.2cm,height=4.2cm]{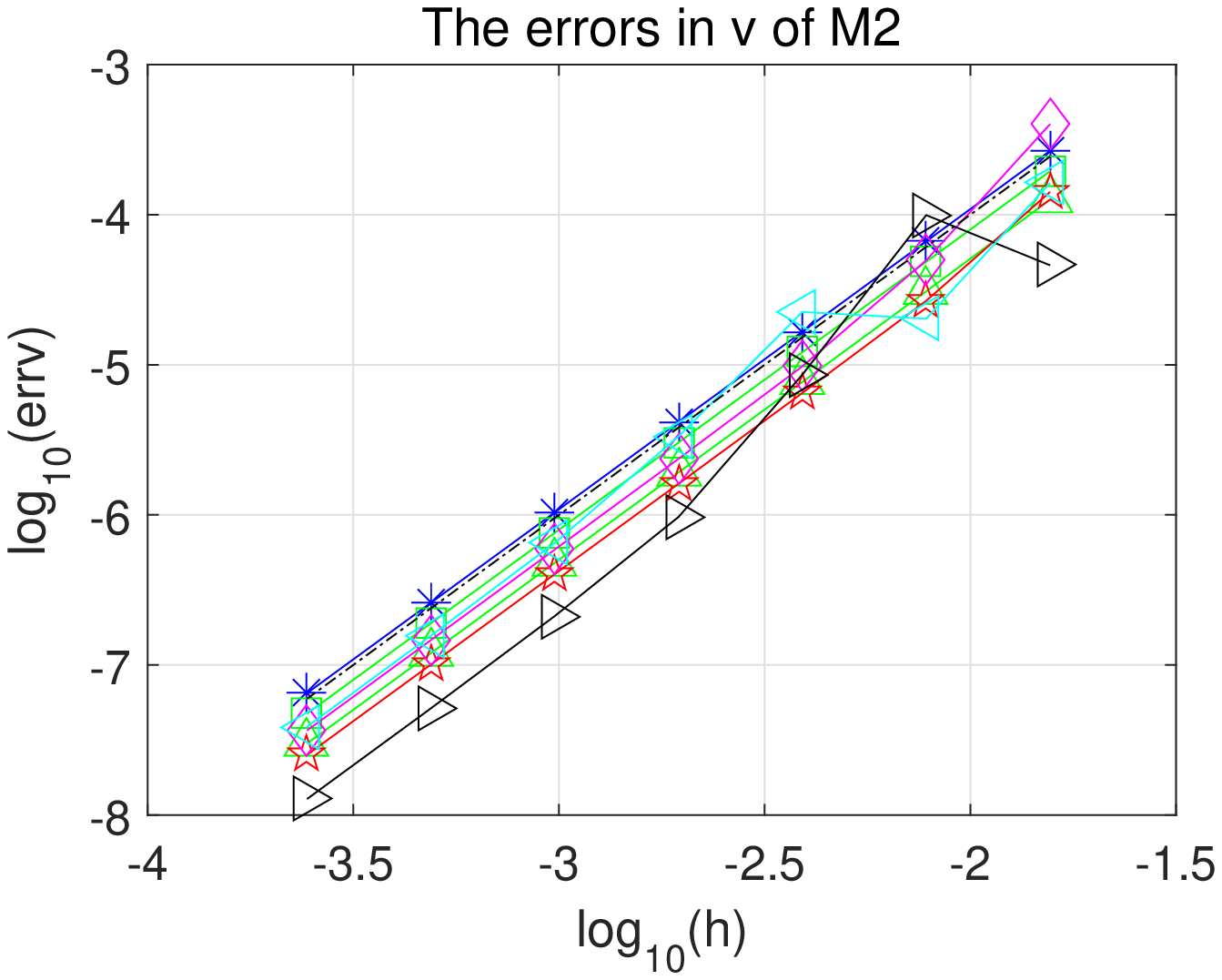}
\caption{The   errors in $x$ (errx) and $v$ (errv) against $h$ for SE and M1-M2 (the slope of the dotted line for SE and M1 is one and for M2 is two).} \label{p4}
\end{figure}

\subsubsection{Convergence}
\begin{figure}[t!]
\centering
\includegraphics[width=4.2cm,height=4.2cm]{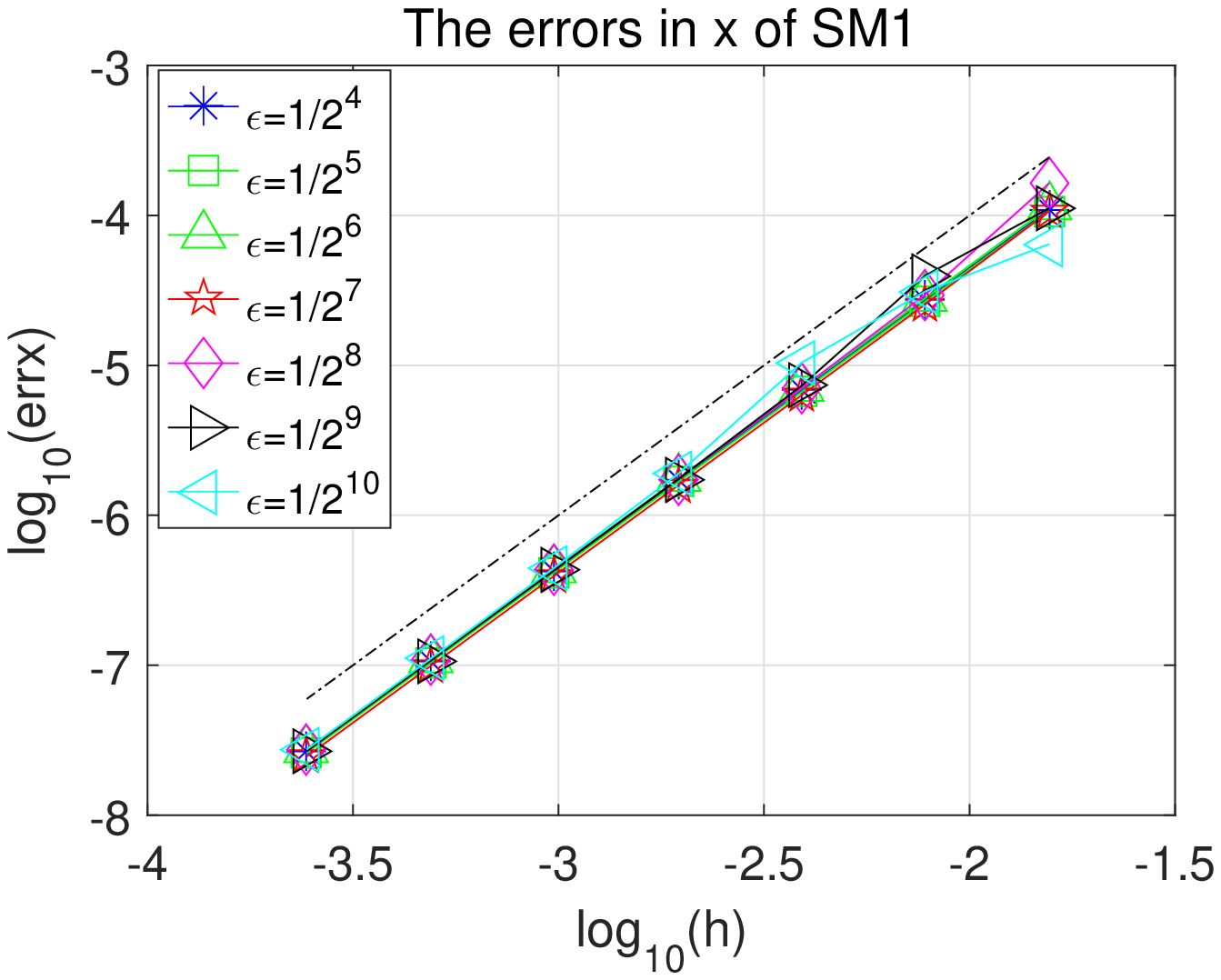}
\includegraphics[width=4.2cm,height=4.2cm]{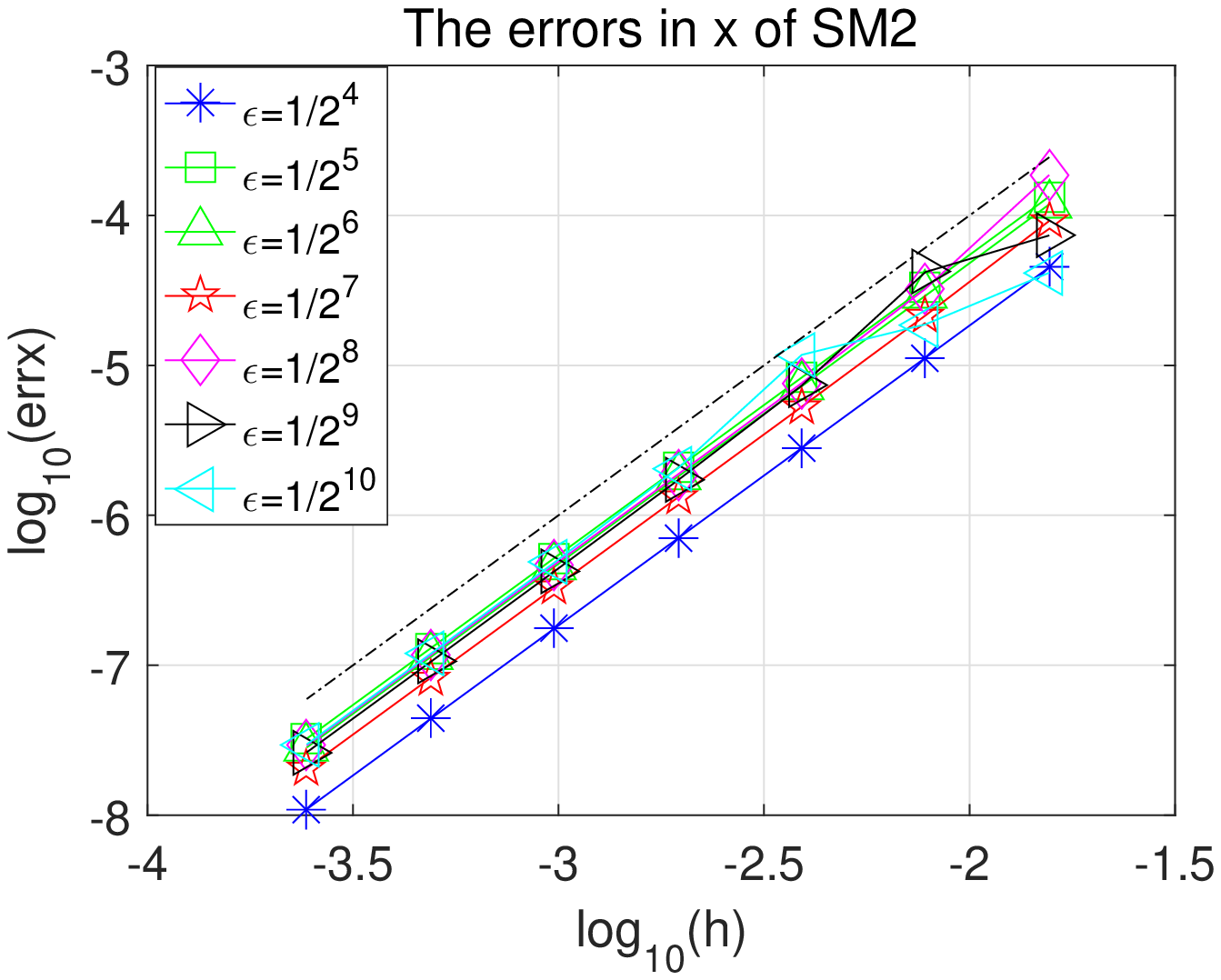}
\includegraphics[width=4.2cm,height=4.2cm]{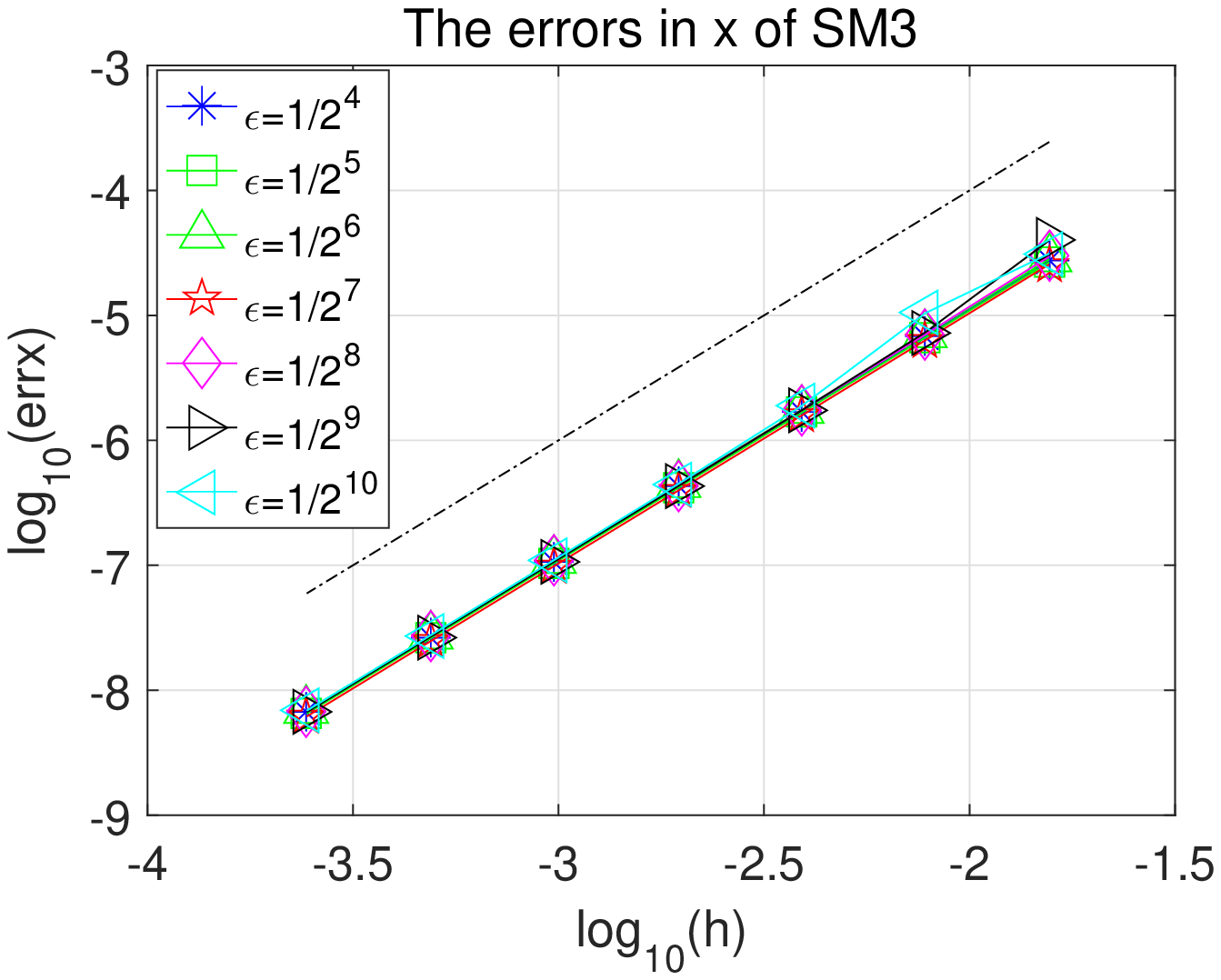}
\includegraphics[width=4.2cm,height=4.2cm]{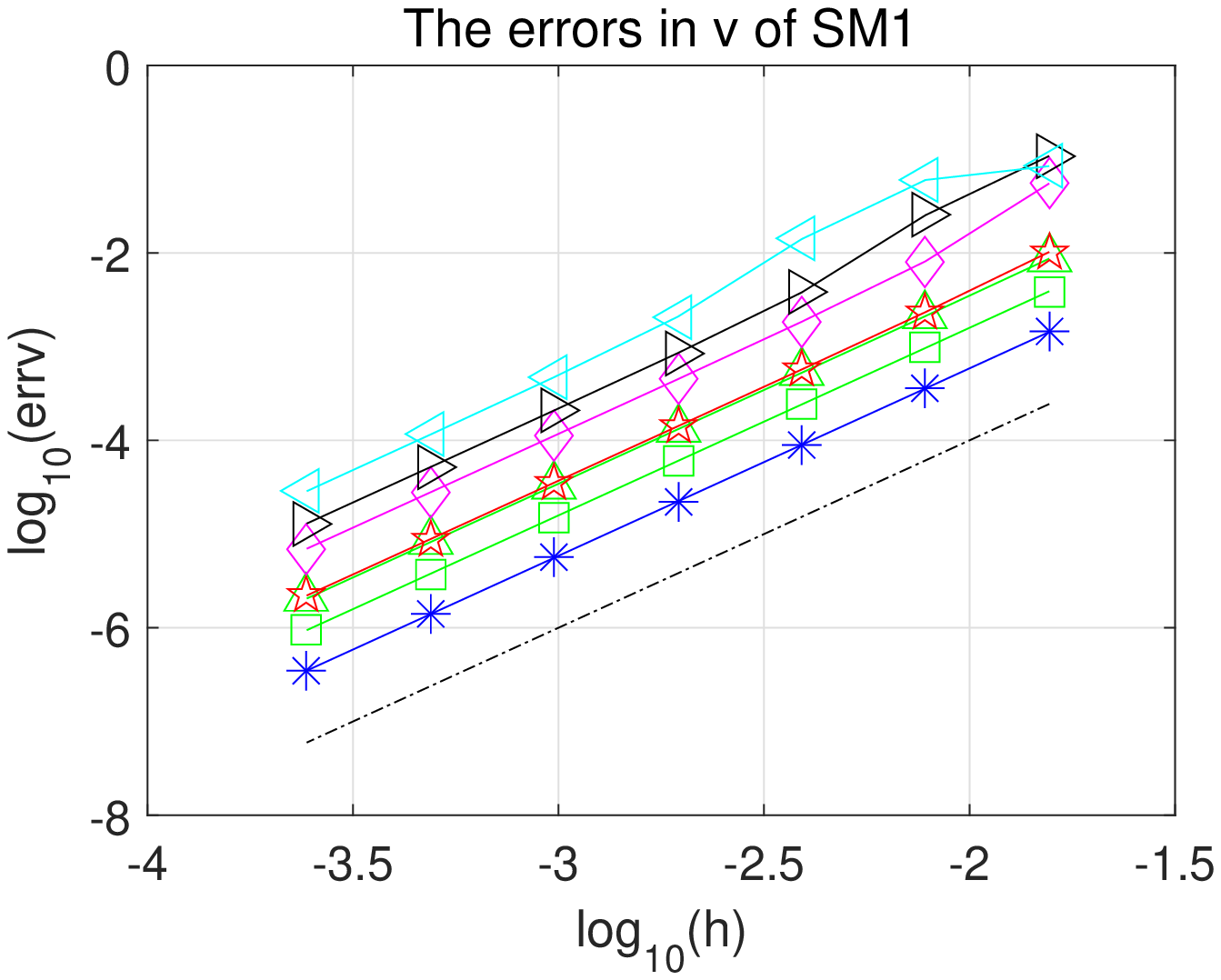}
\includegraphics[width=4.2cm,height=4.2cm]{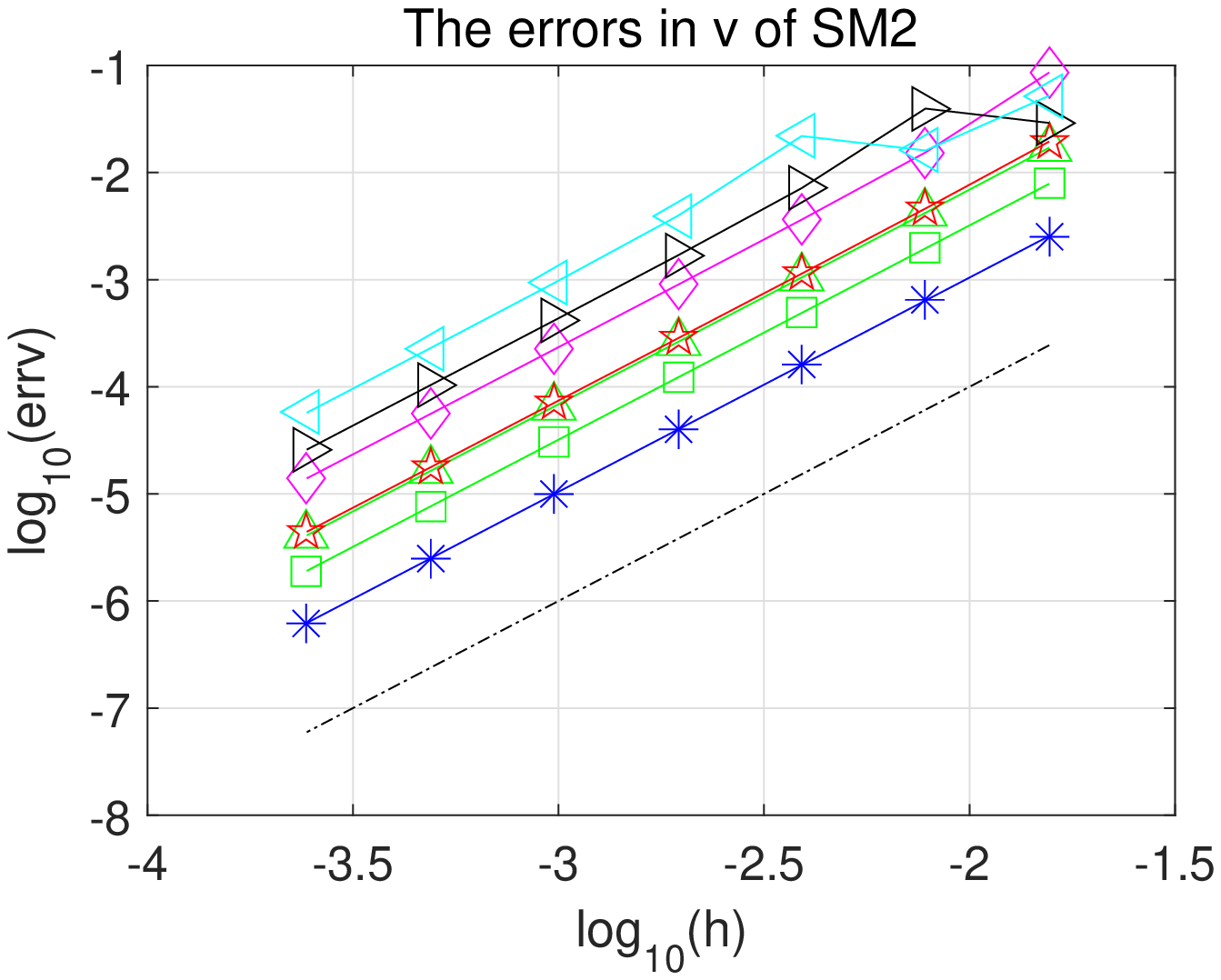}
\includegraphics[width=4.2cm,height=4.2cm]{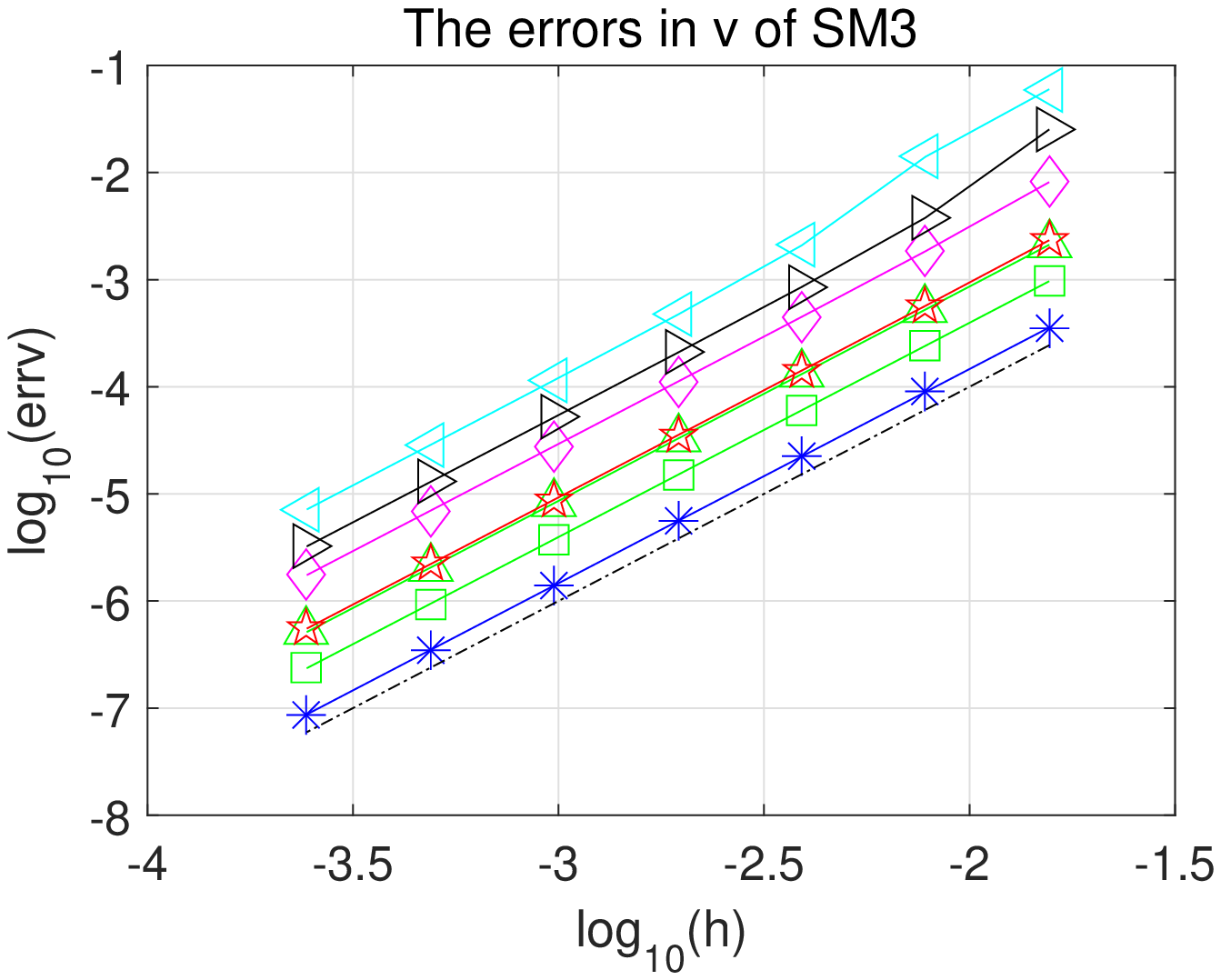}
\caption{The   errors in $x$ (errx) and $v$ (errv) against $h$ for symplectic SM1-SM3  (the slope of the dotted line is two).} \label{p5}
\end{figure}

\begin{figure}[t!]
\centering
\includegraphics[width=4.2cm,height=4.2cm]{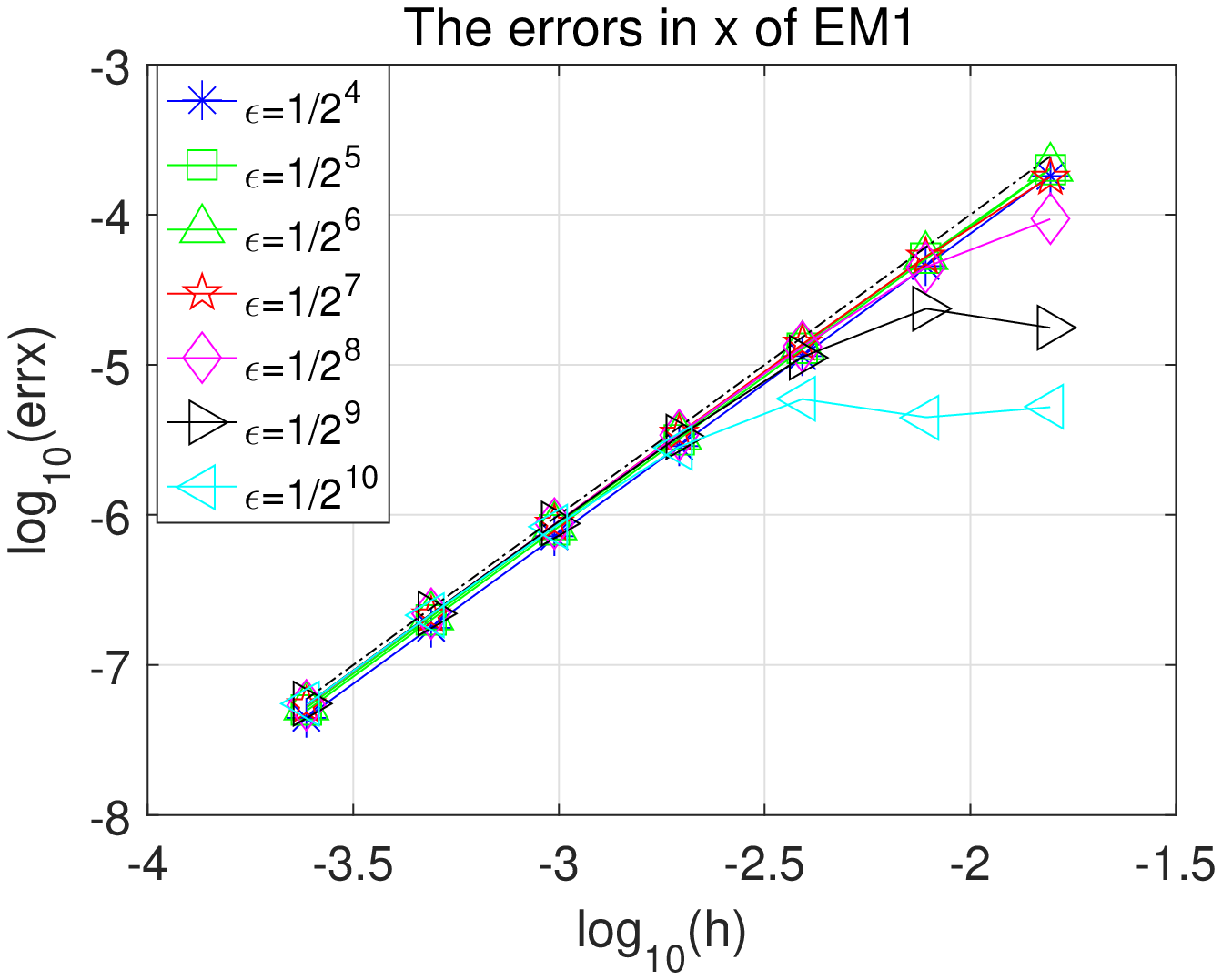}
\includegraphics[width=4.2cm,height=4.2cm]{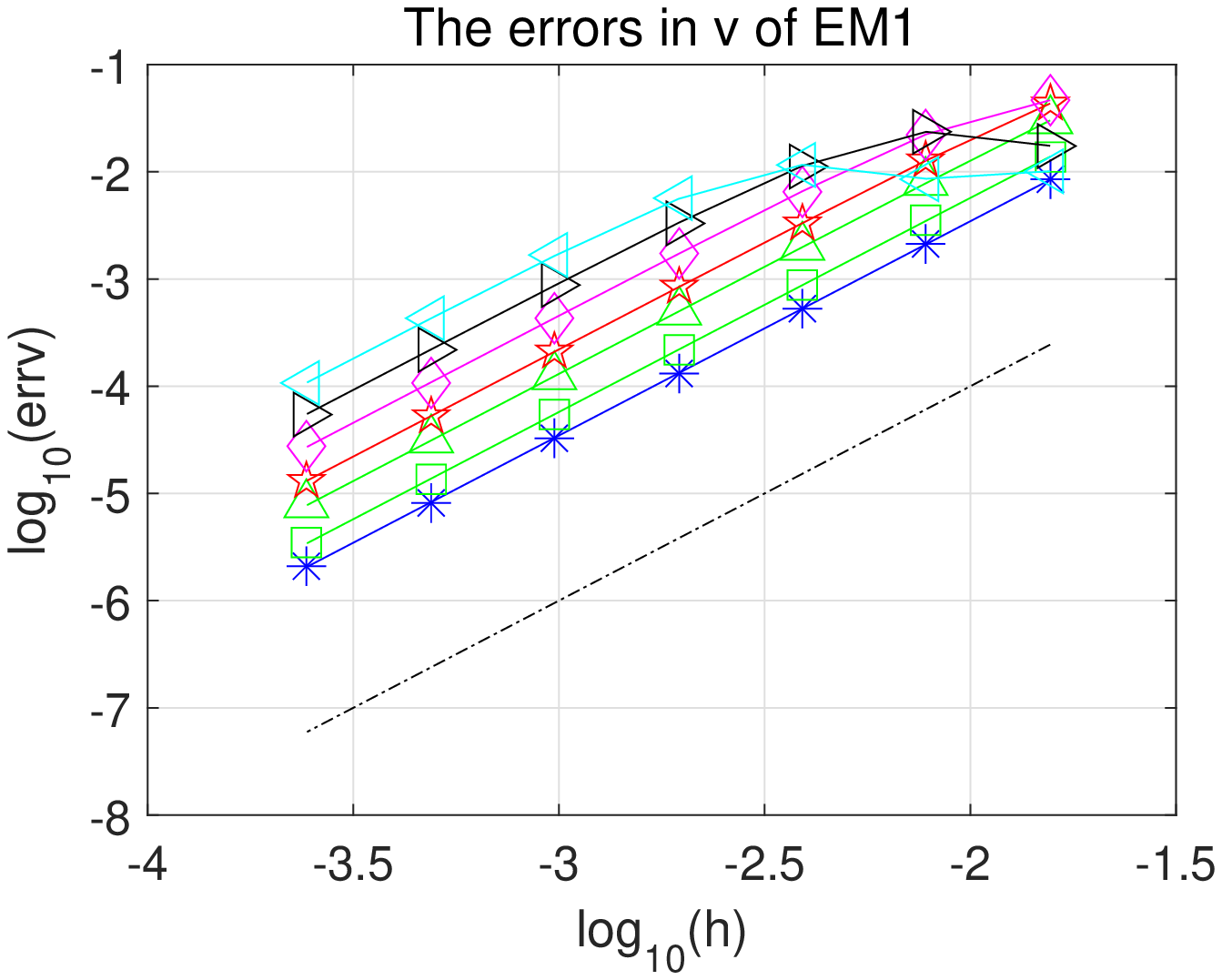}
\caption{The   errors in $x$ (errx) and $v$ (errv) against $h$ for energy-preserving  EM1 (the slope of the dotted line is two).} \label{p6}
\end{figure}

For displaying the results of convergence,
 the problem is solved  on $[0,1]$ with
$h= 1/2^{i}$ for $i=6,\ldots,12$.
 The   global errors
$
   errx:=\frac{\abs{x_n-x(t_n)}}{\abs{x(t_n)}},\ errv:=\frac{\abs{v_n-v(t_n)}}{\abs{v(t_n)}}
$ for different $\varepsilon$ are
shown in Figs. \ref{p4}-\ref{p6}, respectively. It is noted that  we use the result of
standard ODE45 method in MATLAB with an absolute and relative
tolerance equal to $10^{-12}$ as the true solution.  It follows from these results that M1 has   a uniform first-order convergence in both $x$ and $v$ as stated by  \eqref{con M1}.   The other methods only have a uniform second-order error bounds in $x$ as given in Theorem \ref{thm: 4}.

\subsubsection{Efficiency}
\begin{figure}[t!]
\centering
\includegraphics[width=4.2cm,height=4.2cm]{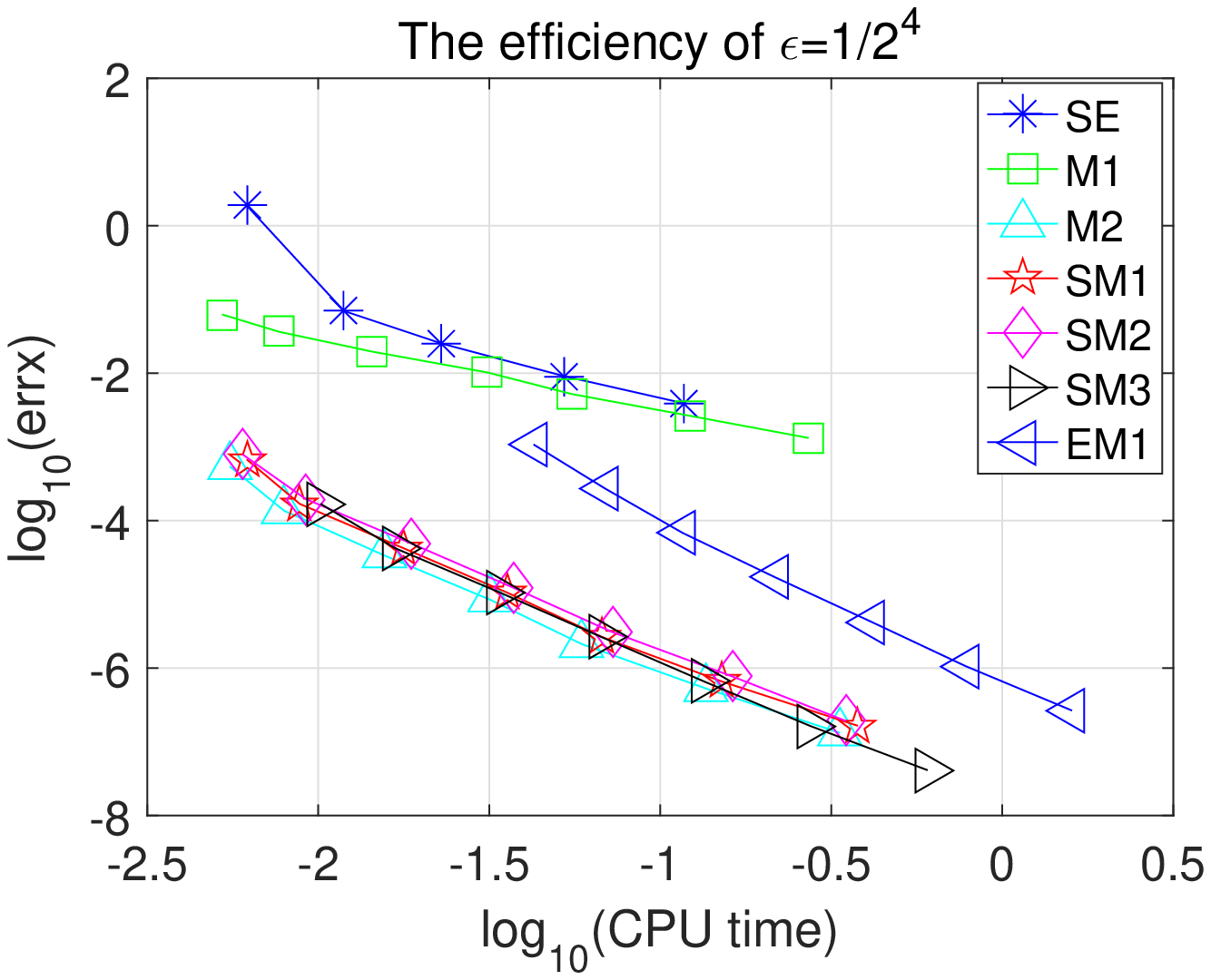}
\includegraphics[width=4.2cm,height=4.2cm]{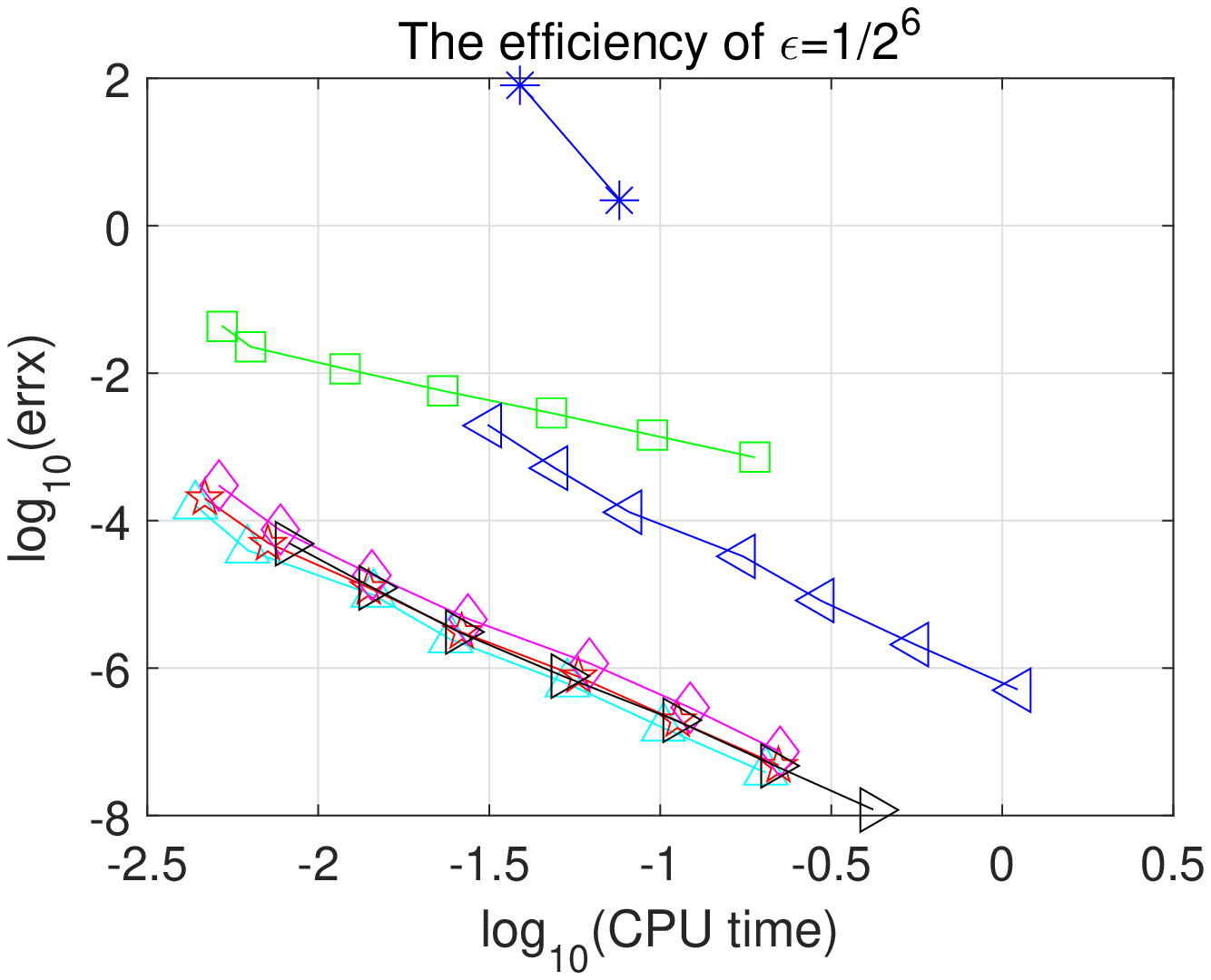}
\includegraphics[width=4.2cm,height=4.2cm]{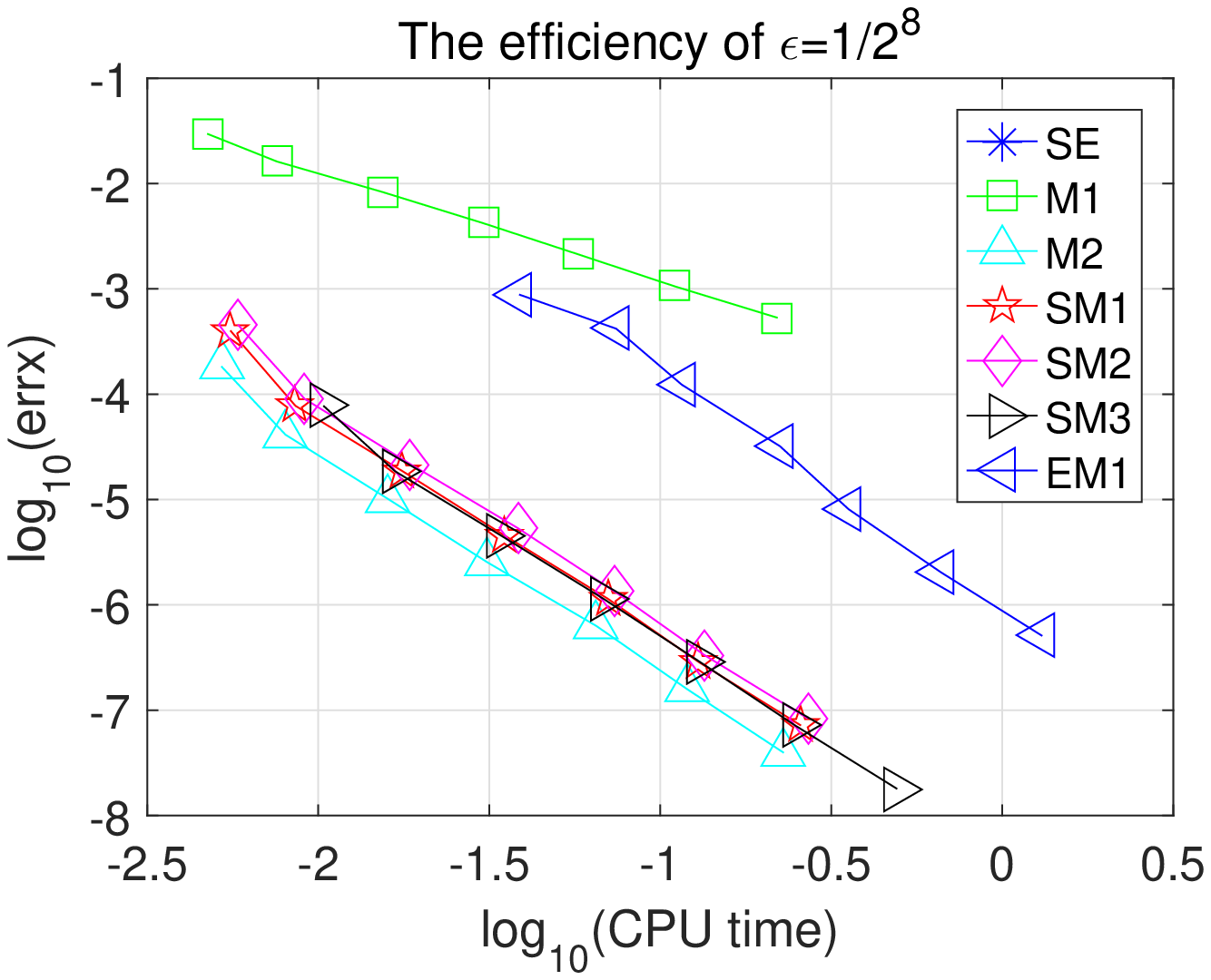}
\caption{The  uniform errors  (err2)  against CPU time.} \label{p7}
\end{figure}

In order to illustrate the efficiency of  the proposed methods, we solve this system till $T=10$. The efficiency of each method (the error err2 versus the CPU time)  is displayed in Fig. \ref{p7}.
This test is conducted  by MATLAB on a laptop
   ThinkPad (CPU: Intel (R) Core (TM) i7-10510U CPU @ 2.30 GHz, Memory: 8 GB, Os: Microsoft Windows 10 with 64bit).
   It can be observed that the computational cost of the new
methods is {cheap} compared with the symplectic Euler
method (Fig.  \ref{p7}).

\begin{figure}[t!]
\centering
\includegraphics[width=4.2cm,height=4.2cm]{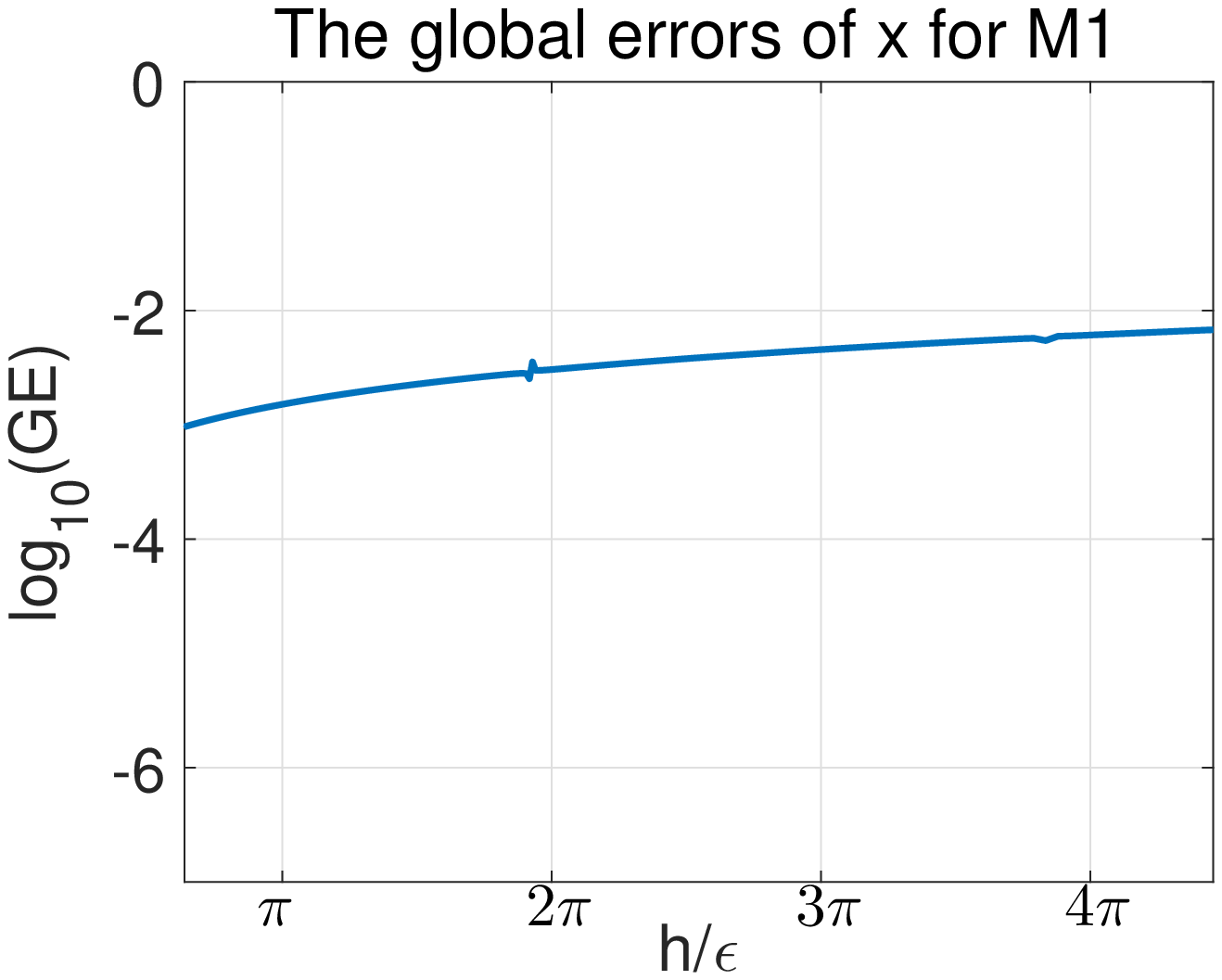}
\includegraphics[width=4.2cm,height=4.2cm]{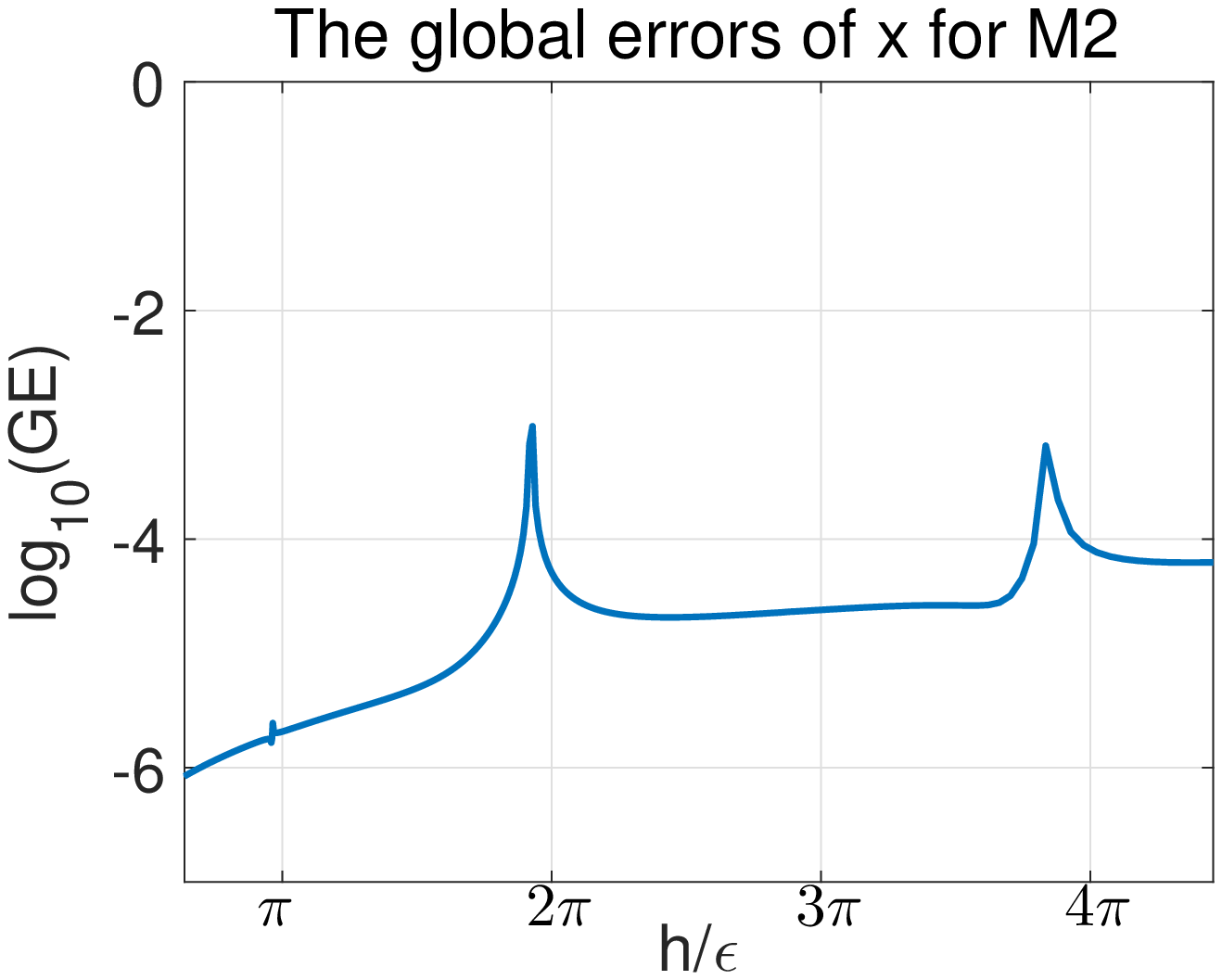}
\includegraphics[width=4.2cm,height=4.2cm]{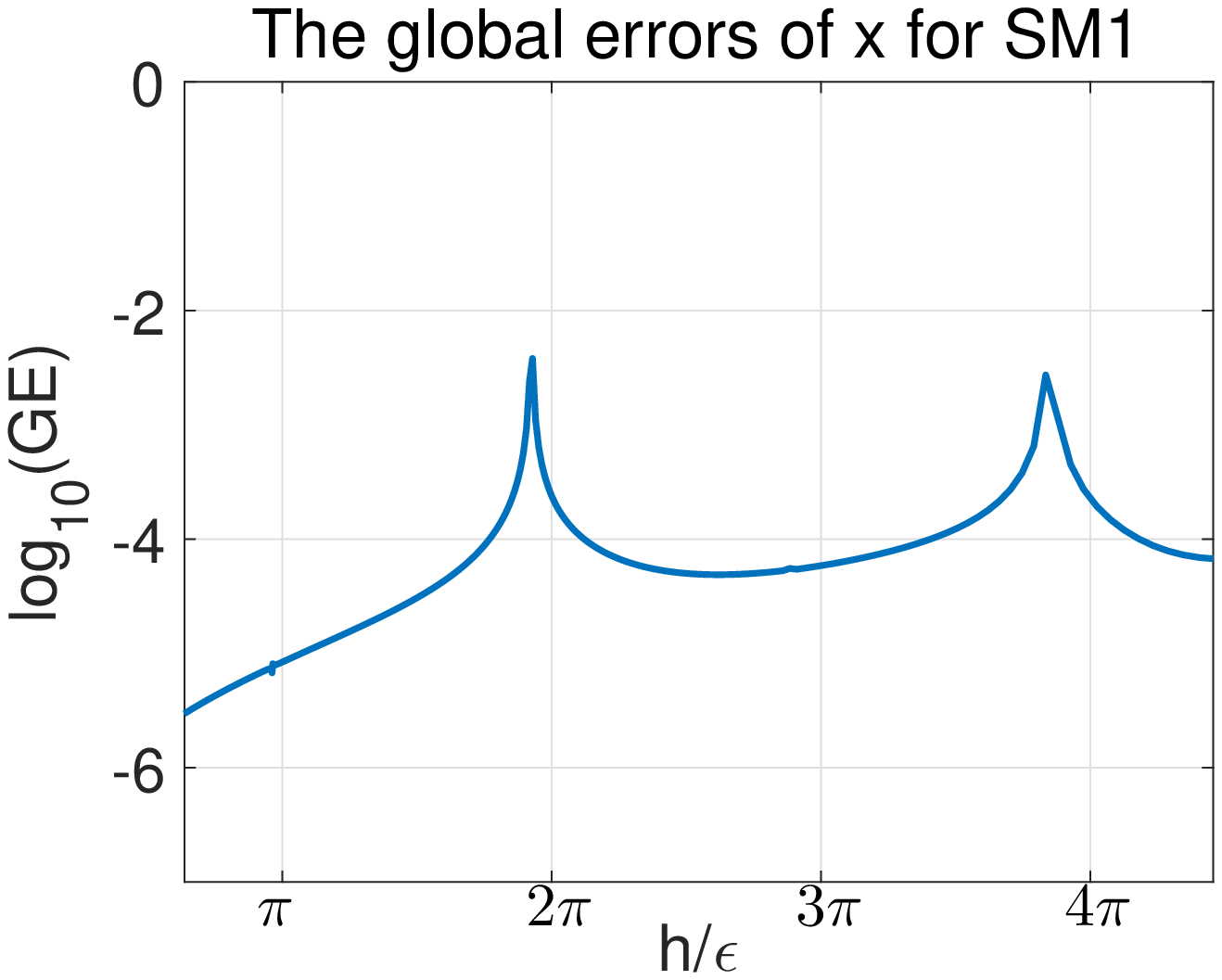}
\includegraphics[width=4.2cm,height=4.2cm]{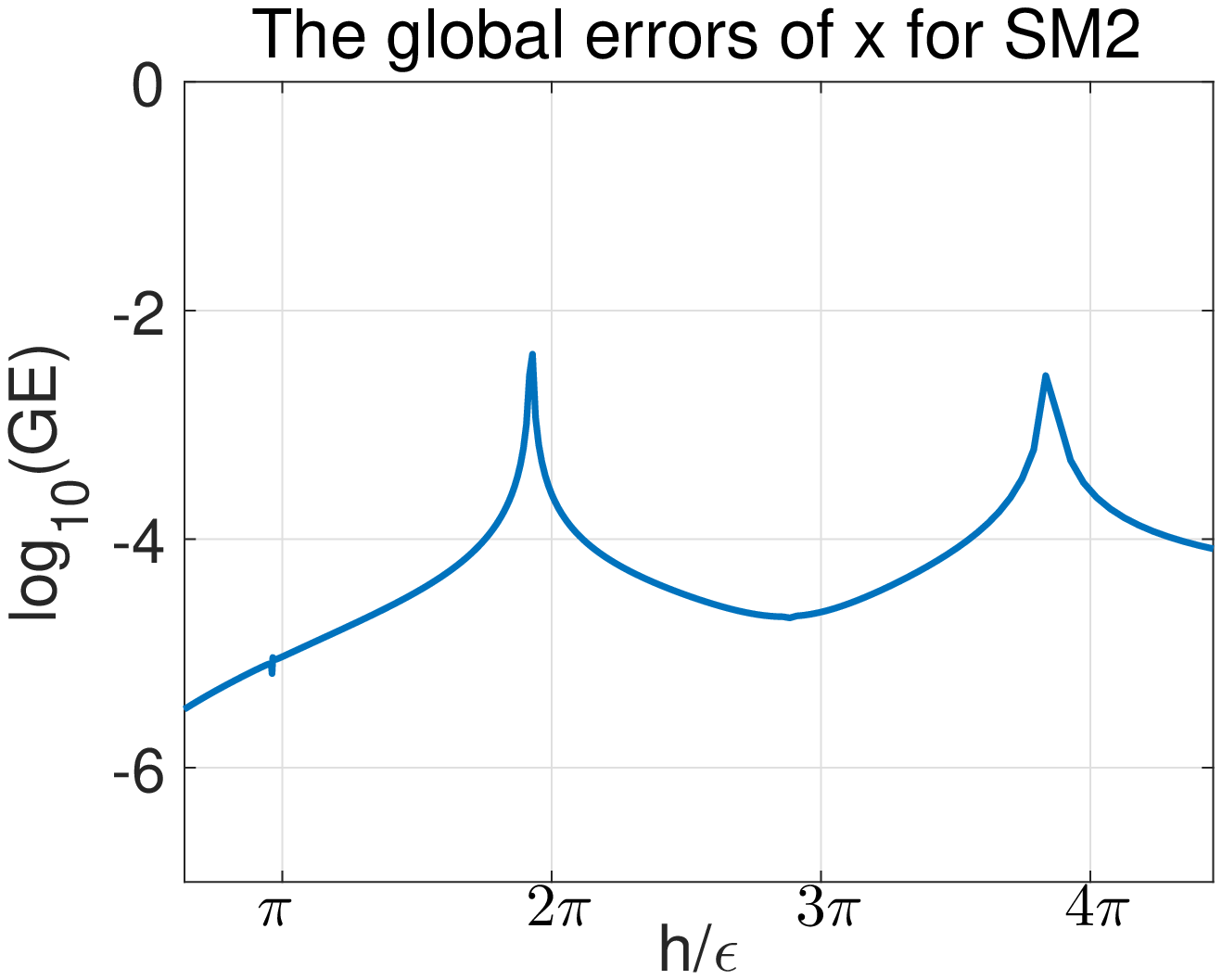}
\includegraphics[width=4.2cm,height=4.2cm]{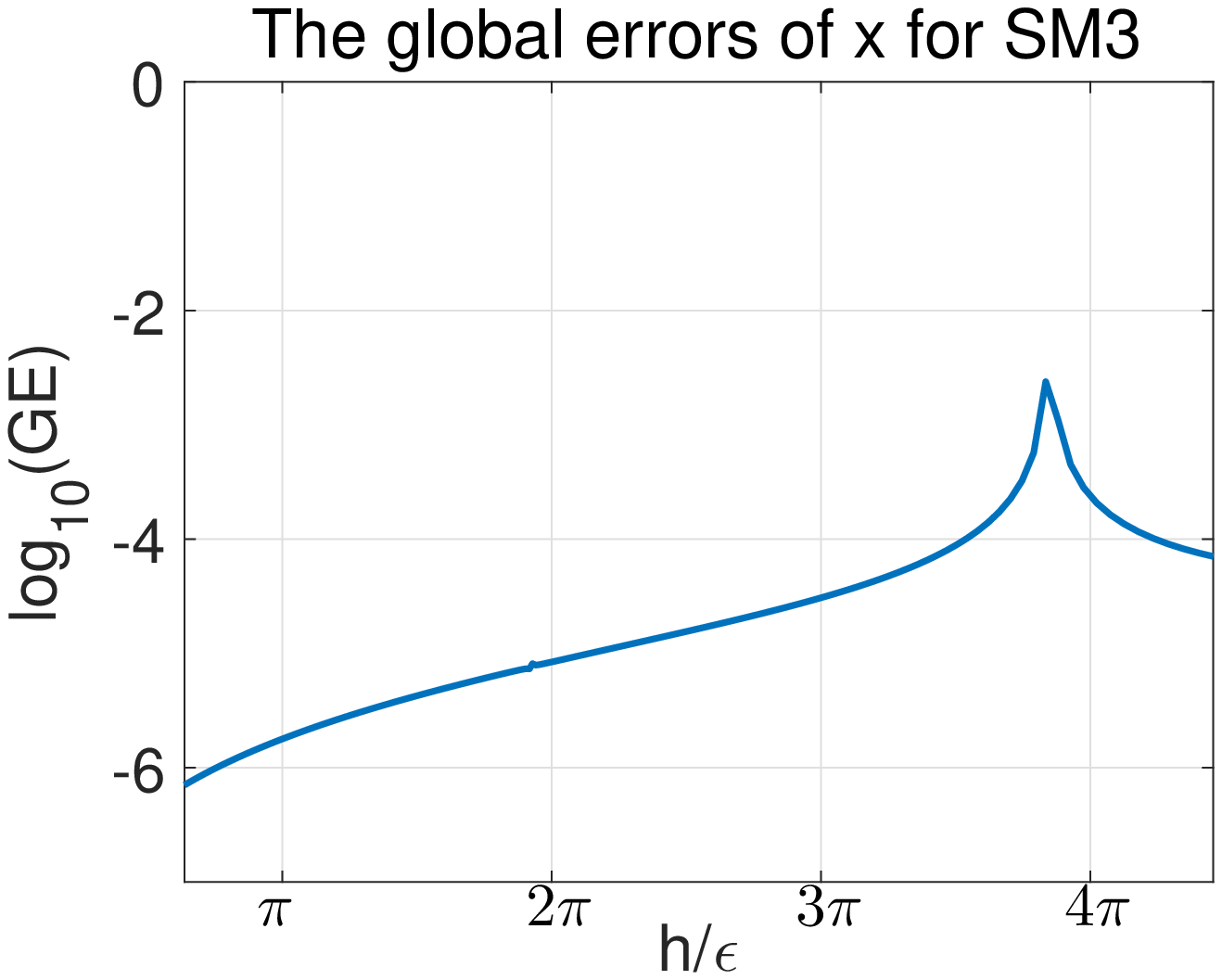}
\includegraphics[width=4.2cm,height=4.2cm]{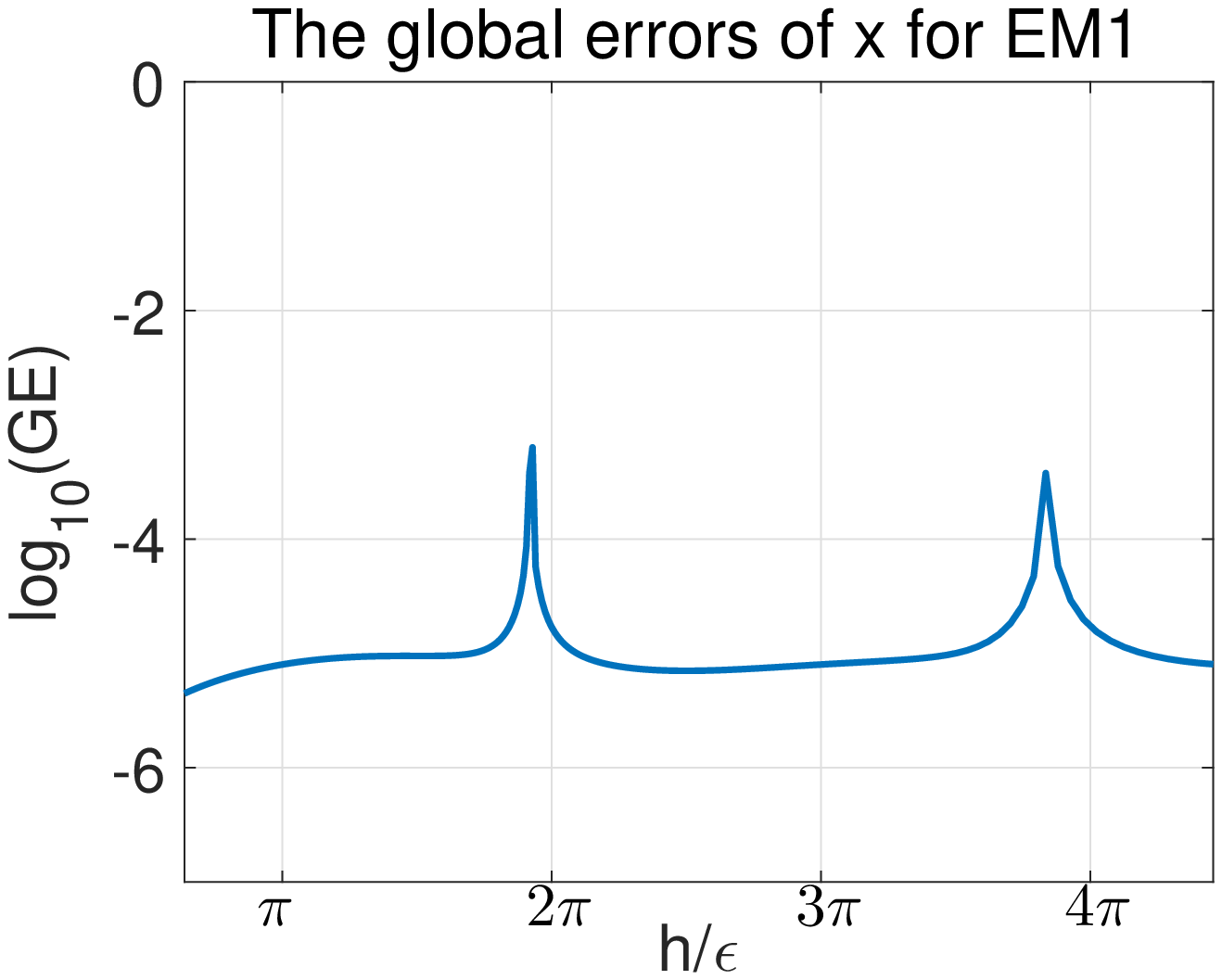}
\caption{The global errors  (GE) of $x$ against $h/\varepsilon$.} \label{p8}
\end{figure}

{\subsubsection{Resonance instability}
 Finally,   we show the resonance instability of the proposed methods. This is done by fixing  $\epsilon=1/2^{10}$ and showing the errors at $T=1$ against $h/\epsilon$  in  Fig.~\ref{p8}.
It can be observed   that M1 gives a very good result but other methods  have a  good behavior for
values of $h/\epsilon$ except integral multiples of $\pi$. SM3 shows a not uniform result close to $4\pi$,   other methods M2, SM1, SM2 and EM1 close to even multiples
of $\pi$. This means that SM3 appears more robust near stepsize resonances and other methods behave very similar away from stepsize resonances.}

In the following three sections, we will prove Theorems \ref{thm: 1}, \ref{thm: 3}-\ref{thm: 4}, respectively. In each proof, we  will firstly consider $d=3$ for brevity and then show that how to extend the analysis   to other  $d$  with some necessary modifications.

\section{Proof of symplecticity (Theorem \ref{thm: 1})}\label{sec: symp methods}
 $\bullet$ \textbf{Transformed system and methods.}

 {Due to}   the  skew-symmetric  matrix
$\tilde{B}$, it is clear that there exists a unitary matrix $P$ and
a   diagonal matrix $\Lambda$ such that $ \tilde{B}=P \Lambda
P^\textup{H}, $ where $\Lambda=
\textmd{diag}(-||\tilde{B}||\mathrm{i},0,||\tilde{B}||\mathrm{i}) $.
 With the linear change of variable
\begin{equation}\label{change of
variable}\tilde{x}(t)= P^\textup{H} x(t),\quad \tilde{v}(t)=
P^\textup{H}
 v(t),\end{equation}
 the system
\eqref{charged-particle sts-cons} can be rewritten as
\begin{equation}\label{necharged-sts-first order}
\begin{array}[c]{ll}
\frac{d}{dt }\left(
  \begin{array}{c}
    \tilde{x} \\
    \tilde{v} \\
  \end{array}
\right)  =\left(
            \begin{array}{cc}
              0 &  I\\
             0  & \tilde{\Omega} \mathrm{i} \\
            \end{array}
          \right) \left(
  \begin{array}{c}
    \tilde{x} \\
    \tilde{v} \\
  \end{array}
\right)+\left(
  \begin{array}{c}
   0 \\
    \tilde{F}(\tilde{x}) \\
  \end{array}
\right),\quad \left(
                \begin{array}{c}
                  \tilde{x}_{0} \\
                  \tilde{v}_{0} \\
                \end{array}
              \right)=\left(
                        \begin{array}{c}
                P^\textup{H} x_{0} \\
                 P^\textup{H} \dot{x}_{0}\\
                        \end{array}
                      \right),
\end{array}
\end{equation}
where $\tilde{\Omega}=
\textmd{diag}(-\tilde{\omega},0,\tilde{\omega})$,
{$\tilde{\omega}=\frac{||\tilde{B}||}{\varepsilon}$,} and
$\tilde{F}(\tilde{x})=P^\textup{H} F(P \tilde{x})=
-\nabla_{\tilde{x}} U(P\tilde{x})$. In this paper, we denote the
vector $x$ by $x=(x^{-1},x^0,x^1)^{\intercal}$ and the same notation
  is used   for all the
vectors in $\RR^3$ or $\CC^3$. According to  \eqref{change of
variable} and the property of the unitary matrix $P$, one has that
\begin{equation}\label{spe cond}\tilde{x}^{-1}=\overline{(\tilde{x}^{1})},\ \
 \tilde{v}^{-1}=\overline{(\tilde{v}^{1})},\ \ \tilde{x}^{0},\
\tilde{v}^{0} \in \RR.\end{equation} The energy of this transformed
system \eqref{necharged-sts-first order} is given by
\begin{equation*}
E(x,v)=\frac{1}{2}\abs{P\tilde{v}}^2+U(P
\tilde{x})=\frac{1}{2}\abs{\tilde{v}}^2+U(P \tilde{x})
:=\tilde{E}(\tilde{x},\tilde{v}).
\end{equation*}
For this transformed  system,  we can modify the schemes of SM1-SM3 accordingly. For example,
 the   scheme   \eqref{AEI-new}  has a transformed form for   \eqref{necharged-sts-first order}
\begin{equation}\label{AEI-trans}
\begin{array}[c]{ll}\tilde{X}_{i}=\tilde{x}_{n}+ c_ih\varphi_1(c_ih\tilde{\Omega}\mathrm{i}) \tilde{v}_{n}+
h^2 \sum\limits_{j=1}^{s}{\alpha}_{ij}(h\tilde{\Omega}\mathrm{i})\tilde{F} (\tilde{X}_j),\ \ i=1,2,\ldots,s,\\
\tilde{x}_{n+1}=\tilde{x}_{n}+ h\varphi_1(h\tilde{\Omega}\mathrm{i})
\tilde{v}_{n}+h^2 \sum\limits_{i=1}^{s}\beta_i(h\tilde{\Omega}\mathrm{i})\tilde{F} (\tilde{X}_i),\\
\tilde{v}_{n+1}=\varphi_0(h\tilde{\Omega}\mathrm{i})\tilde{v}_{n}+h
\sum\limits_{i=1}^{s}\gamma_i(h\tilde{\Omega}\mathrm{i})\tilde{F}
(\tilde{X}_i).
\end{array}
\end{equation}

Denote the transformed method \eqref{AEI-trans} by
\begin{equation}
\begin{array}[c]{ll}\tilde{X}^{J}_{i}=\tilde{x}^{J}_{n}+ c_ih\varphi_1(c_ih\tilde{\Omega}^{J}\mathrm{i}) \tilde{v}^{J}_{n}+
h^2 \sum\limits_{j=1}^{s}{\alpha}_{ij}(h\tilde{\Omega}^{J}\mathrm{i})\tilde{F}^{J}_j,\ \ i=1,2,\ldots,s,\\
\tilde{x}^{J}_{n+1}=\tilde{x}^{J}_{n}+
h\varphi_1(h\tilde{\Omega}^{J}\mathrm{i})
\tilde{v}^{J}_{n}+h^2 \sum\limits_{i=1}^{s}\beta_i(h\tilde{\Omega}^{J}\mathrm{i})\tilde{F}^{J}_i,\\
\tilde{v}^{J}_{n+1}=e^{h\tilde{\Omega}^{J}\mathrm{i}}\tilde{v}^{J}_{n}+h
\sum\limits_{i=1}^{s}\gamma_i(h\tilde{\Omega}^{J}\mathrm{i})\tilde{F}^{J}_i,
\end{array}\label{a}
\end{equation}
where the superscript index $J$ {for $J=-1,0,1$} denotes
the $(J+2)$th entry of a vector or a matrix and $\tilde{F}^{J}_{i}$
denotes  the $(J+2)$th entry of $\tilde{F}(\tilde{X}_{i})$. With the
notation of differential 2-form,  we   need to prove that (see
\cite{hairer2006})
$$
\sum\limits_{J=-1}^{1}dx_{n+1}^{J} \wedge
dp_{n+1}^{J}=\sum_{J=-1}^{1}dx_{n}^{J} \wedge dp_{n}^{J}.$$

We compute
\begin{equation}
\begin{array}[c]{ll}&\sum\limits_{J=-1}^{1}dx_{n+1}^{J} \wedge
dp_{n+1}^{J}=\sum\limits_{J=-1}^{1}d\bar{x}_{n+1}^{J} \wedge
dp_{n+1}^{J}=\sum\limits_{J=-1}^{1}d(\bar{P}\bar{\tilde{x}}_{n+1})^{J}
\wedge
d(P\tilde{p}_{n+1})^{J}\\
&= \sum\limits_{J=-1}^{1}\Big( \sum\limits_{i=-1}^{1}(\bar{P}_{J+2,i+2}{
d}\bar{\tilde{x}}_{n+1}^{i})\Big)\wedge
 \Big(  \sum\limits_{k=-1}^{1}(P_{J+2,k+2}{ d}\tilde{p}_{n+1}^{k})\Big)\\
 &= \sum\limits_{J=-1}^{1}\sum\limits_{i=-1}^{1}\sum\limits_{k=-1}^{1} \bar{P}_{J+2,i+2}  P_{J+2,k+2} ({
d}\bar{\tilde{x}}_{n+1}^{i} \wedge  { d}\tilde{p}_{n+1}^{k})\\
 &=  \sum\limits_{i=-1}^{1}  {
d}\bar{\tilde{x}}_{n+1}^{i} \wedge  { d}\tilde{p}_{n+1}^{i}=
\sum\limits_{J=-1}^{1}  { d}\bar{\tilde{x}}_{n+1}^{J} \wedge  {
d}\tilde{p}_{n+1}^{J},\nonumber
\end{array}\label{td}
\end{equation}
where $P^HP=I$ is used here.
Similarly, one has $ \sum\limits_{J=-1}^{1}dx_{n}^{J}
\wedge dp_{n}^{J}= \sum\limits_{J=-1}^{1}  {
d}\bar{\tilde{x}}_{n}^{J} \wedge { d}\tilde{p}_{n}^{J}.$ Thus we
only need to prove
$ \sum\limits_{J=-1}^{1}  { d}\bar{\tilde{x}}_{n+1}^{J} \wedge  {
d}\tilde{p}_{n+1}^{J}= \sum\limits_{J=-1}^{1}  {
d}\bar{\tilde{x}}_{n}^{J} \wedge { d}\tilde{p}_{n}^{J},
$
i.e.
\begin{equation}\label{imp res}
\begin{array}[c]{ll} \sum\limits_{J=-1}^{1}  { d}\bar{\tilde{x}}_{n+1}^{J} \wedge  {
d}\tilde{v}_{n+1}^{J}-\frac{1}{2} \sum\limits_{J=-1}^{1}  {
d}\bar{\tilde{x}}_{n+1}^{J} \wedge  { d}
(\tilde{\Omega}^{J}\mathrm{i}\tilde{x}_{n+1}^{J})\\
=
\sum\limits_{J=-1}^{1}  { d}\bar{\tilde{x}}_{n}^{J} \wedge {
d}\tilde{v}_{n}^{J}-\frac{1}{2}  \sum\limits_{J=-1}^{1}  {
d}\bar{\tilde{x}}_{n}^{J} \wedge { d}
(\tilde{\Omega}^{J}\mathrm{i}\tilde{x}_{n}^{J}).
\end{array}
\end{equation}

$\bullet$  \textbf{Symplecticity of the transformed  methods.}

In this part, we will prove that the result \eqref{imp res} is true if the
following conditions are satisfied
\begin{equation}
\begin{aligned}&\gamma_{j}(K) -K\beta_{j}(K)
=d_j I, \ \ d_j\in \mathbb{C},\\
&\gamma_{j}(K) [ \bar{\varphi}_{1}(K) -c_{j}\bar{\varphi}_1(c_jK)]
=\beta_{j}(K) [ e^{-K}+K
\bar{\varphi}_{1}(K)-c_jK \bar{\varphi}_1(c_jK)],
\\
&\bar{\beta}_{i}(K) \gamma_{j}(K) -\frac{1}{2}K \bar{\beta}_{i}(K)
\beta_{j}(K) -\bar{{\alpha}}_{ji}(K) [ \gamma_{j}(K)
  -K\beta_{j}(K)   ]\\
 &\qquad =
\beta_{j}(K)  \bar{ \gamma}_{i}(K) +\frac{1}{2}K \beta_{j}(K)
\bar{\beta}_{i}(K) -{\alpha}_{ij}(K) [ \bar{\gamma}_{i}(K)
  +K\bar{\beta}_{i}(K)   ],
\end{aligned}\label{17}
\end{equation}
where $i,j=1,2,\ldots,s,$ and $K=h\tilde{\Omega}\mathrm{i}$.
Here $\bar{\varphi}_1$ denotes the conjugate of
$\varphi_1$ and the same notation is used for other functions.

{In view of} the definition of differential 2-form (see
\cite{hairer2006}), it can be proved that $ \overline{{
d}\bar{\tilde{x}}_{n}^{J} \wedge { d}  \tilde{v}_{n}^{J}} ={ d}
\tilde{x}_{n}^{J} \wedge { d}  \bar{\tilde{v}}_{n}^{J}\ \textmd{and}
\  { d}\bar{\tilde{x}}_{n}^{J} \wedge { d} \tilde{x}_{n}^{J} \in
\mathrm{i} \mathbb{R}.$ In the light of the scheme \eqref{a} and the fact that   any exterior product
$\wedge$ appearing here is real, it is
obtained that
\begin{equation}
\begin{aligned}
   &d\bar{\tilde{x}}_{n+1}^{J}\wedge d\tilde{v}_{n+1}^{J}-\frac{1}{2}d\bar{\tilde{x}}_{n+1}^{J}\wedge d(\tilde{\Omega}^{J}\mathrm{i}\tilde{x}_{n+1}^{J})
=   d\bar{\tilde{x}}_{n}^{J}  \wedge d\tilde{v}_{n}^{J}
 - \frac{1}{2}d\bar{\tilde{x}}_{n}^{J}\wedge  d(\tilde{\Omega}^{J}\mathrm{i}\tilde{x}_{n}^{J})\\
 + &h\sum\limits_{j=1}^{s}[  \gamma_{j}(K^{J})
-K^{J} \beta_{j}(K^{J})            ]
d\bar{\tilde{x}}_{n}^{J} \wedge  d\tilde{F}_{j}^{J} +[he^{K^{J}} \bar{\varphi}_{1}(K^{J})
-\frac{1}{2}h^{2} \tilde{\Omega}^{J}\mathrm{i}\bar{\varphi}_{1}(K^{J}) \varphi_{1}(K^{J})]
d\bar{\tilde{v}}_{n}^{J}  \wedge   d\tilde{v}_{n}^{J}\\
+&h^{2}\sum\limits_{j=1}^{s}[\bar{\varphi}_{1}(K^{J})   \gamma_{j}(K^{J})
- \beta_{j}(K^{J})  e^{-K^{J}}
 -h \tilde{\Omega}^{J}\mathrm{i} \bar{\varphi}_{1}(K^{J}) \beta_{j}(K^{J})]
d\bar{\tilde{v}}_{n}^{J} \wedge  d\tilde{F}_{j}^{J}\\
+&h^{3}\sum\limits_{i,j=1}^{s}[\bar{\beta}_{i}(K^{J})   \gamma_{j}(K^{J})
-\frac{1}{2}h \tilde{\Omega}^{J}\mathrm{i} \bar{\beta}_{i}(K^{J})   \beta_{j}(K^{J})]
d\bar{\tilde{F}}_{i}^{J} \wedge d\tilde{F}_{j}^{J},
\end{aligned}\label{21}
\end{equation} where the fact that $e^{K^{J}} -h
\tilde{\Omega}^{J}\mathrm{i}\varphi_{1}(K^{J})=I$ is used here.

On the other hand, from the {first $s$ equalities} of
\eqref{a}, it follows that
\begin{equation}
d\tilde{x}^{J}_{n} =d\tilde{X}^{J}_{i}- c_ih\varphi_1(c_iK^{J})d\tilde{v}^{J}_{n}-h^2
\sum\limits_{j=1}^{s} \alpha_{ij}(K^{J})d\tilde{F}_{j}^{J},\ \
i=1,2,\ldots,s. \nonumber
\end{equation}
{We then obtain}
\begin{equation}
d\bar{\tilde{x}}^{J}_{n}\wedge d\tilde{F}^{J}_{j}
=d\bar{\tilde{X}}^{J}_{j}\wedge d\tilde{F}^{J}_{j} -
c_jh\bar{\varphi}_1(c_jK^{J})d\bar{\tilde{v}}^{J}_{n}
\wedge d\tilde{F}^{J}_{j}-h^2
\sum\limits_{i=1}^{s}\bar{{\alpha}}_{ji}(K^{J})
d\bar{\tilde{F}}_{i}^{J}\wedge d\tilde{F}^{J}_{j},\ \
j=1,2,\ldots,s. \nonumber
\end{equation}
Inserting this into  \eqref{21} and summing over all $J$ yields
\begin{subequations}
\begin{align}
\sum\limits_{J=-1}^{1}d\bar{\tilde{x}}_{n+1}^{J}\wedge d\tilde{v}_{n+1}^{J}&-\frac{1}{2}\sum\limits_{J=-1}^{1}d\bar{\tilde{x}}_{n+1}^{J}\wedge d(\tilde{\Omega}^{J}\mathrm{i}\tilde{x}_{n+1}^{J})=\sum\limits_{J=-1}^{1} d\bar{\tilde{x}}_{n}^{J}  \wedge
d\tilde{v}_{n}^{J}
- \frac{1}{2}\sum\limits_{J=-1}^{1} d\bar{\tilde{x}}_{n}^{J}\wedge  d(\tilde{\Omega}^{J}\mathrm{i}\tilde{x}_{n}^{J})  \nonumber \\
   + &  h\sum\limits_{j=1}^{s}\sum\limits_{J=-1}^{1} [ \gamma_{j}(K^{J})
-K^{J}\beta_{j}(K^{J})]
d\bar{\tilde{X}}^{J}_{j}\wedge d\tilde{F}^{J}_{j} \label{sta er 0}  \\
  + & h \sum\limits_{J=-1}^{1}  [e^{K^{J}} \bar{\varphi}_{1}(K^{J})
-\frac{1}{2}K^{J}\bar{\varphi}_{1}(K^{J}) \varphi_{1}(K^{J})]
d\bar{\tilde{v}}_{n}^{J}  \wedge   d\tilde{v}_{n}^{J} \label{sta er 1}\\
    + &
h^{2}\sum\limits_{j=1}^{s}\sum\limits_{J=-1}^{1} \Big[\bar{\varphi}_{1}(K^{J})    \gamma_{j}(K^{J})
- \beta_{j}(K^{J})  e^{-K^{J}} -K^{J}  \bar{\varphi}_{1}(K^{J}) \beta_{j}(K^{J})
\label{sta er 2}\\
 &\qquad \qquad \ \ -c_j\bar{\varphi}_1(c_jK^{J})
 [\gamma_{j}(K^{J})
  -K^{J}\beta_{j}(K^{J})]
 \Big]
d\bar{\tilde{v}}_{n}^{J} \wedge  d\tilde{F}_{j}^{J} \label{sta er 3}\\
    + &h^3\sum\limits_{i,j=1}^{s}\sum\limits_{J=-1}^{1}\Big[
\bar{\beta}_{i}(K^{J})   \gamma_{j}(K^{J})
-\frac{1}{2}h \tilde{\Omega}^{J}\mathrm{i} \bar{\beta}_{i}(K^{J})   \beta_{j}(K^{J})\label{sta er 4}\\
   & \qquad \qquad \ \
-\bar{{\alpha}}_{ji}(K^{J})
[  \gamma_{j}(K^{J})
  -K^{J} \beta_{j}(K^{J})   ]
\Big]d\bar{\tilde{F}}_{i}^{J} \wedge d\tilde{F}_{j}^{J}. \label{sta er 5}
 \end{align}
\end{subequations}

$\circ$   Prove that $\eqref{sta er 0}=0$.

Based on the {first $s$ conditions} of \eqref{17},
$\tilde{F}(\tilde{x})=-\nabla_{\tilde{x}} U(P\tilde{x})$ and
\eqref{td}, it can be {verified} that
$d\bar{\tilde{X}}^{J}_{j}\wedge d\tilde{F}^{J}_{j}= d
X^{J}_{j}\wedge dF^{J}_{j}$. Thus, one has
\begin{equation}
\begin{array}
[c]{llllllllll}
&\sum\limits_{J=-1}^{1} [ \gamma_{j}(K^{J})
-K^{J}\beta_{j}(K^{J})            ]
d\bar{\tilde{X}}^{J}_{j}\wedge d\tilde{F}^{J}_{j}\\
=&d_{j}\sum\limits_{J=-1}^{1}   d\bar{\tilde{X}}^{J}_{j}\wedge
d\tilde{F}^{J}_{j} =d_{j}\sum\limits_{J=-1}^{1}   d X^{J}_{j}\wedge
dF^{J}_{j}
=-d_{j}\sum\limits_{J=-1}^{1}  dF^{J}_{j} \wedge d X^{J}_{j}\\
= &-d_{j} \sum\limits_{J=-1}^{1} (\frac{\partial{
F_{j}^{J}(X_{j})}}{\partial{x^{I}}} d X_{j}^{I})  \wedge d X_{j}^{J}
=-d_{j} \sum\limits_{J,I=-1}^{1}(-\frac{\partial^{2}{ U(Px)} }
{\partial{x^{J}} \partial{x^{I}}} )d X_{j}^{I}  \wedge d
X_{j}^{J}=0.
\end{array}\nonumber
\end{equation}

 $\circ$    Prove that $\eqref{sta er 1}=0$.

{Using} the  property  of $\tilde{v}_{n}$, we have
\begin{equation}
 d\bar{\tilde{v}}_{n}^{-1}  \wedge   d\tilde{v}_{n}^{-1}=-d\bar{\tilde{v}}_{n}^{1}  \wedge   d\tilde{v}_{n}^{1}, \  \ \ d\bar{\tilde{v}}_{n}^{0}  \wedge   d\tilde{v}_{n}^{0}=0, \nonumber
\end{equation}
and
\begin{equation}
 e^{K^{1}} \bar{\varphi}_{1}(K^{1})
-\frac{1}{2}K^{1} \bar{\varphi}_{1}(K^{1})
\varphi_{1}(K^{1})
=e^{K^{-1}}
\bar{\varphi}_{1}(K^{-1}) -\frac{1}{2}K^{-1} \bar{\varphi}_{1}(K^{-1})
\varphi_{1}(K^{-1}).\nonumber
 \end{equation}
Therefore, it {follows} that
  \begin{equation}
  \sum\limits_{J=-1}^{1}  [e^{K^{J}} \bar{\varphi}_{1}(K^{J})
-\frac{1}{2}K^{J}\bar{\varphi}_{1}(K^{J}) \varphi_{1}(K^{J})]
d\bar{\tilde{v}}_{n}^{J}  \wedge   d\tilde{v}_{n}^{J}=0.\nonumber
\end{equation}

$\circ$    Prove that \eqref{sta er 2}-\eqref{sta er 5}$=0$.

In the light of  {all the identities after the previous
$s$ ones in}  \eqref{17}, the last two terms \eqref{sta er
2}-\eqref{sta er 5} vanish.

{The  results stated above leads to \eqref{imp res}.}
Then, it can be {verified straightforwardly}  that the
coefficients of SM1-SM3 satisfy \eqref{17}. Therefore, these methods
are symplectic.

$\bullet$  \textbf{Extension of the proof to other $d$.}

For a general $d\geq2,$ since $\tilde{B}$ is  skew-symmetric,  there exists a unitary matrix $P$ and
a   diagonal matrix $\Lambda$ such that $ \tilde{B}=P \Lambda
P^\textup{H}, $ where
\begin{equation} \label{lambda}\Lambda=
\left\{\begin{aligned}
&\textmd{diag}(-\tilde{\omega}_l\mathrm{i},\ldots,-\tilde{\omega}_1\mathrm{i},0,\tilde{\omega}_1\mathrm{i},\ldots,\tilde{\omega}_l\mathrm{i}),\ \ d=2l+1,\\
&\textmd{diag}(-\tilde{\omega}_l\mathrm{i},\ldots,-\tilde{\omega}_1\mathrm{i},\tilde{\omega}_1\mathrm{i},\ldots,\tilde{\omega}_l\mathrm{i}),\ \ \ \ \  d=2l.
\end{aligned}\right.
\end{equation}
For both cases, the above proof can be extended without any difficulty.

 \section{Proof of long-time energy conservation (Theorem \ref{thm: 3} )}\label{sec: conservations}
In this section, we will show the long time near-conservation of
energy along SM2 algorithm. We first  derive
 modulated Fourier expansion (see, e.g. \cite{Hairer00,Hairer2018,lubich19,Wang2021})  with
sufficient many terms for SM2. Then one almost-invariant of the
expansion is studied and  based on which the long-time near
conservation is confirmed.  The proof of other methods can be given
by modifying the operators $\mathcal{L}(hD),\hat{\mathcal{L}}(hD)$ \eqref{new L1L2}
and following the  way given {below.}

$\bullet$ \textbf{Reformulation of SM2.}

Using  symmetry, the  algorithm  SM2 can be expressed in a two-step
form
\begin{equation}\label{2 steps method-0}
\left\{ \begin{array}[c]{ll}
 x_{n+1}-2x_{n}+x_{n-1}=
h(\varphi_1(h \Omega)-\varphi_1(-h \Omega)) v_{n}+\frac{1}{2}h^2(\varphi_1(h \Omega)+\varphi_1(-h \Omega)) F_{n},\\
 x_{n+1}-x_{n-1}=
h(\varphi_1(h \Omega)+\varphi_1(-h \Omega))
v_{n}+\frac{1}{2}h^2(\varphi_1(h \Omega)-\varphi_1(-h \Omega)) F_{n},
\end{array}\right.
\end{equation}
with $F_{n}:=F(x_{n}),$
which yields that
\begin{equation*}
 \begin{array}[c]{ll}
 \alpha (h \Omega)\frac{x_{n+1}-2x_{n}+x_{n-1}}{h^2}= \beta(h \Omega) \Omega \frac{x_{n+1}-x_{n-1}}{2h}  +\gamma
 (h \Omega)  F_{n},
\end{array}
\end{equation*}
where $
 \alpha (\xi)=\frac{\xi}{\varphi_1(\xi)-\varphi_1(-\xi)},\  \
\beta (\xi)=\frac{2}{\varphi_1(\xi)+\varphi_1(-\xi)},\  \ \gamma
(\xi)=\xi\frac{2
\varphi_1(\xi)\varphi_1(-\xi)}{\varphi^2_1(\xi)-\varphi^2_1(-\xi)}.$

 For the transformed
system \eqref{necharged-sts-first order}, it becomes
 \begin{equation}\label{2 steps method-new}
 \begin{array}[c]{ll}
\tilde {\alpha}
(h\tilde{\Omega})\frac{\tilde{x}_{n+1}-2\tilde{x}_{n}+\tilde{x}_{n-1}}{h^2}=
\tilde{\beta}(h\tilde{\Omega})\mathrm{i}
 \tilde{\Omega} \frac{\tilde{x}_{n+1}-\tilde{x}_{n-1}}{2h}  +\tilde{\gamma}
 (h\tilde{\Omega})  \tilde{F}_{n},
\end{array}
\end{equation}
where the coefficient functions are given by $
 \tilde{\alpha} (\xi)=\frac{1}{\sinc^2(\frac{\xi}{2})},\
\tilde{\beta} (\xi)=\frac{1}{\sinc(\xi)},\ \tilde{\gamma}
(\xi)=\xi\csc(\xi)$
with $\sinc(\xi)=\sin(\xi)/\xi$.

Define the operators
\begin{equation}\label{new L1L2}
 \begin{array}[c]{ll}
\mathcal{L}(hD)= \frac{1}{2h\sinc(h\tilde{\Omega})}(e^{hD}- e^{-hD}),\ \
\hat{\mathcal{L}}(hD)=  \tilde{\alpha} (h\tilde{\Omega})
\frac{e^{hD}-2+e^{-hD}}{h^2}
-\tilde{\beta}(h\tilde{\Omega})\mathrm{i}
 \tilde{\Omega} \frac{e^{hD}-e^{-hD}}{2h},\end{array}\end{equation}
 where $D$ is the differential operator. The Taylor expansions of the  operator $\mathcal{L}(hD)$  are
\begin{equation*}
\begin{aligned}\mathcal{L}(hD)=&   \tilde{\Omega}  \csc(h\tilde{\Omega})
(hD)+\frac{1}{6}\tilde{\Omega} \csc(h\tilde{\Omega})
(hD)^3+\cdots,\\
\mathcal{L}(hD+\mathrm{i}h\tilde{\omega})=&\mathrm{i}\textmd{diag}\Big(
 \tilde{\omega},\frac{ \sin(h\tilde{\omega})}{h},\tilde{\omega}\Big)+\textmd{diag}\Big(
 \tilde{\omega}\cot(h\tilde{\omega}),\frac{ \cos(h\tilde{\omega})}{h},\tilde{\omega}\cot(h\tilde{\omega})\Big)(hD)+\cdots,\\
\mathcal{L}(hD-\mathrm{i}h\tilde{\omega})=&-\mathrm{i}\textmd{diag}\Big(
 \tilde{\omega},\frac{ \sin(h\tilde{\omega})}{h},\tilde{\omega}\Big)+\textmd{diag}\Big(
 \tilde{\omega}\cot(h\tilde{\omega}),\frac{ \cos(h\tilde{\omega})}{h},\tilde{\omega}\cot(h\tilde{\omega})\Big)(hD)+\cdots,\\
\mathcal{L}(hD+\mathrm{i}kh\tilde{\omega})=&
\mathrm{i}\sin(kh\tilde{\omega})\tilde{\Omega}\csc(h\tilde{\Omega})+\cos(kh\tilde{\omega})\tilde{\Omega}\csc(h\tilde{\Omega})+\cdots,
\ \ \textmd{where}\ \abs{k}>1.
\end{aligned}
\end{equation*}
The operator $\hat{\mathcal{L}}(hD)$   can be expressed in its
Taylor expansion  as
\begin{equation*}
\begin{aligned}\hat{\mathcal{L}}(hD)=& - \tilde{\Omega}^2 \csc(h\tilde{\Omega})
(\mathrm{i}hD)-\frac{1}{4}\tilde{\Omega}^2
\csc^2\big(\frac{1}{2}h\tilde{\Omega}\big)
(\mathrm{i}hD)^2+\cdots,\\
\hat{\mathcal{L}}(hD+\mathrm{i}h\tilde{\omega})=&\textmd{diag}\Big(-2
 \tilde{\omega}^2,\frac{2
(\cos(h\tilde{\omega})-1)}{h^2},0\Big) \\
&+\textmd{diag}\Big(  \tilde{\omega}^2
(2\cot(h\tilde{\omega})+\csc(h\tilde{\omega})),
\frac{2\sin(h\tilde{\omega})}{h^2},
  \tilde{\omega}^2\csc(h\tilde{\Omega}) \Big)(\mathrm{i}hD)+\cdots,\\
\hat{\mathcal{L}}(hD-\mathrm{i}h\tilde{\omega})=&\textmd{diag}\Big(
0,\frac{2 (\cos(h\tilde{\omega})-1)}{h^2},-2
 \tilde{\omega}^2\Big)\\
&-\textmd{diag}\Big( \tilde{\omega}^2\csc(h\tilde{\Omega}) ,
\frac{2\sin(h\tilde{\omega})}{h^2}, \tilde{\omega}^2
(2\cot(h\tilde{\omega})+\csc(h\tilde{\omega}))\Big) (\mathrm{i}hD)+\cdots,\\
\hat{\mathcal{L}}(hD+\mathrm{i}kh\tilde{\omega})=&
2\sin\big(\frac{1}{2}hk \tilde{\omega}\big)
\tilde{\Omega}^2\csc\big(\frac{1}{2}h\tilde{\Omega}\big)
\csc(h\tilde{\Omega})\sin\big(\frac{1}{2}h(\tilde{\Omega}- k
\tilde{\omega}I)\big)
\\
&-\sin\big(\frac{1}{2}h(\tilde{\Omega}-2 k
\tilde{\omega})\big)\tilde{\Omega}^2\csc\big(\frac{1}{2}h\tilde{\Omega}\big)\csc(h\tilde{\Omega})(\mathrm{i}hD)+\cdots,
\ \ \textmd{where}\ \abs{k}>1.
\end{aligned}
\end{equation*}



$\bullet$ \textbf{Modulated Fourier expansion.}

We first present  the modulated Fourier expansion of the numerical
result $\tilde{x}_{n}$ and $\tilde{v}_{n}$ for solving the transformed
system \eqref{necharged-sts-first order}.

We will   look for smooth coefficient functions
$\tilde{\zeta}_{k}$ and $\tilde{\eta}_{k}$ such that for $t=nh$,  the functions
\begin{equation}
\begin{aligned} &\tilde{x}_{h}(t)= \sum\limits_{|k|<N} \mathrm{e}^{\mathrm{i}k\tilde{\omega}
t}\tilde{\zeta}_{k}(t)+\tilde{R}_{h,N}(t),\quad
 \ \tilde{v}_{h}(t)= \sum\limits_{|k|<N} \mathrm{e}^{\mathrm{i}k\tilde{\omega}
 t}\tilde{\eta}_{k}(t)+\tilde{S}_{h,N}(t)
\end{aligned}
\label{MFE-1}%
\end{equation} yield a small defect $\tilde{R},\tilde{S}$ when they are inserted into the
numerical scheme  \eqref{2 steps method-new}.

\vskip1mm  $\circ$  Construction of the coefficients
functions.

 Inserting  the first expansion of \eqref{MFE-1} into the
two-step form \eqref{2 steps method-new}, expanding the nonlinear
function into its Taylor series and comparing the coefficients of
$\mathrm{e}^{\mathrm{i}k\tilde{\omega} t}$,  we obtain
\begin{equation}\label{atanew}
\begin{aligned}& \hat{\mathcal{L}}(hD)\tilde{\zeta}_0=\tilde{\gamma}
 (h\tilde{\Omega})\Big(\tilde{F}(\tilde{\zeta}_0)+
\sum\limits_{s(\alpha)=0}\frac{1}{m!}\tilde{F}^{(m)}( \tilde{\zeta}_0)( \tilde{\zeta})_{\alpha}\Big),\\
&\hat{\mathcal{L}}(hD+\mathrm{i}hk\tilde{\omega})\tilde{\zeta}_{k}
=\tilde{\gamma}
 (h\tilde{\Omega})
\sum\limits_{s(\alpha)=k}\frac{1}{m!}\tilde{F}^{(m)}( \tilde{\zeta}_0)( \tilde{\zeta})_{\alpha},\quad \abs{k}>0,\\
\end{aligned} %
\end{equation}
where the sum ranges over $m\geq0$,
$s(\alpha)=\sum\limits_{j=1}^{m}\alpha_j$ with
$\alpha=(\alpha_1,\ldots,\alpha_m)$  and  $0<|\alpha_i|<N$,
 and
$(\tilde{\zeta})_{\alpha}$ is an abbreviation for
$(\tilde{\zeta}_{\alpha_1},\ldots,\tilde{\zeta}_{\alpha_m})$. This
formula gives the  modulation system  for the coefficients
$\tilde{\zeta}_{k}$ of the modulated Fourier expansion. Choosing the dominating terms and considering the Taylor expansion of $\hat{\mathcal{L}}$ given above, the
following ansatz of $\tilde{\zeta}_{k}$  can be obtained:
\begin{equation}\label{ansatz}%
\begin{array}{ll}
\dot{\tilde{\zeta}}_0^{\pm1}=\frac{-h^2\tilde{\omega}A(h\tilde{\omega})}{8
 \mathrm{i}  \sin^2(\frac{1}{2}h\tilde{\omega}) }\big(\mathcal{G}^{\pm10}(\cdot)+\cdots\big),\
&\ddot{\tilde{\zeta}}_0^{0}= \mathcal{G}^{00}(\cdot)+\cdots,\\
\tilde{\zeta}^{-1}_1=\frac{ h^3\tilde{\omega}A(h\tilde{\omega})}{-16
\sin^2(\frac{1}{2}h\tilde{\Omega})\sin(h\tilde{\omega}) }
\big(\mathcal{F}_1^{-10}(\cdot)+\cdots\big),\
&\tilde{\zeta}^{0}_1=\frac{
 h^2}{-4  \sin^2(h\tilde{\omega}/2) }
\big(\mathcal{F}_1^{00}(\cdot)+\cdots\big),
\\
\dot{\tilde{\zeta}}^{1}_1=\frac{h^2\tilde{\omega}A(h\tilde{\omega})}{8\mathrm{i}\sin^2(\frac{1}{2}h\tilde{\omega})}
\big(\mathcal{F}_1^{10}(\cdot)+\cdots\big),\
&\dot{\tilde{\zeta}}^{-1}_{-1}=\frac{h^2\tilde{\omega}A(h\tilde{\omega})}{-8\mathrm{i}\sin^2(\frac{1}{2}h\tilde{\omega})}
\big(\mathcal{F}_{-1}^{-10}(\cdot)+\cdots\big),\\
\tilde{\zeta}^{0}_{-1}=\frac{ h^2}{-4\sin^2(h\tilde{\omega}/2) }
\big(\mathcal{F}_{-1}^{00}(\cdot)+\cdots\big), \
&\tilde{\zeta}^{1}_{-1}=\frac{
h^3\tilde{\omega}A(h\tilde{\omega})}{-16
\sin^2(\frac{1}{2}h\tilde{\Omega})\sin(h\tilde{\Omega}) }
\big(\mathcal{F}_{-1}^{10}(\cdot)+\cdots\big),
\\
\tilde{\zeta}_{k}=\frac{
h^3\tilde{\Omega}A(h\tilde{\Omega})}{16\sin(\frac{1}{2}h\tilde{\Omega})
\sin(\frac{1}{2}h(\tilde{\Omega}- k \tilde{\omega}I))
\sin(\frac{1}{2}hk \tilde{\omega}I)}
\big(\mathcal{F}_{k}^{0}(\cdot)+\cdots\big) &\textmd{for}\
\abs{k}>1,
\end{array} %
\end{equation}
where  the dots  stand  for power series in $\sqrt{h}$ and $A(h\tilde{\omega})= 2\sinc^2(\frac{1}{2}h\tilde{\omega})$. In this
paper we truncate the ansatz after the $\mathcal{O}(h^{N+1})$ terms.
On the basis of the second formula of \eqref{2 steps method-0}, one
has \begin{equation}\label{v-x-rea}
 \begin{array}[c]{ll}\tilde{v}_{n}&= \frac{1}{h(\varphi_1(\mathrm{i}h\tilde{\Omega})+\varphi_1(-\mathrm{i}h\tilde{\Omega}))}(\tilde{x}_{n+1}-\tilde{x}_{n-1})
 -\frac{1}{2}h^2\frac{\varphi_1(\mathrm{i}h\tilde{\Omega})-
 \varphi_1(-\mathrm{i}h\tilde{\Omega}) }{h(\varphi_1(\mathrm{i}h\tilde{\Omega})+
 \varphi_1(-\mathrm{i}h\tilde{\Omega}))}\tilde{F} (  \tilde{x}_{n})\\
 &= \frac{1}{2h\sinc(h\tilde{\Omega})}(\tilde{x}_{n+1}-\tilde{x}_{n-1})
 -\frac{1}{2}\mathrm{i}h\tan(\frac{h}{2}\tilde{\Omega})\tilde{F} (  \tilde{x}_{n}).
\end{array}
\end{equation}
Inserting  \eqref{MFE-1} into  \eqref{v-x-rea}, expanding the
nonlinear function   into its Taylor series and comparing the
coefficients of $\mathrm{e}^{\mathrm{i}k\tilde{\omega} t}$,  one
arrives
\begin{equation}\label{atav}
\begin{aligned}& \tilde{\eta}_0=\mathcal{L}(hD)\tilde{\zeta}_0- \frac{1}{2}\mathrm{i}h\tan(\frac{h}{2}\tilde{\Omega})\Big(\tilde{F}(\tilde{\zeta}_0)+
\sum\limits_{s(\alpha)=0}\frac{1}{m!}\tilde{F}^{(m)}( \tilde{\zeta}_0)( \tilde{\zeta})_{\alpha}\Big),\\
&\tilde{\eta}_{k}=\mathcal{L}(hD+\mathrm{i}hk\tilde{\omega})\tilde{\zeta}_{k}
-\frac{1}{2}\mathrm{i}h\tan(\frac{h}{2}\tilde{\Omega})
\sum\limits_{s(\alpha)=k}\frac{1}{m!}\tilde{F}^{(m)}( \tilde{\zeta}_0)( \tilde{\zeta})_{\alpha},\quad \abs{k}>0.\\
\end{aligned} %
\end{equation}
This formula gives the  modulation system  for the coefficients
$\tilde{\eta}_{k}$ of the modulated Fourier expansion by the Taylor expansion of $\mathcal{L}$ and by choosing the dominating terms.

\vskip1mm  $\circ$
Initial values.

For the first-order and second-order differential equations appeared in \eqref{ansatz}, initial values are needed and we derive them as follows.

According to the conditions
$\tilde{x}_{h}(0)=\tilde{x}_0$ and
 $\tilde{v}_{h}(0)=\tilde{v}_0$,
we have
\begin{equation}\label{Initial values-1}%
\begin{aligned}
&\tilde{x}^0_0=\tilde{\zeta}^0_0(0)+\mathcal{O}(\tilde{\omega}^{-1}),\ \
\tilde{x}^{\pm1}_0=\tilde{\zeta}^{\pm1}_0(0)+\mathcal{O}(\tilde{\omega}^{-1}),\\
&\tilde{v}^0_0=\tilde{\eta}^0_0(0)+\mathcal{O}(\tilde{\omega}^{-1})=\dot{\tilde{\zeta}}^0_0(0)+\mathcal{O}(\tilde{\omega}^{-1}),\\
&\tilde{v}^{1}_0=\tilde{\eta}^{1}_0(0)+\tilde{\eta}^{1}_{1}(0)+\mathcal{O}(\tilde{\omega}^{-1})=
\dot{\tilde{\zeta}}^{1}_0(0)+\mathrm{i}\tilde{\omega}\tilde{\zeta}^{1}_{1}(0)+\mathcal{O}(\tilde{\omega}^{-1}),\\
&\tilde{v}^{-1}_0=\tilde{\eta}^{-1}_0(0)+\tilde{\eta}^{-1}_{-1}(0)+\mathcal{O}(\tilde{\omega}^{-1})=
\dot{\tilde{\zeta}}^{-1}_0(0)-\mathrm{i}\tilde{\omega}\tilde{\zeta}^{-1}_{-1}(0)+\mathcal{O}(\tilde{\omega}^{-1}). \end{aligned} %
\end{equation}
Thus  the initial values $\tilde{\zeta}_0^0(0)=\mathcal{O}(1)$ and
$\dot{\tilde{\zeta}}_0^0(0)=\mathcal{O}(1)$ can be derived by
considering the first and third formulae, respectively. According to the second
equation of \eqref{Initial values-1},   one gets the initial value
$\tilde{\zeta}^{\pm1}_0(0)=\mathcal{O}(1)$. It follows from  the
fourth formula that $\tilde {\zeta}_{1}^{1}(0)=\frac{1}{\mathrm{i}
\tilde{\omega}}\big(\tilde{v}^{1}_0-\dot{\tilde{\zeta}}^{1}_0(0)+\mathcal{O}(h)\big)
=\mathcal{O}(\tilde{\omega}^{-1})$, and likewise one has
$\tilde{\zeta}^{-1}_{-1}(0)=\mathcal{O}(\tilde{\omega}^{-1})$.

\vskip1mm  $\circ$  Bounds of the coefficient functions.

 With the ansatz \eqref{ansatz},  we achieve the bounds
 \begin{equation*}
\begin{array}{ll}\dot{\tilde{\zeta}}_0^{\pm1}=\mathcal{O}(h),\  \  \ \
\ddot{\tilde{\zeta}}_0^{0}=\mathcal{O}(1),\ &
\dot{\tilde{\zeta}}_1^{1}=\mathcal{O}(h),\quad  \ \  \dot{\tilde{\zeta}}_{-1}^{-1}=\mathcal{O}(h),  \\
\tilde{\zeta}_{1}^{-1}=\mathcal{O}(h^{\frac{5}{2}}), \ \
\tilde{\zeta}_{1}^{0}=\mathcal{O}(h^2),\ &
\tilde{\zeta}_{-1}^{0}=\mathcal{O}(h^2),  \ \
\tilde{\zeta}_{-1}^{1}=\mathcal{O}(h^{\frac{5}{2}}).
\end{array}
\end{equation*}
 According to the initial
values stated above,   the bounds \begin{equation*}
 \tilde{\zeta}_0^{\pm1}=\mathcal{O}(1), \quad  \ \
 \tilde{\zeta}_0^{0}=\mathcal{O}(1),\quad \ \
\tilde{\zeta}_1^{1}=\mathcal{O}(h),\quad \ \
\tilde{\zeta}_{-1}^{-1}=\mathcal{O}(h),
\end{equation*}
are obtained. Moreover, we have the following results for coefficient functions
$\tilde{\eta}$
\begin{equation}
\begin{array}{ll}
\tilde{\eta}_0^0=\dot{\tilde{\zeta}}_0^0+\mathcal{O}(h),\
  &\tilde{\eta}^{\pm 1}_0=\frac{h\tilde{\omega}}{\sin(h\tilde{\omega})}\dot{\tilde{\zeta}}^{\pm 1}_0+\mathcal{O}(h),\\
\tilde{\eta}_{\pm1}^{0}=\mathrm{i}\tilde{\omega}\sinc(h\tilde{\omega})\tilde{\zeta}_{0}^{\pm1}+\mathcal{O}(h),\
&\tilde{\eta}_1^{\pm 1}= \mathrm{i}\tilde{\omega}
\tilde{\zeta}_{1}^{\pm1}+\mathcal{O}(h),  \
\tilde{\eta}_{-1}^{\pm1}=-\mathrm{i}\tilde{\omega}
\tilde{\zeta}_{-1}^{\pm1}+\mathcal{O}(h).
\end{array}
\label{rea20-num}%
\end{equation}
A
further result is true
$$\tilde{\zeta}_{k}=\mathcal{O}(h ^{ \abs{k}+1 } ),\ \ \ \
\tilde{\eta}_{k}=\mathcal{O}(h^{\abs{k}})\ \ \ \ \textmd{for}\ \
\abs{k}>1.
$$

 \vskip1mm  $\circ$
 Defect.

Define
\begin{equation*}
\begin{aligned}
 \delta_1(t+h)=&\tilde{x}_{h}(t+h)-\tilde{x}_{h}(t)- h\varphi_1(\mathrm{i}h\tilde{\Omega})\tilde{v}_{h}(t)-\frac{1}{2}h^2\varphi_1(\mathrm{i}h\tilde{\Omega})\tilde{F}(\tilde{x}_{h}(t)),\\
  \delta_2(t+h)= &\tilde{v}_{h}(t+h)-e^{
  \mathrm{i}h\tilde{\Omega}}\tilde{v}_{h}(t)-\frac{1}{2}h\varphi_0(\mathrm{i}h\tilde{\Omega})\tilde{F}(\tilde{x}_{h}(t))-
  \frac{1}{2}h \tilde{F}(\tilde{x}_{h}(t+h))
\end{aligned}
\end{equation*}
for $t=nh$. Considering the two-step formulation, it is clear that
$\delta_1(t+h)+\delta_1(t-h)=\mathcal{O}(h^{4})$.  According to the
choice for the initial values, we obtain
$\delta_1(0)=\mathcal{O}(h^{N+2}).$ Therefore, one has
$\delta_1(t)=\mathcal{O}(h^{N+2})+\mathcal{O}(t h^{N+1}).$
{Using} this result and \eqref{v-x-rea}, we have
$\delta_2=\mathcal{O}(h^{N}).$ Then let
$\tilde{R}_{n}=\tilde{x}_{n}-\tilde{x}_{h}(t)$ and
$\tilde{S}_{n}=\tilde{v}_{n}-\tilde{v}_{h}(t).$ With the scheme of
SM2, the error recursion is obtained as follows:
\begin{equation*}
\begin{aligned}
\left(
  \begin{array}{c}
    \tilde{R}_{n+1} \\
     \tilde{S}_{n+1} \\
  \end{array}
\right)=\left(
                                                           \begin{array}{cc}
                                                             I & h\varphi_1(\mathrm{i}h\tilde{\Omega}) \\
                                                             0 & e^{\mathrm{i}h\tilde{\Omega}} \\
                                                           \end{array}
                                                         \right)
\left(
  \begin{array}{c}
    \tilde{R}_{n} \\
     \tilde{S}_{n} \\
  \end{array}
\right)+\frac{1}{2}h\left(
          \begin{array}{c}
            h\varphi_1 \Gamma_{n}\tilde{R}_{n} \\
             \varphi_0 \Gamma_{n}\tilde{R}_{n}+  \Gamma_{n+1}\tilde{R}_{n+1} \\
          \end{array}
        \right)+\left(
                  \begin{array}{c}
                    \delta_1 \\
                    \delta_2 \\
                  \end{array}
                \right),
\end{aligned}
\end{equation*}
where $\Gamma_{n}:=\int_{0}^1 \tilde{F}_x(\tilde{x}_{n}+\tau
\tilde{R}_{n}) d\tau$. Solving this recursion and the application of
a discrete Gronwall inequality  gives \begin{equation*}
 \tilde{R}_{h,N}(t)=\mathcal{O}(t^2h^{N}),\ \ \ \ \ \ \ \ \ \tilde{S}_{h,N}(t)=\mathcal{O}(t^2h^{N}/\eps).
\end{equation*}

Using the relationships shown in  \eqref{relationships}, the numerical solution of SM2 admits the following  modulated Fourier expansion
\begin{equation*}
\begin{aligned} &x_{n}=  \sum\limits_{|k|<N} \mathrm{e}^{\mathrm{i}k\tilde{\omega} t}\zeta_{k}(t)+\mathcal{O}(t^2h^{N}),\
\ v_{n}= \sum\limits_{|k|<N} \mathrm{e}^{\mathrm{i}k\tilde{\omega}
t}\eta_{k}(t)+\mathcal{O}(t^2h^{N}/\eps),
\end{aligned}
\end{equation*}
where  $\zeta_{k}=P\tilde{\zeta}_{k}$ and $\eta_{k}=P\tilde{\eta}_{k}.$
Moreover, we have $\zeta_{-k}=\overline{\zeta_{k}}$ and
$\eta_{-k}=\overline{\eta_{k}}$.


$\bullet$ \textbf{An almost-invariant.}

Denote
$$\vec{\tilde{\zeta}}=\big(\tilde{\zeta}_{-N+1},\cdots,\tilde{\zeta}_{-1},
\tilde{\zeta}_{0},\tilde{\zeta}_{1},\cdots,\tilde{\zeta}_{N-1}\big).$$
An almost-invariant of the modulated Fourier expansion \eqref{MFE-1} is
given as follows.

According to the above analysis, it is
deduced that
\begin{equation*}
\begin{aligned}
& \tilde{\gamma}^{-1}
 (h\tilde{\Omega})\hat{\mathcal{L}}(hD) \tilde{x}_{h}=
\tilde{F}(\tilde{x}_{h})+\mathcal{O}(h^{N}),
\end{aligned}
\end{equation*}
where we use  the   denotations $\tilde{x}_{h}=\sum\limits_{
|k|<N}\tilde{x}_{h,k}$ with  $
\tilde{x}_{h,k}=\mathrm{e}^{\mathrm{i}k\tilde{\omega}
t}\tilde{\zeta}_{k}.$ Multiplication of this result with $P$ yields
\begin{equation*}
\begin{aligned}
& P\tilde{\gamma}^{-1}
 (h\tilde{\Omega})\hat{\mathcal{L}}(hD)P^\textup{H}
P\tilde{x}_{h}=P\tilde{\gamma}^{-1}
 (h\tilde{\Omega})\hat{\mathcal{L}}(hD)P^\textup{H} x_{h}\\
 &=
P\tilde{F}(\tilde{x}_{h})+\mathcal{O}(h^{N+2})=
F(x_{h})+\mathcal{O}(h^{N}),
\end{aligned}
\end{equation*}
where $x_{h}=\sum\limits_{ |k|<N}x_{h,k}$  with
$x_{h,k}=\mathrm{e}^{\mathrm{i}k\tilde{\omega} t}\zeta_k.$ Rewrite
the equation in terms of $x_{h,k}$ and then one has
\begin{equation*}
\begin{aligned} &P\tilde{\gamma}^{-1}
 (h\tilde{\Omega})\hat{\mathcal{L}}(hD)P^\textup{H} x_{h,k}=-\nabla_{x_{-k}}\mathcal{U}(\vec{x})+\mathcal{O}(h^{N}),
\end{aligned}
\end{equation*}
where $\mathcal{U}(\vec{x})$ is defined as
\begin{equation*}\mathcal{U}(\vec{x})=U(x_{h,0})+
\sum\limits_{s(\alpha)=0}\frac{1}{m!}U^{(m)}(x_{h,0}) (x_h)_{\alpha}
\end{equation*}
with $\vec{x}=\big(x_{h,-N+1},\cdots,x_{h,-1},
x_{h,0},x_{h,1},\cdots,x_{h,N-1}\big).$
 Multiplying this equation with $
\big(\dot{x}_{h,-k}\big)^\intercal$ and summing up yields
\begin{equation*}
\begin{aligned} & \sum\limits_{|k|<N}
 \big(\dot{x}_{h,-k}\big)^\intercal P\tilde{\gamma}^{-1}
 (h\tilde{\Omega})\hat{\mathcal{L}}(hD)P^\textup{H}
x_{h,k}+\frac{d}{dt}\mathcal{U}(\vec{x})=\mathcal{O}(h^{N}).
\end{aligned}
\end{equation*}
Denoting $\vec{\zeta}=\big(\zeta_{-N+1},\cdots,\zeta_{-1},
\zeta_{0},\zeta_{1},\cdots,\zeta_{N-1}\big)$ and switching to the
quantities  $\zeta^k$,  we obtain
\begin{equation}
\begin{aligned}
 \mathcal{O}(h^{N}) =& \sum\limits_{|k|<N}
 \big(\dot{\zeta}_{-k}-\textmd{i} k \tilde{\omega} \zeta_{-k}\big)^\intercal
P\tilde{\gamma}^{-1}
 (h\tilde{\Omega})\hat{\mathcal{L}}(hD+\mathrm{i}hk\tilde{\omega})P^\textup{H} \zeta_k+\frac{d}{dt}\mathcal{U}(\vec{\zeta})\\
= & \sum\limits_{|k|<N}
 \big(\dot{ \overline{\zeta}}_{k}-\textmd{i} k \tilde{\omega}
 \overline{\zeta_{k}}\big)^\intercal P\tilde{\gamma}^{-1}
 (h\tilde{\Omega})
\hat{\mathcal{L}}(hD+\mathrm{i}hk\tilde{\omega})P^\textup{H}
\zeta_k+\frac{d}{dt}\mathcal{U}(\vec{\zeta})\\
= & \sum\limits_{|k|<N}
 \big(\dot{ \overline{\tilde{\zeta}}}_{k}-\textmd{i} k \tilde{\omega}
 \overline{\tilde{\zeta}_{k}}\big)^\intercal P^\textup{H}P\tilde{\gamma}^{-1}
 (h\tilde{\Omega})
\hat{\mathcal{L}}(hD+\mathrm{i}hk\tilde{\omega})P^\textup{H} P
\tilde{\zeta}_k+\frac{d}{dt}\mathcal{U}(\vec{\zeta})\\
= & \sum\limits_{|k|<N}
 \big( \dot{\overline{\tilde{\zeta}}}_{k}-\textmd{i} k \tilde{\omega}
 \overline{\tilde{\zeta}_{k}}\big)^\intercal \tilde{\gamma}^{-1}
 (h\tilde{\Omega})
\hat{\mathcal{L}}(hD+\mathrm{i}hk\tilde{\omega})
\tilde{\zeta}_k+\frac{d}{dt}\mathcal{U}(\vec{\zeta}).
\end{aligned}
\label{duu-new1}%
\end{equation}
  In what follows, we   show that  the right-hand side of
\eqref{duu-new1} is the total derivative of an expression that
depends only on $\tilde{\zeta}_k$ and derivatives
thereof. 
Consider   $k = 0$ and then it follows that
$$\hat{\mathcal{L}}(hD) \tilde{\zeta}_0=\textmd{i}hM_1\dot{\tilde{\zeta}}_0+h^2M_2 \ddot{\tilde{\zeta}}_0+\textmd{i}h^3M_3
\ddot{\tilde{\zeta}}_0+\cdots, \ \textmd{where} \ M_k\in
\mathbb{R}^{3\times 3}\ \textmd{for}\ k=1,2,\ldots.$$ By the   formulae given  on p. 508 of
\cite{hairer2006}, we know that
$\textmd{Re} \big( \dot{\overline{\tilde{\zeta}}}_{0}\big)^\intercal
 \hat{\mathcal{L}}(hD) \tilde{\zeta}_0$ is a total derivative. For
$k\neq0$, in the light of
\begin{equation*}
\hat{\mathcal{L}}(hD+\mathrm{i}hk\tilde{\omega})
\tilde{\zeta}_k=N_0\tilde{\zeta}_k+\textmd{i}hN_1\dot{\tilde{\zeta}}_k+h^2N_2
\ddot{\tilde{\zeta}}_k+\textmd{i}h^3N_3
\ddot{\tilde{\zeta}}_k+\cdots, \  \textmd{where} \ N_k\in
\mathbb{R}^{3\times 3}\ \textmd{for}\ k=0,1,\ldots,\end{equation*} it
is easy to check that $\textmd{Re} \big(
\dot{\overline{\tilde{\zeta}}}_{k}\big)^\intercal
\tilde{\gamma}^{-1}
 (h\tilde{\Omega})\hat{\mathcal{L}}(hD+\mathrm{i}hk\tilde{\omega}) \tilde{\zeta}_k$ and
$\textmd{Re} \big(\textmd{i} k \tilde{\omega}
\overline{\tilde{\zeta}_{k}}\big)^\intercal \tilde{\gamma}^{-1}
 (h\tilde{\Omega})\hat{\mathcal{L}}(hD+\mathrm{i}hk\tilde{\omega}) \tilde{\zeta}_k$ are both total
derivatives. Therefore, there exists a function
$\mathcal{E}$ such
 that
$\frac{d}{dt}\mathcal{E}[\vec{\tilde{\zeta}}](t)=\mathcal{O}(h^{N})$.
{It follows from an integration that}
\begin{equation}\label{Epro}
\mathcal{E}[\vec{\tilde{\zeta}}](t)=\mathcal{E}[\vec{\tilde{\zeta}}](0)+\mathcal{O}(th^{N}).
\end{equation}

{On the basis of}  the previous analysis,  the
construction of $\mathcal{E}$ is derived as follows
\begin{equation*}
\begin{aligned}\mathcal{E}[\vec{\tilde{\zeta}}](t_n)=&
\frac{1}{2}(\dot{\overline{\tilde{\zeta}}}^0)
 ^\intercal \frac{2\sinc(\frac{1}{2}h\tilde{\Omega})}{
\varphi_1(\mathrm{i}h\tilde{\Omega})\varphi_1(-\mathrm{i}h\tilde{\Omega})+\varphi_1(-\mathrm{i}h\tilde{\Omega})\varphi_1(\mathrm{i}h\tilde{\Omega})}\dot{\tilde{\zeta}}^0
\\
&+\frac{1}{2} \frac{\tilde{\omega}}{h} h\tilde{\omega}
\frac{2\sinc^2(\frac{1}{2}h\tilde{\omega})}{
\varphi_1(\mathrm{i}h\tilde{\omega})\varphi_1(-\mathrm{i}h\tilde{\omega})+\varphi_1(-\mathrm{i}h\tilde{\omega})\varphi_1(\mathrm{i}h\tilde{\omega})}\big(\abs{
\tilde{\zeta}_1^1}^2+\abs{
\tilde{\zeta}_{-1}^{-1}}^2\big)+\mathcal{U}(\vec{\zeta})+\mathcal{O}(h^2)\\
=& \frac{1}{2}\abs{ \dot{\tilde{\zeta}}_0^0}^2
+\frac{1}{2}\tilde{\omega}^2
 \big(\abs{ \tilde{\zeta}_1^1}^2+\abs{
\tilde{\zeta}_{-1}^{-1}}^2\big)+U(
P^\textup{H}\tilde{\zeta}^0)+\mathcal{O}(h).
\end{aligned}
\end{equation*}


$\bullet$ \textbf{Long-time  near-conservation.}

Considering the result of
$\mathcal{E}$ and the
 relationship \eqref{rea20-num} between
$\tilde{\zeta}$ and $\tilde{\eta}$, we obtain
\begin{equation}\label{HIW}
\begin{aligned}   \mathcal{E}
[\vec{\tilde{\zeta}}](t_n) &=\frac{1}{2}\abs{ \dot{\tilde{\zeta}}_0^0}^2
+\frac{1}{2}\tilde{\omega}^2
 \big(\abs{ \tilde{\zeta}_1^1}^2+\abs{
\tilde{\zeta}_{-1}^{-1}}^2\big)+U( P^\textup{H}\tilde{\zeta}^0)+\mathcal{O}(h)\\
&=\frac{1}{2}\abs{ \tilde{\eta}_0^0}^2 +\frac{1}{2}
 \big(\abs{ \tilde{\eta}_1^1}^2+\abs{
\tilde{\eta}_{-1}^{-1}}^2\big)+U( P^\textup{H}\tilde{\zeta}^0)+\mathcal{O}(h).
\end{aligned}
\end{equation}
 We are now in a position to show  the long-time
conservations of SM2.

In terms of    the bounds of the coefficient functions, one arrives at
 \begin{equation}\label{HIE}\begin{aligned} &E(x_{n},v_{n})=\tilde{E}(\tilde{x}_{n},\tilde{v}_{n})=\frac{1}{2}\big(\abs{ \tilde{\eta}_0^0}^2  +  \abs{ \tilde{\eta}_1^1}^2+\abs{
\tilde{\eta}_{-1}^{-1}}^2\big)+U(P^\textup{H}\tilde{\zeta}^0)
+\mathcal{O}(h).
\end{aligned}\end{equation}
A comparison between   \eqref{HIW} and \eqref{HIE} gives
$
\mathcal{E} [\vec{\tilde{\zeta}}](t_n)=E(x_{n},v_{n}) +\mathcal{O}(h).
$
 Based on \eqref{Epro} and this result, the
statement \eqref{long EP} is  easily obtained by considering $nh^N\leq1$ and
\begin{equation*}
\begin{array}{ll}
E(x_{n},v_{n})&=\mathcal{E} [\vec{\tilde{\zeta}}](t_n)+\mathcal{O}(h)=\mathcal{E} [\vec{\tilde{\zeta}}](t_{n-1})+\mathcal{O}(h^{N+1})+\mathcal{O}(h)\\
&=\mathcal{E} [\vec{\tilde{\zeta}}](t_{n-2})+2\mathcal{O}(h^{N+1})+\mathcal{O}(h)=\cdots\\
&=\mathcal{E} [\vec{\tilde{\zeta}}](t_{0})+n\mathcal{O}(h^{N+1})+\mathcal{O}(h)=E(x_{0},v_{0})+\mathcal{O}(h).
\end{array}
\end{equation*}

%

$\bullet$  \textbf{Extension of the proof to other $d$.}

According to  the scheme \eqref{lambda} of $\Lambda$, the operators $\mathcal{L}(hD)$  and $\hat{\mathcal{L}}(hD)$ as well as   their properties
can be changed accordingly. The modulated Fourier expansions are modified as
\begin{equation*}
\begin{aligned} &\tilde{x}_{h}(t)= \sum\limits_{k\in\mathcal{N}^*} \mathrm{e}^{\mathrm{i}(k \cdot \tilde{\Omega})
t}\tilde{\zeta}_{k}(t)+\tilde{R}_{h,N}(t),\quad
 \ \tilde{v}_{h}(t)= \sum\limits_{k\in\mathcal{N}^*} \mathrm{e}^{\mathrm{i}(k \cdot \tilde{\Omega})
 t}\tilde{\eta}_{k}(t)+\tilde{S}_{h,N}(t),
\end{aligned}
\end{equation*}
 where $k=(k_1,\ldots,k_l),\ \tilde{\Omega}=(\tilde{\omega}_1,\ldots,\tilde{\omega}_l),\
k\cdot
\tilde{\Omega}=k_1\tilde{\omega}_1+\cdots+k_l\tilde{\omega}_l. $ The
set $\mathcal{N}^*$ is defined as follows. For  the resonance module
$\mathcal{M}= \{k\in \mathbb{Z}^{l}:\ k\cdot \tilde{\Omega}=0\}$, we
let $\mathcal{K}$ be a set of representatives of the equivalence
classes in $\mathbb{Z}^l\backslash \mathcal{M}$ which are chosen
such that for each $k\in\mathcal{K}$ the sum $|k|$ is minimal in the
equivalence class $[k] = k +\mathcal{M},$ and that with
$k\in\mathcal{K}$, also $-k\in\mathcal{K}.$ We denote, for the
positive integer $N$, $\mathcal{N}=\{k\in\mathcal{K}:\ |k|\leq N\},\
\mathcal{N}^*=\mathcal{N}\backslash \{(0,\ldots,0)\}. $ Then the
almost-invariant  can be modified accordingly and the
{long-time} near conservation can be proved.

\section{Proof of  convergence (Theorem \ref{thm: 4})}\label{sec: UA peoperty}
In this section, we discuss the convergence of the algorithms. The proof will be firstly given for M1-M2 and EM1 by using the averaging technique  and then presented for SM1-SM3 by using modulated Fourier expansion.
\subsection{Proof for M1-M2 and EM1}
The proof will be given for EM1 and it can be adapted to M1-M2 easily.

$\bullet$ \textbf{Re-scaled system and methods.}

In order to establish the uniform error bounds,  the strategy developed in \cite{Chartier,zhao2020} will be used in the proof.
This means that the time re-scaling
$\tau:=t/\eps$ is considered and
$ \dot{q}(\tau)=\eps \dot{x}(t),\ \dot{w}(\tau)=\eps \dot{v}(t),
$
where the notations $q(\tau):=x(t),\   w(\tau):=v(t)$ are used.
Then   the convergent analysis will be given for the following   long-time problem
\begin{equation}\label{model}
\begin{split}
\dot{q}(\tau)=\eps w(\tau),\quad\dot{w}(\tau)=\tilde{B} w(\tau) +\eps F(q(\tau)),\  \tau\in[0, \frac{T}{\eps}],\\
q(0)=q_0:=x_0,\  w(0)=w_0:=v_0,
\end{split}
\end{equation}
where $\dot{q}$ (resp. $\dot{w}$) is referred to the
derivative of $q$ (resp. $w$) with respect to $\tau$.
The solution of this system satisfies $\|q\|_{L^\infty(0,T/\eps)}+\|w\|_{L^\infty(0,T/\eps)}\lesssim 1$
and for solving (\ref{model}), the method EM1 becomes
\begin{equation}\label{scheme}
\begin{split}
q_{n+1}=&q_{n}+
\eps \Delta \tau\varphi_1(\Delta \tau\tilde{B})w_{n}+\frac{\Delta\tau^2\eps^2}{2} \varphi_2(\Delta\tau\tilde{B}) \int_{0}^1
F\big(\rho q_{n}+(1-\rho)q_{n+1}\big) d\rho,\\
w_{n+1}=&\e^{\Delta\tau\tilde{B}}w_{n}+\Delta\tau\eps\varphi_1(\Delta\tau\tilde{B}) \int_{0}^1
F\big(\rho q_{n}+(1-\rho)q_{n+1}\big) d\rho,\quad 0\leq n<\frac{T}{\eps\Delta\tau}.
\end{split}
\end{equation}
where $\Delta \tau$ is the time step $\Delta \tau=\tau_{n+1}-\tau_n$ and $q_{n}\approx q(\tau_n),\ w_n\approx w(\tau_n)$ is
the numerical solution.

$\bullet$ \textbf{Local truncation errors.}

Based on  (\ref{scheme}),   the local truncation errors $\xi_n^q$ and $\xi_{n}^w$ for $0\leq n<\frac{T}{\eps\Delta\tau}$  are defined as
\begin{equation}\label{local error}
\begin{split}
q(\tau_{n+1})=&q(\tau_n)+\Delta\tau\eps \varphi_1(\Delta\tau \tilde{B})w(\tau_n)\\
&+ \frac{\Delta\tau^2\eps^2}{2}\varphi_2(\Delta\tau \tilde{B}) \int_0^1F\left(\rho q(\tau_n)+(1-\rho)q(\tau_{n+1})\right)d\rho+\xi_{n}^q,\\
w(\tau_{n+1})=&\e^{\Delta\tau\tilde{B}}w(\tau_n)+\Delta\tau\eps \varphi_1(\Delta\tau \tilde{B}) \int_0^1F\left(\rho q(\tau_n)+(1-\rho)q(\tau_{n+1})\right)d\rho+\xi_{n}^w.
  \end{split}
\end{equation}
 By this result   and the variation-of-constants formula of \eqref{model}, we compute
\[
\begin{aligned}\xi_{n}^w
=&\varepsilon \Delta\tau \int_{0}^1
 e^{(1-\sigma)\Delta\tau \tilde{B}} F(q(\tau_n+\Delta\tau\sigma))  d\sigma\\&- \varepsilon \Delta\tau \varphi_{1}(
\Delta\tau \tilde{B}) \int_{0}^1
F\big(q(\tau_{n})+\sigma(q(\tau_{n+1})-q(\tau_{n})) \big) d\sigma\\
=&\varepsilon  \sum\limits_{j=0}^{1}\Delta\tau^{j+1}\varphi_{j+1}(h
\tilde{B}) \hat F^{(j)}(\tau_n)
-\varepsilon \Delta\tau\varphi_1(\Delta\tau \tilde{B})F(q(\tau_{n}))\\
&-\varepsilon^2 \Delta\tau^2\varphi_1(\Delta\tau \tilde{B})\int_{0}^{1}[\sigma \dfrac{\partial{F}}{\partial q}(q(\tau_{n})) w(\tau_{n})]d\sigma+\mathcal{O}(\varepsilon^2 \Delta\tau^3)\\
=&\varepsilon \Delta\tau^{2}  \varphi_{2}(\Delta\tau \tilde{B}) \hat F^{(1)}(\tau_n)
-\frac{1}{2}\varepsilon^2 \Delta\tau^2\varphi_1(\Delta\tau \tilde{B})
\dfrac{\partial{F}}{\partial q}(q(\tau_{n})) w(\tau_{n})+\mathcal{O}(\varepsilon^2 \Delta\tau^3),
\end{aligned}
\]
where $\hat{F}%
(\xi)=F(q(\xi))$ and $\hat{F}^{(j)}$ denotes the $j$th derivative of
$\hat{F}$ with respect to $\tau$. By this definition, it follows that
$$\hat F^{(1)}(\tau_n)=\dfrac{\partial{F}}{\partial q}(q(\tau_{n}))\frac{\textmd{d} q}{\textmd{d} \tau}(\tau_{n})
=\dfrac{\partial{F}}{\partial q}(q(\tau_{n}))\varepsilon w(\tau_{n}).$$ Consequently, the
local error becomes
\begin{equation}\label{zeta est}
\begin{aligned}&\xi_{n}^w=&\varepsilon^2 \Delta\tau^{2}  \big(\varphi_{2}(\Delta\tau
\tilde{B}) -\frac{1}{2}   \varphi_1(\Delta\tau
\tilde{B})\big) \dfrac{\partial{F}}{\partial q}(q(\tau_{n}))
w(\tau_{n})+\mathcal{O}(\varepsilon^2 \Delta\tau^3)=\mathcal{O}(\varepsilon^2 \Delta\tau^3),
\end{aligned}
\end{equation}
where the result $\varphi_{2}(\Delta\tau
\tilde{B}) -\frac{1}{2}   \varphi_1(\Delta\tau
\tilde{B})=\mathcal{O}(\Delta\tau)$ is used here.
 Similarly, we obtain \begin{equation}\label{xi est}\xi_{n}^q=\mathcal{O}(\varepsilon^3 \Delta\tau^3).\end{equation}

\begin{rem}
It is noted that for M1, the local truncation errors are
\begin{equation}\label{xi estM1}\xi_{n}^w=\mathcal{O}(\varepsilon^2 \Delta\tau^2),\ \ \xi_{n}^q=\mathcal{O}(\varepsilon^2 \Delta\tau^2).\end{equation}
\end{rem}

$\bullet$ \textbf{Error bound.}

In this part, we will first prove {the  boundedness of
EM1:} there exists a generic constant $\Delta\tau_0>0$ independent of $\eps$
and $n$, such that for $0<\Delta\tau\leq \Delta\tau_0$, {the following
inequalities are true:}
 \begin{equation}\label{bounded}
|q_{n}|\leq \|q\|_{L^\infty(0,T/\eps)}+1,\quad |w_{n}|\leq \|w\|_{L^\infty(0,T/\eps)}+1,\quad
0\leq n\leq \frac{T}{\eps\Delta\tau}.
\end{equation}
 For $n=0$, (\ref{bounded}) is obviously true. Then we assume (\ref{bounded}) is true up to some $0\leq m<\frac{T}{\eps\Delta\tau}$, and we shall show that  (\ref{bounded}) holds for $m+1$.

Define the error of the scheme
$$e_n^{q}:=q(\tau_{n})-q_{n},\quad e_n^{w}:=w(\tau_{n})-w_{n},\quad 0\leq n< \frac{T}{\eps\Delta\tau}.$$
For $n\leq m$, subtracting (\ref{local error}) from the scheme (\ref{scheme}) leads to
\begin{equation}\label{lm err}
 e_{n+1}^q=e_{n}^q+\Delta\tau\eps \varphi_{1}(\Delta\tau \tilde{B}) e_{n}^w+\eta_{n}^q+\xi_{n}^q,\
 e_{n+1}^w= \e^{\Delta\tau\tilde{B}}e_{n}^w+\eta_{n}^w+\xi_{n}^w,\ 0\leq n\leq m,
 \end{equation}
 where we use the following notations
 \begin{align*}
   \eta_{n}^q=&\frac{\Delta\tau^2\eps^2}{2}\varphi_{2}(\Delta\tau \tilde{B})\int_0^1\left[F\left(\rho q(\tau_n)+(1-\rho)q(\tau_{n+1})
   \right)-F\left(\rho q_n+(1-\rho)q_{n+1}
   \right)\right]d\rho,\\
   \eta_{n}^w=&\Delta\tau\eps \varphi_{1}(\Delta\tau \tilde{B})\int_0^1\left[F\left(\rho q(\tau_n)+(1-\rho)q(\tau_{n+1})
   \right)-F\left(\rho q_n+(1-\rho)q_{n+1}
   \right)\right]d\rho.
 \end{align*}
From the induction assumption of the boundedness, it follows that
 \begin{align}\label{eta}
|\eta_{n}^q|\lesssim \Delta\tau^2\eps^2\left(|e_{n}^q|+|e_{n+1}^q|\right),\quad
|\eta_{n}^w|\lesssim \Delta\tau\eps\left(|e_{n}^q|+|e_{n+1}^q|\right),\quad 0\leq n<m.
 \end{align}

Taking the absolute value ({Euclidean} norm) on both
sides of (\ref{lm err}) and using (\ref{eta}), {we have}
 \begin{align*}
   |e_{n+1}^q|+|e_{n+1}^w|-|e_{n}^q|-|e_{n}^w|\lesssim &\Delta\tau\eps \left(|e_{n}^w|+|e_{n}^q|+|e_{n+1}^q|\right)+|\xi_{n}^q|
   +|\xi_{n}^w|,\quad 0\leq n\leq m.
 \end{align*}
Summing them up for $0\leq n\leq m$ gives
 $$|e_{m+1}^q|+|e_{m+1}^w|\lesssim \Delta\tau\eps \sum_{n=0}^{m}\left(|e_{n}^w|+|e_{n}^q|+|e_{n+1}^q|\right)+
 \sum_{n=0}^{m}\left(|\xi_{n}^q|
   +|\xi_{n}^w|\right).
 $$
In the light of the truncation errors in (\ref{zeta est}) and the fact that $m\Delta\tau\eps\lesssim 1$, one has
$$ |e_{m+1}^q|+|e_{m+1}^w|\lesssim \Delta\tau\eps \sum_{n=0}^{m}\left(|e_{n}^w|+|e_{n}^q|+|e_{n+1}^q|\right)+\eps\Delta\tau^2,$$
{and then  by Gronwall's inequality arrives at}
 \begin{equation}\label{condi con}|e_{m+1}^q|+|e_{m+1}^w|\lesssim \eps\Delta\tau^2,\quad 0\leq m<\frac{T}{\eps\Delta\tau}. \end{equation}
Meanwhile, concerning
$$|q_{m+1}|\leq |q(\tau_{m+1})|+|e_{m+1}^q|,\quad |w_{m+1}|\leq |w(\tau_{m+1})|+|e_{m+1}^w|,$$
 there exists a generic constant $\Delta\tau_0>0$ independent of $\eps$ and $m$, such that for $0<\Delta\tau\leq \Delta\tau_0$, (\ref{bounded})
 holds for $m+1$. {This} completes the induction.

Now we rewrite \eqref{lm err} as \begin{equation*}
 e_{n+1}^q=e_{n}^q+\Delta\tau \varphi_{1}(\Delta\tau \tilde{B}) (\eps e_{n}^w)+\eta_{n}^q+\xi_{n}^q,\
 (\eps e_{n+1}^w)= \e^{\Delta\tau\tilde{B}}(\eps e_{n}^w)+\eps\eta_{n}^w+\eps\xi_{n}^w.
 \end{equation*}
Following the same way as stated above, it is arrived that
 \begin{equation}\label{condi con new}|e_{m+1}^q|+|\eps e_{m+1}^w|\lesssim \eps^2 \Delta\tau^2,\quad 0\leq m<\frac{T}{\eps\Delta\tau}. \end{equation}

\begin{rem}
We note that for M1, the global  error given in \eqref{condi con} becomes
 \begin{equation*}|e_{m+1}^q|+|e_{m+1}^w|\lesssim \eps\Delta\tau,\quad 0\leq m<\frac{T}{\eps\Delta\tau}, \end{equation*}
 which proves the result \eqref{con M1} of M1.
\end{rem}

$\bullet$  \textbf{Proof of the results for the methods applied to \eqref{charged-particle sts-cons}.}

By considering the grids in the $t$ variable and $\tau$ variable, it is obtained that $h=\eps \Delta \tau$.
This shows that for the original system \eqref{charged-particle sts-cons} and the re-scaled system \eqref{model},
$x(t_n)=q(\tau_n )$ and $v(t_n)=w(\tau_n )$. Moreover, by comparing \eqref{EAVF} with \eqref{scheme}, we know that the numerical solution $x_{n}, v_{n}$ of \eqref{EAVF} is identical to $q_{n}, w_{n}$ of \eqref{scheme}. Therefore, the result
\eqref{condi con new} yields the uniform error bound in $x$ given in \eqref{con M2} and also shows the non-uniform error in $v$  of \eqref{con M2}.

\subsection{Proof for SM1-SM3}
For SM1-SM3, the above proof cannot be applied since  their local truncation errors  will lose a factor of $\eps$ in (\ref{zeta est}) and (\ref{xi est}).   This motivates us to consider modulated Fourier expansions (see, e.g. \cite{Hairer00,Hairer2018,lubich19,hairer2006,Wang2020}) for analysis in this part. The proof will be briefly shown for SM2 and it can be   modified   for the other two methods easily.

$\bullet$  \textbf{Decomposition of the numerical solution.}
Now we turn back to the SM2 given in \eqref{2 steps method-new}  and consider its modulated Fourier expansion \eqref{MFE-1}.
In order to derive the convergence, we need to  explicitly present the results of $\tilde{\zeta}_{k}$ and $\tilde{\eta}_{k}$ with $\abs{k}\leq1$.
In the light of  \eqref{atanew} and the properties of $\hat{\mathcal{L}}(hD)$, we obtain
\begin{equation}\label{ansatz-num}%
\begin{array}{ll}
\dot{\tilde{\zeta}}_0^{\pm1}=\frac{-h^2\tilde{\omega}A(h\tilde{\omega})}{8
 \mathrm{i}  \sin^2(\frac{1}{2}h\tilde{\omega}) }\Big(\tilde{F}(\tilde{\zeta}_0)+\tilde{F}''(\tilde{\zeta}_0)(\tilde{\zeta}_{1},\tilde{\zeta}_{-1})\Big)_{\pm1},\
&\ddot{\tilde{\zeta}}^0_{0}= \Big(\tilde{F}(\tilde{\zeta}_0)+\tilde{F}''(\tilde{\zeta}_0)(\tilde{\zeta}_{1},\tilde{\zeta}_{-1})\Big)_0,\\
\tilde{\zeta}^{-1}_1=\frac{ h^3\tilde{\omega}A(h\tilde{\omega})}{-16
\sin^2(\frac{1}{2}h\tilde{\Omega})\sin(h\tilde{\omega}) }
\big(\tilde{F}'(\tilde{\zeta}_0)\tilde{\zeta}_{1}\big)_{-1},\
&\tilde{\zeta}^{0}_1=\frac{
 h^2}{-4  \sin^2(h\tilde{\omega}/2) }
\big(\tilde{F}'(\tilde{\zeta}_0)\tilde{\zeta}^{1}\big)_0,
\\
\dot{\tilde{\zeta}}_{1}^1=\frac{h^2\tilde{\omega}A(h\tilde{\omega})}{8\mathrm{i}\sin^2(\frac{1}{2}h\tilde{\omega})}
\big(\tilde{F}'(\tilde{\zeta}_0)\tilde{\zeta}_{1}\big)_1,\
&\dot{\tilde{\zeta}}_{-1}^{-1}=\frac{h^2\tilde{\omega}A(h\tilde{\omega})}{-8\mathrm{i}\sin^2(\frac{1}{2}h\tilde{\omega})}
\big(\tilde{F}'(\tilde{\zeta}_0)\tilde{\zeta}_{-1}\big)_{-1},\\
\tilde{\zeta}^{0}_{-1}=\frac{ h^2}{-4\sin^2(h\tilde{\omega}/2) }
\big(\tilde{F}'(\tilde{\zeta}_0)\tilde{\zeta}_{-1}\big)_0, \
&\tilde{\zeta}^{1}_{-1}=\frac{
h^3\tilde{\omega}A(h\tilde{\omega})}{-16
\sin^2(\frac{1}{2}h\tilde{\Omega})\sin(h\tilde{\Omega}) }
\big(\tilde{F}'(\tilde{\zeta}_0)\tilde{\zeta}_{-1}\big)_1.
\end{array} %
\end{equation}
 Then the following results \begin{equation}
\begin{array}{ll}
\tilde{\eta}_0^0=\dot{\tilde{\zeta}}_0^0+\mathcal{O}(h),\
  &\tilde{\eta}^{\pm 1}_0=\frac{h\tilde{\omega}}{\sin(h\tilde{\omega})}\dot{\tilde{\zeta}}^{\pm 1}_0+\mathcal{O}(h),\\
\tilde{\eta}_{\pm1}^{0}=\mathrm{i}\tilde{\omega}\sinc(h\tilde{\omega})\tilde{\zeta}^{0}_{\pm1}+\mathcal{O}(h),\
&\tilde{\eta}_1^{\pm 1}= \mathrm{i}\tilde{\omega}
\tilde{\zeta}_{1}^{\pm1}+\mathcal{O}\Big(h
\abs{\mathrm{i}\tan(\frac{h}{2}\tilde{\Omega})}  \Big),  \\
\tilde{\eta}_{-1}^{\pm1}=-\mathrm{i}\tilde{\omega}
\tilde{\zeta}_{-1}^{\pm1}+\mathcal{O}\Big(h
\abs{\mathrm{i}\tan(\frac{h}{2}\tilde{\Omega})}  \Big)
\end{array}
\label{rea20-num}%
\end{equation} are easily arrived by  considering \eqref{atav} as well as the property
of $\mathcal{L}(hD)$.

$\bullet$  \textbf{Decomposition of the exact solution.}
  Following the result given in  \cite{lubich19},
the exact solution of   \eqref{necharged-sts-first order} admits the
following expansion
\begin{equation}
\begin{aligned} &\tilde{x}(t)=  \sum\limits_{|k|\leq 1} \mathrm{e}^{\mathrm{i}k\tilde{\omega}
t}\tilde{\mu}_k(t)+\tilde{d}_{\tilde{x}}(t),\ \ \  \ \tilde{v}(t)=
\sum\limits_{|k|\leq 1} \mathrm{e}^{\mathrm{i}k\tilde{\omega}
t}\tilde{\nu}_k(t)+\tilde{d}_{\tilde{v}}(t),
\end{aligned}
\label{MFE-exact}%
\end{equation}
where the defects are bounded by $
\tilde{d}_{\tilde{x}}(t)=\mathcal{O}(\tilde{\omega}^{-2}),\
\tilde{d}_{\tilde{v}}(t)=\mathcal{O}(\tilde{\omega}^{-2}/\eps).$
The  functions   $\tilde{\mu}^k$  are given by
\begin{equation}\label{ansatz-exact}%
\begin{array}{ll}
\dot{\tilde{\mu}}_0^{\pm1}=\frac{1}{
 \mp \mathrm{i}  \tilde{\omega}}\big( \tilde{F}_{\pm1}(\tilde{\mu}^0)+\tilde{F}_{\pm1}''(\tilde{\mu}^0)(\tilde{\mu}^{1},\tilde{\mu}^{-1})\big),\
&\ddot{\tilde{\mu}}^0_{0}=\tilde{F}_{0}(\tilde{\mu}^0)+\tilde{F}_{0}''(\tilde{\mu}^0)(\tilde{\mu}^{1},\tilde{\mu}^{-1}),\\
\tilde{\mu}^{-1}_1=\frac{ 1}{-2
\tilde{\omega}^2}\tilde{F}_{-1}'(\tilde{\mu}^0)\tilde{\mu}^{1},\
&\tilde{\mu}^{0}_1=\frac{ 1}{-
\tilde{\omega}^2}\tilde{F}_{0}'(\tilde{\mu}^0)\tilde{\mu}^{1},\\
\dot{\tilde{\mu}}_{1}^1=\frac{1}{\mathrm{i} \tilde{\omega}}
\tilde{F}_{1}'(\tilde{\mu}^0)\tilde{\mu}^{1},\
&\dot{\tilde{\mu}}_{-1}^{-1}=\frac{1}{-\mathrm{i} \tilde{\omega}}
\tilde{F}_{-1}'(\tilde{\mu}^0)\tilde{\mu}^{-1},\\
\tilde{\mu}^{0}_{-1}=\frac{ 1}{-
\tilde{\omega}^2}\tilde{F}_{0}'(\tilde{\mu}^0)\tilde{\mu}^{-1}, \
&\tilde{\mu}^{1}_{-1}=\frac{ 1}{-2
\tilde{\omega}^2}\tilde{F}_{1}'(\tilde{\mu}^0)\tilde{\mu}^{-1},
\end{array} %
\end{equation}
and the  functions   $\tilde{\nu}^k$ are
\begin{equation}
\begin{aligned} &\tilde{\nu}_0=\dot{\tilde{\mu}}_0,   \quad \
\tilde{\nu}_{\pm 1}=\pm \mathrm{i}\tilde{\omega}\tilde{\mu}_{\pm
1}+\dot{\tilde{\mu}}_{\pm 1}=\pm
\mathrm{i}\tilde{\omega}\tilde{\mu}_{\pm
1}+\mathcal{O}(\tilde{\omega}^{-1}).
\end{aligned}
\label{rea20-exact}%
\end{equation}
   The initial values are the same as those of the numerical solutions.

$\bullet$  \textbf{Proof of the convergence.}

Looking closely to the equations \eqref{ansatz-exact} and
\eqref{ansatz-num}, which determine the modulated Fourier expansion
coefficients, it is obtained that
$\tilde{x}^{\ast}(t)-\tilde{x}_h(t)=\mathcal{O}(h^2).$ Similarly,
with  \eqref{rea20-exact} and \eqref{rea20-num}, one has
$\tilde{v}^{\ast}(t)-\tilde{v}_h(t)=\mathcal{O}(h^2).$ According to
the above results and the defects of modulated Fourier expansions, we have the following diagram:
\begin{equation*}
\xymatrix{
\textmd{Exact solution }  (\tilde{x}(nh),\tilde{v}(nh)) \ar[d]                  &\textmd{Numerical solution }  (\tilde{x}_n,\tilde{v}_n ) \ar[d]   \\
\textmd{Modulated Fourier expansion }
(\tilde{x}^{\ast},\tilde{v}^{\ast})\quad
\ar[u]_{(\mathcal{O}(h^2),\mathcal{O}(h^2/\eps))}\ar[r]^{(\mathcal{O}(h^2),\mathcal{O}(h^2))}
&\ \quad\textmd{Modulated Fourier expansion }
(\tilde{x}_h,\tilde{v}_h
)\ar[l]\ar[u]_{(\mathcal{O}(h^2),\mathcal{O}(h^2/\eps))} }\end{equation*}
The error bounds  $$\tilde{x}(nh)-\tilde{x}_{n}=  \mathcal{O}(h^2),\quad \
\tilde{v}(nh)-\tilde{v}_n=  \mathcal{O}(h^2/\eps)$$ are immediately obtained on the basis of
 this diagram. This obviously yields $$x(nh)-x_{n}=  \mathcal{O}(h^2),\quad \
v(nh)-v_n=  \mathcal{O}(h^2/\eps)$$
and the proof is complete.

\section{Conclusions} \label{sec:conclusions}
Structure-preserving algorithms constitute an interesting and
important class of numerical methods. Furthermore,  algorithms with uniformly
errors of highly oscillatory systems have received a
great deal of attention. In this paper, we have formulated and
analysed  some   structure-preserving algorithms with uniform error bound
for solving nonlinear highly oscillatory Hamiltonian systems.
 Two kinds of algorithms with uniform error bound were given to preserve the symplecticity and energy, respectively.  All the theoretical results were supported by a numerical experiment and were proved in detail.

 Last but not least, it is noted that all the algorithms and analysis  are also suitable to the non-highly oscillatory system \eqref{charged-particle sts-cons}  with $\eps=1$. Meanwhile, there are some issues brought by this paper which can be researched further. For the system \eqref{charged-particle sts-cons} with a matrix $\tilde{B}(x)$ depending on $x$, how to modify the methods and extend the analysis to get higher order uniformly accurate structure-preserving algorithms? This point will be considered in future.  Another issue for future exploration is the extension and application of the methods presented in this paper to the  Vlasov equations under strong magnetic field \cite{VP1,Zhao}.

\section*{Acknowledgements}

This work was supported by the Natural Science Foundation of China (NSFC) under grant 11871393; International Science and Technology Cooperation Program of Shaanxi Key Research \& Development Plan under grant 2019KWZ-08. The first author is grateful to Christian Lubich for his valuable
comments on Theorem \ref{thm: 3} as well as its proof, which was done in part at UNIVERSITAT T\"{U}BINGEN when the first author worked there as a postdoctoral researcher (2017-2019, supported by the Alexander von Humboldt Foundation).


\end{document}